\documentclass[10pt]{amsart}
\usepackage{amsfonts, amsmath, amssymb, amsthm, amscd, enumerate}

\usepackage{euscript,amstext}
\usepackage{mathrsfs}
\usepackage{indentfirst}

\usepackage[utf8]{inputenc}
\allowdisplaybreaks

\theoremstyle{plain}
\newtheorem{thm}{Theorem}[section]
\newtheorem{theorem}[thm]{Theorem}

\newtheorem{lemma}[thm]{Lemma}
\newtheorem{proposition}[thm]{Proposition}

\theoremstyle{remark}
\newtheorem{rem}[thm]{Remark}

\theoremstyle{definition}

\newtheorem{definition}[thm]{Definition}

\newcommand{\rB}{B}

\newcommand{\rM}{M}

\newcommand{\rZ}{Z}

\begin{document}

\title
{Families torsion and Morse functions for covering spaces}

\author{Guangxiang Su}
\thanks{Chern Institute of Mathematics and LPMC, Nankai University, Tianjin, 300071, China. \\
e-mail: guangxiangsu@nankai.edu.cn. Supported by NSFC 11571183}

\begin{abstract}
In this paper we prove the Cheeger-M\"{u}ller theorem for $L^2$-analytic torsion form under the assumption that there exists a fiberwise Morse function and the Novikov-Shubin invariant is positive.
\end{abstract}
\maketitle

\renewcommand{\theequation}{\thesection.\arabic{equation}}
\setcounter{equation}{0}

\section{Introduction}
\setcounter{equation}{0}
Let $F$ be a unitary flat vector bundle on a closed Riemannian
manifold $M$. In \cite{RS}, Ray and Singer defined an analytic
torsion associated to $(M,F)$ and proved that it does not depend on
the Riemannian metric on $M$. Moreover, they conjectured that this
analytic torsion coincides with the classical Reidemeister torsion
defined using a triangulation on $M$ (cf. \cite{Mi}). This
conjecture was later proved in the celebrated papers of Cheeger
\cite{C} and M\"{u}ller \cite{Mu1}. Later, M\"{u}ller in \cite{Mu2} generalized
this result to the case where $F$ is a unimodular flat
vector bundle on $M$. In \cite{BZ1}, inspired by the considerations
of Quillen \cite{Q}, at almost the same time of \cite{Mu2}, Bismut and Zhang reformulated the above
Cheeger-M\"{u}ller theorem as an equality between the Reidemeister
and Ray-Singer metrics defined on the determinant of cohomology, and
proved an extension of it to the case of general flat vector bundle
over $M$. In \cite{BZ1}, they make use of the Witten deformation \cite{W} of the de Rham complex.

Let $\rZ \to \rM \xrightarrow{\pi} \rB $ be a fiber bundle with connected closed fibers $Z_{x}=\pi^{-1}(x)$ and $F$ be a flat complex vector bundle on $M$ with a flat connection $\nabla^{F}$ and a Hermitian metric $h^{F}$. Let $T^{H}M$ be a horizontal distribution for the fiber bundle and $g^{TZ}$ be a vertical Riemannian metric. Then in  \cite{BL} Bismut and Lott introduced the torsion form $\mathcal{T}(T^{H}M, g^{TZ},h^{F})\in \Omega(B)$ (cf. \cite[(3.118)]{BL}). Assuming that there exists a fiberwise $G$-invariant Morse function with a fibrewise Thom-Smale gradient field, Bismut and Goette in \cite{BG} proved the Cheeger-M\"{u}ller type theorem for the equivariant Bismut-Lott torsion form.

The $L^2$-analytic torsion was first introduced by Carey, Lott and Mathai (\cite{CM,M,L}), under the assumptions that the $L^2$-Betti numbers vanish and that certain technical ``determinant class condition" is satisfied. The later condition is satisfied if the Novikov-Shubin invariants are positive. In \cite{CFM}, Carey, Farber and Mathai showed that the condition on the vanishing of the $L^2$-Betti numbers can be relaxed. In \cite{BFKM}, Burghelea, Friedlander, Kappeller and Mcdonald [8] proved the equality between the $L^2$-Reidemeister torsion and $L^2$-Ray-Singer torsion for unitary representations, under the ``determinant class condition". In \cite{Zhang}, under the framework of \cite{BCFM}, Zhang proved the the equality between the $L^2$-Reidemeister torsion and $L^2$-Ray-Singer torsion for arbitrary flat vector bundle and without the ``determinant class condition".

In \cite{Gong;CoverTorsion}, Gong and Rothenberg defined the $L^2$-analytic torsion form
and proved that the torsion form is smooth,
under the condition that the Novikov-Shubin invariant is at least half of the dimension of the base manifold. In \cite{Schick;NonCptTorsionEst}, Azzali, Goette and Schick proved
that the integrand defining the $L^2$-analytic torsion form,
as well as several other invariants related to the signature operator,
converges provided the Novikov-Shubin invariants are positive (or of determinant class and $L^{2}$-acyclic).
However, they did not prove the smoothness of the $ L ^2 $-analytic torsion form.
To consider transgression formula, they had to use weak derivatives. In \cite{SS}, So and Su proved that under the condition that the Novikov-Shubin invariants are positive the $L^2$-analytic torsion is a smooth form. Then it is natural to define the combinatorial $L^2$-torsion form and study the Cheeger-M\"{u}ller theorem for the $L^2$-analytic torsion form and the combinatorial $L^2$-torsion form. In this paper, we will use the methods in \cite{BG} and \cite{Zhang} to prove the theorem.

The rest of the paper is organized as follows. In Section 2, we recall the definition of the $L^2$-analytic torsion form. In Section 3, we recall the family of Thom-Smale complexes and the $L^2$ torsion form of it. In Section 4, we prove the main theorem of this paper. In Section 5, we extend \cite[Section 10]{BG} to the $L^2$-case and prove Theorem \ref{t4.7}. In Section 6, we prove Theorem \ref{t3}. In Appendixes, we generalize some results of \cite{B,BLe} to the current case, which are needed in this paper.

\section{$L^2$-analytic torsion form}
\setcounter{equation}{0}
In this section we will recall the definition of the $L^2$-analytic torsion form.

Let $X\to M\xrightarrow{\pi} S$ to be a fiber bundle with oriented closed connected fiber $X_{s}=\pi^{-1}(s)$. We assume that $S$ is compact. Let $g^{TX}$ be a Riemannian metric on $X$. Let $\widetilde{X}$ be a Galois covering of $X$ with covering group $\Gamma$. Let $\widetilde{X}\to \widetilde{M}\xrightarrow{\widetilde{\pi}} S$ be a fiber bundle. Let $g^{T\widetilde{X}}$ be the lifted Riemannian metric. We assume that there is a Riemannian metric $g^{TS}$ on $S$. Choose a horizontal subbundle $T^{H}\widetilde{M}$, so that
\begin{align}\label{2.1}
 T\widetilde{M}=T\widetilde{X}\oplus T^{H}\widetilde{M}.
 \end{align} 
 Let $P^{T\widetilde{X}}:T\widetilde{M}\to T\widetilde{X}$ be the projection. Then we have a Riemannian metric $g^{T\widetilde{M}}$ on $\widetilde{M}$ defined by
\begin{align}
g^{T\widetilde{M}}=g^{T\widetilde{X}}\oplus\widetilde{\pi}^* g^{TS}.
\end{align}

Let $\nabla^{T\widetilde{M},L}$ be the Levi-Civita connection on $(T\widetilde{M}, g^{T\widetilde{M}})$. Let $
\nabla^{T\widetilde{X}}=P^{T\widetilde{X}}\nabla^{T\widetilde{M},L}$ be the connection on $T\widetilde{X}$. Let $\nabla^{T\widetilde{M}}=\widetilde{\pi}^*\nabla^{TS}\oplus \nabla^{T\widetilde{X}}$ be the connection on $T\widetilde{M}$.

Put
\begin{align}\nonumber
\widetilde{S}=\nabla^{T\widetilde{M},L}-\nabla^{T\widetilde{M}}.
\end{align}

Observe that 
\begin{align}\label{2.2}
T^H\widetilde{M}\cong \widetilde{\pi}^*TS.
\end{align}
By (\ref{2.1}) and (2.2), we have the identification of bundles of algebras
\begin{align}\label{2.4}
\Lambda^{\bullet}(T^*\widetilde{M})\cong \pi^*\Lambda^{\bullet}(T^*S)\widehat{\otimes}\Lambda^{\bullet}(T^*\widetilde{X}).
\end{align}

Let $F\to M$ be flat vector bundle on $M$ with flat connection $\nabla^F$ and Hermitian metric $g^F$. Let $\widetilde{F}\to \widetilde{M}$ be the lifted bundle of $F$ with lifted flat connection $\nabla^{\widetilde{F}}$ and lifted metric $g^{\widetilde{F}}$. Set $\omega(\nabla^{\widetilde{F}},g^{\widetilde{F}})=(g^{\widetilde{F}})^{-1}\nabla^{\widetilde{F}}g^{\widetilde{F}}$ and
\begin{align}\nonumber
\nabla^{\widetilde{F},u}=\nabla^{\widetilde{F}}+{1\over 2}\omega\left(\nabla^{\widetilde{F}},g^{\widetilde{F}}\right). 
\end{align}

We equip $\Omega^{\bullet}(\widetilde{X},\widetilde{F}|_{\widetilde{X}})$ with the Hermitian product such that if $s,s'\in\Omega^{\bullet}_{c}(\widetilde{X},\widetilde{F}|_{\widetilde{X}})$,
\begin{align}\label{ad2.5}
\left\langle s,s'\right\rangle=\int_{\widetilde{X}}\langle s,s'\rangle_{\Lambda^{\bullet}(T^*\widetilde{X})\widehat{\otimes}\widetilde{F}}dv_{\widetilde{X}}.
\end{align}
Let $g^{\Omega_{(2)}^{\bullet}(\widetilde{X},\widetilde{F}|_{\widetilde{X}})}$ be the corresponding $L^2$ metric. Let $(\Omega^{\bullet}_{(2)}(\widetilde{X},\widetilde{F}|_{\widetilde{X}}),d^{\widetilde{X}})$ be the fiberwise de Rham complex of $L^{2}$ smooth forms along the fibers $\widetilde{X}$ with values in $\widetilde{F}|_{\widetilde{X}}$, equipped with the fiberwise de Rham operator $d^{\widetilde{X}}$. Let $H^{\bullet}_{(2)}(\widetilde{X},\widetilde{F}|_{\widetilde{X}})$ be the reduced $L^2$-cohomology of the complex $(\Omega^{\bullet}_{(2)}(\widetilde{X},\widetilde{F}|_{\widetilde{X}}),d^{\widetilde{X}})$. Let $g^{H^{\bullet}_{(2)}(\widetilde{X},\widetilde{F}|_{\widetilde{X}})}_{L_2}$ be the induced metric on $H^{\bullet}_{(2)}(\widetilde{X},\widetilde{F}|_{\widetilde{X}})$.

Let $\Omega^{\bullet}(S,\Omega_{(2)}^{\bullet}(\widetilde{X},\widetilde{F}|_{\widetilde{X}}))$ be the space of smooth sections of $\Lambda^{\bullet}(T^*S)\widehat{\otimes}\Omega^{\bullet}_{(2)}(\widetilde{X},\widetilde{F}|_{\widetilde{X}})$ on $S$. Using (\ref{2.4}), we have the identification of $\bf Z$-graded vector spaces,
\begin{align}
\Omega^{\bullet}_{(2)}\left(\widetilde{M},\widetilde{F}\right)\cong \Omega^{\bullet}\left(S,\Omega_{(2)}^{\bullet}\left(\widetilde{X},\widetilde{F}|_{\widetilde{X}}\right)\right).
\end{align}

The operator $d^{\widetilde{M}}$ acting on $\Omega^{\bullet}_{(2)}(\widetilde{M},\widetilde{F})$ has degree 1 and is such that $d^{\widetilde{M},2}=0$. Let $\nabla^{\Omega_{(2)}^{\bullet}(\widetilde{X},\widetilde{F}|_{\widetilde{X}})}$ be the connection on $\Omega^{\bullet}_{(2)}(\widetilde{X},\widetilde{F}|_{\widetilde{X}})$, such that if $U\in TS$ and $s\in \Omega^{\bullet}_{(2)}(\widetilde{X},\widetilde{F}|_{\widetilde{X}})$, then
\begin{align}\label{nn1}
\nabla_{U}^{\Omega^{\bullet}_{(2)}(\widetilde{X},\widetilde{F}|_{\widetilde{X})}}s=L_{U^H}s,
\end{align}
where $U^H\in C^{\infty}(\widetilde{M},T^H\widetilde{M})$ such that $\widetilde{\pi}_* U^H=U$.

Let $\nabla^{H^{\bullet}_{(2)}(\widetilde{X},\widetilde{F}|_{\widetilde{X}})}$ be the induced connection on $H^{\bullet}_{(2)}(\widetilde{X},\widetilde{F}|_{\widetilde{X}})$.

If $U,V$ are smooth sections of $TS$, set 
\begin{align}\label{ad1}
T^H (U,V)=-P^{T\widetilde{X}}[U^H,V^H]. 
\end{align}
Then $\widetilde{A}'=d^{\widetilde{M}}$ can be decomposed as $\widetilde{A}'=d^{\widetilde{X}}+\nabla^{\Omega_{(2)}^{\bullet}(\widetilde{X},\widetilde{F}|_{\widetilde{X}})}+i_{T^H}$. Let $\widetilde{A}''$ be the adjoint of the superconnection $\widetilde{A}'$ with respect to the metric $g^{\Omega^{\bullet}(\widetilde{X},\widetilde{F}|_{\widetilde{X}})}$. Let $\nabla^{\Omega_{(2)}^{\bullet}(\widetilde{X},\widetilde{F}|_{\widetilde{X}}),*}$ be the connection on $\Omega_{(2)}^{\bullet}(\widetilde{X},\widetilde{F}|_{\widetilde{X}})$ which is adjoint to $\nabla^{\Omega_{(2)}^{\bullet}(\widetilde{X},\widetilde{F}|_{\widetilde{X}})}$ with respect to $g^{\Omega^{\bullet}(\widetilde{X},\widetilde{F}|_{\widetilde{X}})}$. Then we have $\widetilde{A}''=d^{\widetilde{X},*}+\nabla^{\Omega^{\bullet}(\widetilde{X},\widetilde{F}|_{\widetilde{X}}),*}-T^H\wedge $.

Set 
\begin{align}\nonumber
D^{\widetilde{X}}=d^{\widetilde{X}}+d^{\widetilde{X},*},\ D^{\widetilde{X},2}=\left(d^{\widetilde{X}}+d^{\widetilde{X},*}\right)^2=d^{\widetilde{X}}d^{\widetilde{X},*}+d^{\widetilde{X},*}d^{\widetilde{X}}. 
\end{align}
Let $D^{\widetilde{X},2}_i=D^{\widetilde{X},2}|_{\Omega^i(\widetilde{X},\widetilde{F}|_{\widetilde{X}})}$ and  $\widetilde{P}_{{\ker}D_i^{\widetilde{X},2}}$ be the orthogonal projection onto ${\rm ker}D^{\widetilde{X},2}_i$

\begin{definition} 
The analytic Novikov-Shubin invariant (\cite{Lu}, \cite{NS1}, \cite{NS2}) is defined by 
\begin{align}\nonumber
\alpha_i=\sup\left\{\beta\geq 0|{\rm Tr}_{\Gamma}\left(\exp\left(-tD_{i}^{\widetilde{X},2}\right)\right)-{\rm Tr}_{\Gamma}\left(\widetilde{P}_{{\rm ker}D^{\widetilde{X},2}_i}\right)=\mathcal{O}\left(t^{-{\beta}}\right)\right\},\ i=0,\cdots,{\rm dim}\widetilde{X}.
\end{align}
\end{definition}

Set 
\begin{align}
\widetilde{A}={1\over 2}(\widetilde{A}''+\widetilde{A}'),\ \ \widetilde{B}={1\over 2}(\widetilde{A}''-\widetilde{A}').
\end{align}

For $t>0$, set
\begin{align}
g_{t}^{T\widetilde{X}}={{g^{T\widetilde{X}}}\over{t}}.
\end{align}

Let $N$ be the number operator of $\Omega_{(2)}^{\bullet}(\widetilde{X},\widetilde{F}|_{\widetilde{X}})$, i.e. $N$ acts by multiplication by $k$ on $\Omega^{k}_{(2)}(\widetilde{X},\widetilde{F}|_{\widetilde{X}})$. One verifies easily that 
\begin{align}
g_t^{\Omega_{(2)}^{\bullet}(\widetilde{X},\widetilde{F}|_{\widetilde{X}})}=t^{N-n/2}g^{\Omega_{(2)}^{\bullet}(\widetilde{X},\widetilde{F}|_{\widetilde{X}})}.
\end{align}

Let $\widetilde{A}''_t$ be the adjoint of $\widetilde{A}'$ with respect to $g_t^{\Omega_{(2)}^{\bullet}(\widetilde{X},\widetilde{F}|_{\widetilde{X}})}$, then we have
\begin{align}
\widetilde{A}''_t=t^{-N}\widetilde{A}'' t^N.
\end{align}
Set 
\begin{align}
\widetilde{A}_t={1\over 2}(\widetilde{A}''_t +\widetilde{A}'),\ \ \widetilde{B}_t={1\over 2}(\widetilde{A}''_t-\widetilde{A}').
\end{align}
For $t>0$, set 
\begin{align}
\widetilde{C}'_t=t^{N/2}\widetilde{A}' t^{-N/2},\ \ \widetilde{C}''_t=t^{-N/2}\widetilde{A}'' t^{N/2}.
\end{align}
Then $\widetilde{C}'_t$ is a flat superconnection on $\Omega_{(2)}^{\bullet}(\widetilde{X},\widetilde{F}|_{\widetilde{X}})$, and $\widetilde{C}''_t$ is its adjoint with respect to $g^{\Omega_{(2)}^{\bullet}(\widetilde{X},\widetilde{F}|_{\widetilde{X}})}$. Set
\begin{align}\label{edp2.15}
\widetilde{C}_t={1\over 2}(\widetilde{C}''_t +\widetilde{C}'_t),\ \ \widetilde{D}_t={1\over 2}(\widetilde{C}''_t-\widetilde{C}'_t).
\end{align}
Then we have
\begin{align}\label{ad2.16}
\widetilde{C}_t=t^{N/2}\widetilde{A}_t t^{-N/2},\ \ \widetilde{D}_t=t^{N/2}\widetilde{B}_t t^{-N/2}.
\end{align}

For $t>0$, let $\psi_{t}:\Lambda^{\bullet}(T^*S)\to \Lambda^{\bullet}(T^*S)$ be given by
\begin{align}
\psi_t \omega=t^{{\rm deg}\omega/2}\omega.
\end{align}
We have 
\begin{align}\label{ad2.18}
\widetilde{C}_t=\psi^{-1}_t \sqrt{t}\widetilde{A}\psi_t,\ \ \widetilde{D}_t=\psi^{-1}_t \sqrt{t}\widetilde{B}\psi_t.
\end{align}

We fix a square root $i^{1/2}$ of $i$. Our formulas will not depend on the choice of the square root. Let $\varphi: \Lambda^{\bullet}(T^*S)\to \Lambda^{\bullet}(T^*S)$ be the linear map such that for all homogeneous $\omega\in\Lambda^{\bullet}(T^*S)$,
\begin{align}
\varphi\omega =(2i\pi)^{-{\rm deg}(\omega)/2}\omega.
\end{align}

Set $h(x)=xe^{x^2}$. For $t>0$, set 
\begin{align}
h^{\wedge}\left(\widetilde{A}',g_{t}^{\Omega_{(2)}^{\bullet}(\widetilde{X},\widetilde{F}|_{\widetilde{X}})}\right)=\varphi{\rm Tr}_{\Gamma,s} \left[{N\over 2}h'(\widetilde{B}_t)\right],
\end{align}
where $\Gamma$ means the $\Gamma$-trace (cf. \cite{A}) and $s$ denotes the supertrace in the sense of \cite{Q}.

Put 
\begin{align}
\chi'(\widetilde{F})=\sum_{j=0}^{m}(-1)^j j{\rm dim}_{\Gamma}\left(H^{j}_{(2)}\left(\widetilde{X},\widetilde{F}|_{\widetilde{X}}\right)\right),
\ \chi(\widetilde{F})={\rm dim}(F)e\left(X\right).
\end{align}
Here ${\rm dim}_{\Gamma}$ means the $\Gamma$ dimension. Then $\chi'(F)$ and $\chi(F)$ are locally constant functions on $S$.

In the following of this paper, we will assume that the analytic Novikov-Shubin invariants of $(X,F)$ are positive. 
\begin{definition}
The $L^2$-analytic torsion form is defined by
\begin{multline}\label{ad2.22}
\mathcal{T}_{L^2,h}(T^H \widetilde{M},g^{T\widetilde{X}},\nabla^{\widetilde{F}},g^{\widetilde{F}})=-\int_{0}^{+\infty}\left[h^{\wedge}\left(\widetilde{A}',g_{t}^{\Omega_{(2)}^{\bullet}(\widetilde{X},\widetilde{F}|_{\widetilde{X}})}\right)-{1\over 2}\chi'(\widetilde{F})h'(0)\right.\\
-\left.\left({n\over 4}\chi(\widetilde{F})-{1\over 2}\chi'(\widetilde{F})\right)h'\left(i\sqrt{t}/2\right)\right]{{dt}\over{t}}.
\end{multline}
\end{definition}
By \cite{SS}, the $L^2$-analytic torsion form $\mathcal{T}_{L^2,h}(T^H \widetilde{M},g^{T\widetilde{X}},\nabla^{\widetilde{F}},g^{\widetilde{F}})\in\Omega^{\bullet}(S)$ is a smooth form.

\section{A family of Thom-Smale gradient vector fields}
\setcounter{equation}{0}

Let $X$ be a compact manifold of dimension $n$. Let $f:X\to \mathbb{\bf R}$ be a Morse function. Let $B$ be the set of critical points of $f$,
\begin{align}
B=\{x\in X, df(x)=0\}.
\end{align}
If $x\in B$, recall that the index ${\rm ind}(x)$ is such that the quadratic form $d^2 f(x)$ on $T_x X$ has signature $(n-{\rm ind}(x),{\rm ind}(x))$.

Let $f:M\to {\bf R}$ be a smooth function. We assume that $f$ is Morse along every fiber $X$.

Let $h^{TX}$ be a metric on $TX$. Let $\nabla f\subset TX$ be the gradient field of $f$ along the fiber $X$ with respect to $h^{TX}$. Then $\nabla f\in TX$. We make the fundamental assumption that $Y=-\nabla f$ is Thom-Smale along every fiber $X$.

Let ${\bf B}$ be the zero set of $Y$, i.e. the set of fiberwise critical points of $f$. Let ${\bf B}^i$ be the set of critical points of $f$ which have index $i$ along fibers $X$. Then ${\bf B}$, ${\bf B}^i$ are finite covers of $S$. We denote by $B$, $B^i$ the corresponding fibers.

Let $\widetilde{f}$ be the lifted of $f$ to $\widetilde{M}$, then the restriction of $\widetilde{f}$ to the fiber $\widetilde{X}$ is Morse. Let $\widetilde{\bf B}$ be the zero set of $\nabla \widetilde{f}$ and $\widetilde{B}$ be its fiber. Let $\widetilde{B}^i$ be the set of critical points of $\widetilde{f}$ which have index $i$ along fibers of $\widetilde{X}$. For $x\in \widetilde{\bf B}$, let $\widetilde{W}^u (x)$ and $\widetilde{W}^s (x)$ be the unstable and stable cells. If $x\in \widetilde{\bf B}$, set 
$$T_x \widetilde{X}^u=T_x \widetilde{W}^u (x),\ \ T_x \widetilde{X}^s=T_x \widetilde{W}^s (x).$$
Let $T\widetilde{X}^u$ and $T\widetilde{X}^s$ be the vector bundles on $\widetilde{\bf B}$ with fibers $T_x \widetilde{X}^u$ and $T_x \widetilde{X}^s$ respectively, and we have
\begin{align}\label{s3.2}
T\widetilde{X}|_{\widetilde{\bf B}}=T\widetilde{X}^u\oplus T\widetilde{X}^s.
\end{align}
Let $o^u$, $o^s$ be the ${\bf Z}_2$-lines on $\widetilde{B}$, which are the orientation lines of $T\widetilde{X}^u$, $T\widetilde{X}^s$.

Let $(F,\nabla^F)$ be a complex flat vector bundle over $X$ with flat connection $\nabla^F$, and let $(F^*,\nabla^{F^*})$ be the corresponding dual flat vector bundle carrying the flat connection $\nabla^{F^*}$ and the dual Hermitian metric $g^{F^*}$. Let $(\widetilde{F},\nabla^{\widetilde{F}})$ denote the lifted flat vector bundle over $\widetilde{X}$ obtained as the pullback of $(F,\nabla^F)$ through the covering map. Let $g^{\widetilde{F}}$ be the naturally lifted Hermitian metric on $\widetilde{F}$. Let $(\widetilde{F}^*,\nabla^{\widetilde{F}^*})$ and $g^{\widetilde{F}^*}$ denote the corresponding lifted objects on $\widetilde{X}$.

Set 
\begin{align}\nonumber
C_{\bullet}(W^u,\widetilde{F}^*)=\bigoplus_{x\in\widetilde{B}}\widetilde{F}^*_{x}\otimes o^u_x,
\end{align}
\begin{align}\label{3.3}
C_{i}(W^u,\widetilde{F}^*)=\bigoplus _{x\in\widetilde{B}^i}\widetilde{F}^*_{x}\otimes o^u_x,
\end{align}
where $\bigoplus$ is meant in the $L^2$ sense. Then we have the $L^2$-Thom-Smale complex $(C_{\bullet}(W^u,\widetilde{F}^*),\partial)$ (cf. \cite{Zhang}).

Let $(C^{\bullet}(W^u,\widetilde{F}),\widetilde{\partial})$ be the complex dual to the complex $(C_{\bullet}(W^u,\widetilde{F}^*),\partial)$. By (\ref{3.3}), we get
\begin{align}\nonumber
C^{\bullet}(W^u,\widetilde{F})=\bigoplus_{x\in \widetilde{B}}\widetilde{F}_x\otimes o^u_x,
\end{align}
\begin{align}
C^i(W^u,\widetilde{F})=\bigoplus_{x\in\widetilde{B}^i}\widetilde{F}_x\otimes o^u_x,
\end{align}
where $\bigoplus$ is meant in the $L^2$ sense.

Let $l^2(\Gamma)$ denote the Hilbert space obtained through the $L^2$-completion of the group algebra of $\Gamma$ with respect to the canonical trace on it. Then $C^i(W^u,\widetilde{F})$ is a Hilbert space, which is isomorphic to the direct sum of $n_i$ copies of $l^2(\Gamma)$, where $n_i=\# B^i$. Then we call $(C^{\bullet}(W^u,\widetilde{F}),\widetilde{\partial}))$ the $L^2$-Thom-Smale cochain complex.

Then $C^{\bullet}(W^u,\widetilde{F})$ is a flat ${\bf Z}$-graded vector bundle on $S$. Let $\nabla^{C^{\bullet}(W^u,\widetilde{F})}$ be the corresponding flat connection on $C^{\bullet}(W^u,\widetilde{F})$. Then 
\begin{align}\nonumber
\widetilde{A}^{C^{\bullet}(W^u,\widetilde{F})'}=\widetilde{\partial}+\nabla^{C^{\bullet}(W^u,\widetilde{F})}
\end{align}
 is a flat superconnection of total degree 1 on $C^{\bullet}(W^u,\widetilde{F})$ and we have
\begin{align}
H^{\bullet}_{(2)}\left(C^{\bullet}\left(W^u,\widetilde{F}\right),\widetilde{\partial}\right)\cong H^{\bullet}_{(2)}\left(\widetilde{X},\widetilde{F}|_{\widetilde{X}}\right).
\end{align}
Let $g^{H^{\bullet}_{(2)}(\widetilde{X},\widetilde{F}|_{\widetilde{X}})}_{C^{\bullet}({W}^u,\widetilde{F})}$ be the metric on $H^{\bullet}_{(2)}(\widetilde{X},\widetilde{F}|_{\widetilde{X}})$ induced by the isomorphism.

Let $\widetilde{\partial}^*$ be the adjoint of $\widetilde{\partial}$ and $\nabla^{C^\bullet(W^u,\widetilde{F}),*}$ be the adjoint of $\nabla^{C^\bullet(W^u,\widetilde{F})}$. For $t>0$, set
\begin{align}\nonumber
\widetilde{A}^{C^\bullet(W^u,\widetilde{F})''}_t=t\widetilde{\partial}^* +\nabla^{C^{\bullet}(W^u,\widetilde{F}),*}, {C}^{C^\bullet(W^u,\widetilde{F})'}_t=\sqrt{t}\widetilde{\partial}+\nabla^{C^{\bullet}(W^u,\widetilde{F})},
\end{align}
\begin{align}
{C}^{C^\bullet(W^u,\widetilde{F})''}_t=\sqrt{t}\widetilde{\partial}^*+\nabla^{C^\bullet(W^u,\widetilde{F}),*},\ {B}^{C^\bullet(W^u,\widetilde{F})}_t={1\over 2}\left(\widetilde{A}^{C^\bullet(W^u,\widetilde{F})''}_t-\widetilde{A}^{C^{\bullet}(W^u,\widetilde{F})'}\right).
\end{align}
Then as \cite[Definition 1.17]{BG}, we define
\begin{align}\nonumber
h^{\wedge}\left(\widetilde{A}^{C^{\bullet}(W^u,\widetilde{F})'},g_t^{C^\bullet(W^u,\widetilde{F})}\right)=\varphi{\rm Tr}_{\Gamma,s}\left[{N^{C^\bullet(W^u,\widetilde{F})}\over 2}h'\left(B^{C^\bullet(W^u,\widetilde{F})}_t\right)\right].
\end{align}

Set 
\begin{align}\nonumber
\widetilde{D}^c=\widetilde{\partial}+\widetilde{\partial}^*,\ \widetilde{D}^{c,2}=\left(\widetilde{\partial}+\widetilde{\partial}^*\right)^2=\widetilde{\partial}\widetilde{\partial}^*+\widetilde{\partial}^*\widetilde{\partial},\ \widetilde{D}^{c,2}_i=\widetilde{D}^{c,2}|_{C^i(W^u,\widetilde{F})}. 
\end{align}
Let $\widetilde{P}_{{\rm ker}\widetilde{D}^{c,2}_i}$ be the orthogonal projection onto ${\rm ker}\widetilde{D}^{c,2}_i$. 

\begin{definition}
The cellular Novikov-Shubin invariant (cf. \cite{Lu}) is defined by
\begin{align}\nonumber
\alpha^c_i=\sup\left\{\beta\geq 0|{\rm Tr}_{\Gamma}\left(\exp\left(-t\widetilde{D}^{c,2}_i\right)\right)-{\rm Tr}_{\Gamma}\left(\widetilde{P}_{{\rm ker}\widetilde{D}_i^{c,2}}\right)=\mathcal{O}\left(t^{-{\beta}}\right)\right\},\ i=0,\cdots,{\rm dim}\widetilde{X}.
\end{align}
\end{definition}

\begin{theorem}(\cite[Theorem 2.68]{Lu})
Let $M$ be a cocompact free proper $G$-manifold without boundary and with $G$-invariant Riemannian metric. Then the cellular and the analytic spectral density functions are dilatationally equivalent and the cellular and analytic Novikov-Shubin invariants agree in each dimension $p$. 
\end{theorem}

Since we assume that $\alpha_i>0, i=0,\cdots,{\rm dim}\widetilde{X}$, then by \cite{E,GS,Lu}, we have $\alpha_i^c>0, i=0,\cdots, {\rm dim}\widetilde{X}$. Set
\begin{align}\nonumber
\chi^c (\widetilde{F})=\sum_{i=0}^n (-1)^i i {\rm dim}_{\Gamma}\left(C^i (W^u,\widetilde{F})\right).
\end{align}

As \cite[Proposition 1.27]{BG}, we have 
\begin{proposition}\label{edp3.2}
\begin{align}
h^{\wedge}\left(\widetilde{A}^{C^{\bullet}(W^u,\widetilde{F})'},g_t^{C^\bullet(W^u,\widetilde{F})}\right)=
{1\over 2}\chi^c (\widetilde{F})+\mathcal{O}(t),\ as\ t\to 0.
\end{align}
\end{proposition}

\begin{proof}
The follows similarly as (\ref{su4.18}), Theorem \ref{det4.5} and (\ref{de4.34}) below. 
\end{proof}

Since $\alpha^c_i>0, i=0,\cdots,{\rm dim}\widetilde{X}$, using Duhamel expansion developed in \cite[Theorem 4.1]{Schick;NonCptTorsionEst}, one has
\begin{theorem}\label{edt3.3}
There exists $\gamma>0$ such that
\begin{align}\nonumber
h^{\wedge}\left(\widetilde{A}^{C^{\bullet}(W^u,\widetilde{F})'},g_t^{C^\bullet(W^u,\widetilde{F})}\right)={1\over 2}\chi'(\widetilde{F})+\mathcal{O}(t^{-\gamma}),\ as\ t\to \infty.
\end{align}
\end{theorem}

\begin{proof}
As in \cite{Schick;NonCptTorsionEst}, set 
$$\theta(t)={\rm Tr}_{\Gamma}\left[\exp \left(-t\widetilde{D}^{c,2}\right)\right]-{\rm Tr}_{\Gamma}\left[\widetilde{P}_{{\rm ker}\widetilde{D}^{c,2}}\right].$$
Then by the same proof of \cite[Lemma 4.5]{Schick;NonCptTorsionEst}, there exists a monotone decaying function $\bar{s}=\bar{s}(t)$ and $\alpha>0$ such that
\begin{align}\label{r23.8}
\theta\left({t\bar{s}(t)}\right)\cdot \left({1\over{\bar{s}(t)}}\right)^{{{\rm dim} B}\over 2}\leq t^{-\alpha},\ {\rm as}\ t\to \infty.
\end{align}
Choose $T$ such that $\bar{s}(T)< {1\over{{\rm dim}B+1}}$. Then as \cite[Lemma 4.2]{Schick;NonCptTorsionEst}, for $k\geq 0$, there exists a constant $C$ such that for all $s>0$, $t>T$, 
\begin{align}\nonumber
\left\|\left(\sqrt{t}\widetilde{D}^c\right)^k e^{-st\widetilde{D}^{c,2}}\right\|_{\rm op}\leq C s^{-{c\over 2}},\ \ for\ k\geq 0,
\end{align}
\begin{align}\nonumber
\left\| e^{-st \widetilde{D}^{c,2}}-\widetilde{P}_{{\rm ker}\widetilde{D}^{c,2}}\right\|_{\Gamma,1}=\theta(st),
\end{align}
\begin{align}\label{r23.9}
\left\|\left(\sqrt{t}\widetilde{D}^c\right)^k e^{-st\widetilde{D}^{c,2}}\right\|_{\Gamma,1}\leq C s^{-{k\over 2}}\theta\left({st\over 2}\right),\ \ for\ k\geq 1,
\end{align}
where $\|A\|_{\Gamma,1}={\rm Tr}_{\Gamma}(|A|)$ (cf. Definition \ref{Ade5}).

As (\ref{edp2.15}) and (\ref{ad2.16}), one has
$$h^{\wedge}\left(\widetilde{A}^{C^{\bullet}(W^u,\widetilde{F})'},g_t^{C^\bullet(W^u,\widetilde{F})}\right)=\varphi {\rm Tr}_{\Gamma,s}\left[{N^{C^\bullet(W^u,\widetilde{F})}\over 2}h'\left({{C}^{C^\bullet(W^u,\widetilde{F})''}_t-{C}^{C^\bullet(W^u,\widetilde{F})'}_t}\over 2\right)\right],$$
where
$${{{C}^{C^\bullet(W^u,\widetilde{F})''}_t-{C}^{C^\bullet(W^u,\widetilde{F})'}_t}\over 2}={1\over 2}\sqrt{t}\left(\widetilde{\partial}^*-\widetilde{\partial}\right)+{1\over 2}\left(\nabla^{C^\bullet(W^u,\widetilde{F}),*}-\nabla^{C^\bullet(W^u,\widetilde{F})}\right).$$

Set 
$$\mathbb{X}_{t}^{C^\bullet(W^u,\widetilde{F})}={{{C}^{C^\bullet(W^u,\widetilde{F})''}_t-{C}^{C^\bullet(W^u,\widetilde{F})'}_t}\over 2},\ \Omega^{C^\bullet(W^u,\widetilde{F})}={1\over 2}\left(\nabla^{C^\bullet(W^u,\widetilde{F}),*}-\nabla^{C^\bullet(W^u,\widetilde{F})}\right),$$
$$\omega^{C^\bullet(W^u,\widetilde{F})}={1\over 2}\left(\widetilde{\partial}^*-\widetilde{\partial}\right).$$

We denote the standard $n$-simplex by
$$\Delta^n=\left\{ (s_0,\dots,s_n)\in [0,1]^{n+1}|s_0+\cdots s_n=1\right\}$$
and the standard volume form on $\Delta^n$ by $d^n (s_0,\dots,s_n)$, so that $\Delta^n$ has total volume $1\over{n!}$.

Since
$$\mathbb{X}_{t}^{C^\bullet(W^u,\widetilde{F})}=\sqrt{t}\omega^{C^\bullet(W^u,\widetilde{F})}+\Omega^{C^\bullet(W^u,\widetilde{F})},$$
we have
$$\left(\mathbb{X}_{t}^{C^\bullet(W^u,\widetilde{F})}\right)^2=-t \widetilde{D}^{c,2}+\sqrt{t} \omega^{C^\bullet(W^u,\widetilde{F})}\Omega^{C^\bullet(W^u,\widetilde{F})}+\sqrt{t}\Omega^{C^\bullet(W^u,\widetilde{F})}\omega^{C^\bullet(W^u,\widetilde{F})}+\left(\Omega^{C^\bullet(W^u,\widetilde{F})}\right)^2.$$
Then by Duhamel's principle, we have
\begin{multline}\label{r23.11}
\exp\left(\left(\mathbb{X}_{t}^{C^\bullet(W^u,\widetilde{F})}\right)^2\right)=\sum_{n=0}^{{\rm dim}B}\int_{\Delta^n}e^{-s_0 t \widetilde{D}^{c,2}}\left(\sqrt{t} \omega^{C^\bullet(W^u,\widetilde{F})}\Omega^{C^\bullet(W^u,\widetilde{F})}+\right.\\
\left.\sqrt{t}\Omega^{C^\bullet(W^u,\widetilde{F})}\omega^{C^\bullet(W^u,\widetilde{F})}+\left(\Omega^{C^\bullet(W^u,\widetilde{F})}\right)^2\right)e^{-s_1 t \widetilde{D}^{c,2}}\cdots \\
\cdots \left(\sqrt{t} \omega^{C^\bullet(W^u,\widetilde{F})}\Omega^{C^\bullet(W^u,\widetilde{F})}+\sqrt{t}\Omega^{C^\bullet(W^u,\widetilde{F})}\omega^{C^\bullet(W^u,\widetilde{F})}+\left(\Omega^{C^\bullet(W^u,\widetilde{F})}\right)^2\right)e^{-s_n t\widetilde{D}^{c,2}}\\
d^n (s_0,\dots,s_n). 
\end{multline}
Then by (\ref{r23.8}), (\ref{r23.9}) and (\ref{r23.11}), proceeding as in \cite{Schick;NonCptTorsionEst}, 
we have 
\begin{multline}\label{r23.11}
\lim_{t\to \infty}\left(\mathbb{X}_{t}^{C^\bullet(W^u,\widetilde{F})}\right)^k \exp\left(\left(\mathbb{X}_{t}^{C^\bullet(W^u,\widetilde{F})}\right)^2\right)\\=\widetilde{P}_{{\rm ker}\widetilde{D}^c}\left(\Omega^{C^\bullet(W^u,\widetilde{F})} \widetilde{P}_{{\rm ker}\widetilde{D}^c}\right)^k \exp\left(\left(\Omega^{C^\bullet(W^u,\widetilde{F})} \widetilde{P}_{{\rm ker}\widetilde{D}^c}\right)^2\right),
\end{multline}
in $\|\cdot\|_{\Gamma,1}$, for $k=0,1,2$.

Since 
\begin{align}\label{r23.12}
\widetilde{P}_{{\rm ker}\widetilde{D}^c}\Omega^{C^\bullet (W^u,\widetilde{F})}\widetilde{P}_{{\rm ker}\widetilde{D}^c}={1\over 2}\left(\nabla^{H^\bullet_{(2)}\left(\widetilde{X},\widetilde{F}|_{\widetilde{X}}\right),*}-\nabla^{H^\bullet_{(2)}\left(\widetilde{X},\widetilde{F}|_{\widetilde{X}}\right)}\right).
\end{align}
Then by (\ref{r23.11}) and (\ref{r23.12}), we get the theorem. 

\end{proof}

\begin{definition}
By Proposition \ref{edp3.2} and Theorem \ref{edt3.3}, the $L^2$-combinatorial torsion form is defined by 
\begin{multline}
T_{L^2,h}(\widetilde{A}^{C^{\bullet}(W^u,\widetilde{F})'},g^{C^{\bullet}(W^u,\widetilde{F})})=-\int_{0}^{+\infty}\left[h^{\wedge}\left(\widetilde{A}^{C^{\bullet}(W^u,\widetilde{F})'},g_t^{C^\bullet(W^u,\widetilde{F})}\right)-{1\over 2}\chi'(\widetilde{F})h'(0)\right.\\
-\left.\left({1\over 2}\chi^c (\widetilde{F})-{1\over 2}\chi'(\widetilde{F})\right)h'\left(i\sqrt{t}/2\right)\right]{{dt}\over{t}}.
\end{multline}
\end{definition}

As in \cite{SS}, one can also define the $m$-th Hilbert-Schmidt norm $\|\cdot\|_{{\rm HS}\ m}$ and $m$-th operator norm $\|\cdot\|_{{\rm op}\ m}$ for this case. Since $\alpha_i^c>0, i=0,\cdots,{\rm dim}\widetilde{X}$, as in \cite[Theorem 4.4]{SS}, we have that (\ref{r23.11}) holds in $\|\cdot\|_{{\rm HS}\ m}$. Then under the condition of positive Novikov-Shubin invariant, $T_{L^2,h}(\widetilde{A}^{C^{\bullet}(W^u,\widetilde{F})'},g^{C^{\bullet}(W^u,\widetilde{F})})\in \Omega(S)$ is a smooth form.

\begin{definition}
Let $\widetilde{{\bf P}}^{\infty}: \Omega^{\bullet}_{(2)}(\widetilde{M},\widetilde{F})\to \Omega^{\bullet}(S,C^{\bullet}(W^u,\widetilde{F}))$ be given by
\begin{align}
\widetilde{{\bf P}}^{\infty}\alpha=\sum_{x\in\widetilde{B}}\widetilde{W}^u(x)^*\int_{\overline{\widetilde{W}^u(x)}}\alpha.
\end{align}
\end{definition}
Then the map $\widetilde{{\bf P}}^{\infty}$ is a quasi-isomorphism of ${\bf Z}$-graded mapping $(\Omega_{(2)}^{\bullet}(\widetilde{M},\widetilde{F}),d^{\widetilde{M}})$ into the Thom-Smale complex $(\Omega^{\bullet}(S,C^{\bullet}(W^u,\widetilde{F})),\widetilde{A}^{C^{\bullet}(W^u,\widetilde{F})'})$.

\section{The main theorem}
\setcounter{equation}{0}

In this section, we will prove the main theorem of this paper.

Let $g_{\ell}$ be a path of metrics on $H^{\bullet}_{(2)}(\widetilde{X},\widetilde{F}|_{\widetilde{X}})$ connecting $g^{H^{\bullet}_{(2)}(\widetilde{X},\widetilde{F}|_{\widetilde{X}})}_{C^{\bullet}({W}^u,\widetilde{F})}$ and $g^{H^{\bullet}_{(2)}(\widetilde{X},\widetilde{F}|_{\widetilde{X}})}_{L_2}$, let $\nabla_{\ell}^{H^{\bullet}_{(2)}(\widetilde{X},\widetilde{F}|_{\widetilde{X}}),*}$ be adjoint of $\nabla^{H^{\bullet}_{(2)}(\widetilde{X},\widetilde{F}|_{\widetilde{X}})}$ with respect to $g_{\ell}$, let $$B_{\ell}={1\over 2}\left(\nabla_{\ell}^{H^{\bullet}_{(2)}(\widetilde{X},\widetilde{F}|_{\widetilde{X}}),*}-\nabla^{H^{\bullet}_{(2)}(\widetilde{X},\widetilde{F}|_{\widetilde{X}})}\right).$$
As \cite[Definition 1.10]{BG}, we define
\begin{align}\label{ad4.1}
\widetilde{h}_{L^2}\left(\nabla^{H^{\bullet}_{(2)}(\widetilde{X},\widetilde{F}|_{\widetilde{X}})},g_{\ell}\right)=\int_{0}^{1}\varphi {\rm Tr}_{\Gamma,s}\left[{1\over 2}\left(g_{\ell}\right)^{-1}{{\partial g_{\ell}}\over{\partial {\ell}}}h'(B_{\ell})\right]d{\ell}.
\end{align}
By the same argument in \cite[Theorem 1.11]{BG}, the class of $\widetilde{h}_{L^2}\left(\nabla^{H^{\bullet}_{(2)}(\widetilde{X},\widetilde{F}|_{\widetilde{X}})},g_{\ell}\right)$ in $\Omega^{\bullet}(S)/d\Omega^{\bullet}(S)$ only depends on $g^{H^{\bullet}_{(2)}(\widetilde{X},\widetilde{F}|_{\widetilde{X}})}_{C^{\bullet}({W}^u,\widetilde{F})}$ and $g^{H^{\bullet}_{(2)}(\widetilde{X},\widetilde{F}|_{\widetilde{X}})}_{L_2}$. We denote the class by 
$$\widetilde{h}_{L^2}\left(\nabla^{H^{\bullet}_{(2)}(\widetilde{X},\widetilde{F}|_{\widetilde{X}})},g^{H^{\bullet}_{(2)}(\widetilde{X},\widetilde{F}|_{\widetilde{X}})}_{C^{\bullet}({W}^u,\widetilde{F})},g^{H^{\bullet}_{(2)}(\widetilde{X},\widetilde{F}|_{\widetilde{X}})}_{L_2}\right)\in \Omega^{\bullet}(S)/d\Omega^{\bullet}(S).$$

Let $T'^{H}\widetilde{M}$, $g'^{T\widetilde{X}}$ and $g'^{\widetilde{F}}$ be another triple of data, then the following anomaly formula holds (\cite[Lemma 4.9]{SS}),
\begin{multline}
\mathcal{T}_{L^2,h}(T'^H \widetilde{M},g'^{T\widetilde{X}},\nabla^{\widetilde{F}},g'^{\widetilde{F}})-\mathcal{T}_{L^2,h}(T^H \widetilde{M},g^{T\widetilde{X}},\nabla^{\widetilde{F}},g^{\widetilde{F}})=\\
\int_{X}\widetilde{e}\left(TX,\nabla^{TX},\nabla'^{TX}\right)h\left(\nabla^F,g^F\right)+\int_{X}e\left(TX,\nabla'^{TX}\right)\widetilde{h}\left(\nabla^F,g^F,g'^F\right)\\
-\widetilde{h}_{L^2}\left(\nabla^{H^{\bullet}_{(2)}(\widetilde{X},\widetilde{F}|_{\widetilde{X}})},g_{L^2}^{H^{\bullet}_{(2)}(\widetilde{X},\widetilde{F}|_{\widetilde{X}})},g_{L^2}^{'H^{\bullet}_{(2)}(\widetilde{X},\widetilde{F}|_{\widetilde{X}})}\right)\ {\rm in}\ \Omega^{\bullet}(S)/d\Omega^{\bullet}(S).
\end{multline}
By the same argument of \cite[Theorem 1.31]{BG}, we also have
\begin{multline}
T_{L^2,h}(\widetilde{A}^{C^{\bullet}(W^u,\widetilde{F})'},g^{C^{\bullet}(W^u,\widetilde{F})})-T_{L^2,h}(\widetilde{A}^{C^{\bullet}(W^u,\widetilde{F})'},g'^{C^{\bullet}(W^u,\widetilde{F})})\\
=\widetilde{h}_{L^2}\left(\nabla^{C^\bullet(W^u,\widetilde{F})},g'^{C^\bullet(W^u,\widetilde{F})},g^{C^\bullet(W^u,\widetilde{F})}\right)\\-\widetilde{h}_{L^2}\left(\nabla^{H^{\bullet}_{(2)}(\widetilde{X},\widetilde{F}|_{\widetilde{X}})},g'^{H^{\bullet}_{(2)}(\widetilde{X},\widetilde{F}|_{\widetilde{X}})}_{C^{\bullet}(W^u,\widetilde{F})},g^{H^{\bullet}_{(2)}(\widetilde{X},\widetilde{F}|_{\widetilde{X}})}_{C^{\bullet}(W^u,\widetilde{F})}\right)\ {\rm in}\ \Omega^{\bullet}(S)/d\Omega^{\bullet}(S).
\end{multline}

Let $\psi\left(TX,\nabla^{TX}\right)$ be the Mathai-Quillen current on $TX$, which is defined as in \cite[Definition 6.7]{BG}. Let $h(\nabla^F,g^F)$ be defined in \cite[(1.34)]{BL}. By \cite[Section 7.1]{BG}, we see that $h(\nabla^F, g^F)(\nabla f)^*\psi(TX,\nabla^{TX})$ is smooth on $M$.

Consider the ${\bf Z}_{2}$-graded vector bundle $TX|_{\bf B}=TX^s|_{\bf B}\oplus TX^u|_{\bf B}$ over the manifold ${\bf B}$, where $TX^s$, $TX^u$ are defined in analogue with $T\widetilde{X}^s$, $T\widetilde{X}^u$. To avoid any ambiguity, let us just say that $TX^s$ is the even part of $TX|_{\bf B}$, and $TX^u|_{\bf B}$ is the corresponding odd part. We define the form ${}^0I(TX|_{\bf B},\nabla^{TX|_{\bf B}})$ as in \cite[(4.64)-(4.65)]{BG}. Let ${}^0I(TX|_{\bf B})$ be the corresponding cohomology class.

Now we state the main theorem of this paper.

\begin{thm}\label{main}
The following identity holds in $\Omega^{\bullet}(S)/d\Omega^{\bullet}(S)$,
\begin{multline}
\mathcal{T}_{L^2,h}\left(T^H\widetilde{M},g^{T\widetilde{X}},\nabla^{\widetilde{F}},g^{\widetilde{F}}\right)-T_{L^2,h}\left(\widetilde{A}^{C^{\bullet}({W}^u,\widetilde{F})'},g^{C^{\bullet}({W}^u,\widetilde{F})}\right)\\
+\widetilde{h}_{L^2}\left(\nabla^{H^{\bullet}_{(2)}(\widetilde{X},\widetilde{F}|_{\widetilde{X}})},g^{H^{\bullet}_{(2)}(\widetilde{X},\widetilde{F}|_{\widetilde{X}})}_{C^{\bullet}({W}^u,\widetilde{F})},g^{H^{\bullet}_{(2)}(\widetilde{X},\widetilde{F}|_{\widetilde{X}})}_{L_2}\right)\\=-\int_{X}h(\nabla^{F},g^{F})(\nabla f)^*\psi (TX,\nabla^{TX})
+\sum_{x\in B}(-1)^{{\rm ind}(x)}{\rm dim}F\ {}^{0}I(T_x X|_{{\bf B}}).
\end{multline}
\end{thm}

\begin{proof}
By the discussions in \cite[Chapter 7]{BG}, as in \cite[Chapter 9]{BG}, we will assume that $g^{TX}=h^{TX}$ and the metrics $g^{TX}$, $g^F$ satisfy the assumptions in \cite[Section 9.1]{BG}.

Following \cite[(9.6)]{BG}, set
$$\widetilde{\chi}'^{-}(F)=\sum_{x\in B}{\rm rk}(F)(-1)^{{\rm ind}(x)}{\rm dim}\left(T_x X^u|_{\bf B}\right),$$
$$\widetilde{\chi}'^{+}(F)=\sum_{x\in B}{\rm rk}(F)(-1)^{{\rm ind}(x)}{\rm dim}\left(T_x X^s|_{\bf B}\right),$$
$$\chi'(F)=\sum_{i=0}^{{\rm dim}X}(-1)^i i {\rm dim}\left(H^i (X,F|_X)\right),$$

$${\rm Tr}_{s}[f]=\sum_{x\in B}(-1)^{{\rm ind}(x)}{\rm rk}(F)f(x).$$

Recall that the Bismut-Lott torsion form is defined by 
\begin{multline}
\mathcal{T}_{h}(T^H M,g^{TX},\nabla^{F},g^{F})=-\int_{0}^{+\infty}\left[h^{\wedge}\left(A',g_{t}^{\Omega^{\bullet}(X,F|_{X})}\right)-{1\over 2}\chi'(F)h'(0)\right.\\
-\left.\left({n\over 4}\chi(F)-{1\over 2}\chi'(F)\right)h'\left(i\sqrt{t}/2\right)\right]{{dt}\over{t}}.
\end{multline}
By \cite[Theorem 0.1]{BG}, we only need to prove in $\Omega^{\bullet}(S)/d\Omega^{\bullet}(S)$,
\begin{multline}\label{n3}
\mathcal{T}_{L^2,h}\left(T^H\widetilde{M},g^{T\widetilde{X}},\nabla^{\widetilde{F}},g^{\widetilde{F}}\right)-T_{L^2,h}\left(\widetilde{A}^{C^{\bullet}({W}^u,\widetilde{F})'},g^{C^{\bullet}({W}^u,\widetilde{F})}\right)\\
+\widetilde{h}_{L^2}\left(\nabla^{H^{\bullet}_{(2)}(\widetilde{X},\widetilde{F}|_{\widetilde{X}})},g^{H^{\bullet}_{(2)}(\widetilde{X},\widetilde{F}|_{\widetilde{X}})}_{C^{\bullet}({W}^u,\widetilde{F})},g^{H^{\bullet}_{(2)}(\widetilde{X},\widetilde{F}|_{\widetilde{X}})}_{L_2}\right)\\
=\mathcal{T}_{h}\left(T^H {M},g^{T{X}},\nabla^{{F}},g^{{F}}\right)-T_{h}\left(A^{C^{\bullet}({W}^u,{F})'},g^{C^{\bullet}({W}^u,{F})}\right)\\
+\widetilde{h}\left(\nabla^{H^{\bullet}({X},{F}|_{{X}})},g^{H^{\bullet}({X},{F}|_{{X}})}_{C^{\bullet}({W}^u,{F})},g^{H^{\bullet}({X},{F}|_{{X}})}_{L_2}\right).
\end{multline}

For $T\geq 0$, let $g^{F}_{T}$ be the metric on $F$,
\begin{align}
g^{F}_{T}=e^{-2Tf}g^{F}.
\end{align}
Let $g^{\widetilde{F}}_T$ be the metric on $\widetilde{F}$,
\begin{align}\nonumber
g^{\widetilde{F}}_{T}=e^{-2T\widetilde{f}}g^{\widetilde{F}}.
\end{align}

By the anomaly formulas of Bismut-Lott analytic torsion form and $L^2$-analytic torsion form, (\ref{n3}) is equivalent to the claim that for any $T\geq 0$, in $\Omega^{\bullet}(S)/d\Omega^{\bullet}(S)$, we have
\begin{multline}
\mathcal{T}_{L^2,h}\left(T^H\widetilde{M},g^{T\widetilde{X}},\nabla^{\widetilde{F}},g_T^{\widetilde{F}}\right)-\mathcal{T}\left(T^H {M},g^{T{X}},\nabla^{{F}},g_T^{{F}}\right)\\-\left(T_{L^2,h}\left(\widetilde{A}^{C^{\bullet}({W}^u,\widetilde{F})'},g^{C^{\bullet}({W}^u,\widetilde{F})}\right)-T_h\left(A^{C^{\bullet}({W}^u,{F})'},g^{C^{\bullet}({W}^u,{F})}\right)\right)\\
+\left(\widetilde{h}_{L^2}\left(\nabla^{H_{(2)}^{\bullet}(\widetilde{X},\widetilde{F}|_{\widetilde{X}})},g^{H_{(2)}^{\bullet}(\widetilde{X},\widetilde{F}|_{\widetilde{X}})}_{C^{\bullet}({W}^u,\widetilde{F})},g^{H_{(2)}^{\bullet}(\widetilde{X},\widetilde{F}|_{\widetilde{X}})}_{L_2,T}\right)\right.\\-\left.\widetilde{h}\left(\nabla^{H^{\bullet}({X},{F}|_{{X}})},g^{H^{\bullet}({X},{F}|_{{X}})}_{C^{\bullet}({W}^u,{F})},g^{H^{\bullet}({X},{F}|_{{X}})}_{L_2,T}\right)\right)
=0.
\end{multline}
So we need only to prove in $\Omega^{\bullet}(S)/d\Omega^{\bullet}(S)$, 
\begin{multline}\label{4.6}
\lim_{T\to +\infty}\left\{\mathcal{T}_{L^2,h}\left(T^H\widetilde{M},g^{T\widetilde{X}},\nabla^{\widetilde{F}},g_T^{\widetilde{F}}\right)-\mathcal{T}\left(T^H {M},g^{T{X}},\nabla^{{F}},g_T^{{F}}\right)\right.
\\-\left(T_{L^2,h}\left(\widetilde{A}^{C^{\bullet}({W}^u,\widetilde{F})'},g^{C^{\bullet}({W}^u,\widetilde{F})}\right)-T_h\left(A^{C^{\bullet}({W}^u,{F})'},g^{C^{\bullet}({W}^u,{F})}\right)\right)\\
+\widetilde{h}_{L^2}\left(\nabla^{H_{(2)}^{\bullet}(\widetilde{X},\widetilde{F}|_{\widetilde{X}})},g^{H_{(2)}^{\bullet}(\widetilde{X},\widetilde{F}|_{\widetilde{X}})}_{C^{\bullet}({W}^u,\widetilde{F})},g^{H_{(2)}^{\bullet}(\widetilde{X},\widetilde{F}|_{\widetilde{X}})}_{L_2,T}\right)\\
\left.-\widetilde{h}\left(\nabla^{H^{\bullet}({X},{F}|_{{X}})},g^{H^{\bullet}({X},{F}|_{{X}})}_{C^{\bullet}({W}^u,{F})},g^{H^{\bullet}({X},{F}|_{{X}})}_{L_2,T}\right)\right\}
=0.
\end{multline}

We define $\widetilde{C}_{t,T}$ and $\widetilde{D}_{t,T}$ as in (\ref{edp2.15}) with respect to $g^{\widetilde{F}}_T$, 
then we need the following $L^2$-extension of \cite[Theorems 9.7 and 9.8]{BG}.

\begin{thm}\label{t3}
There exists $\delta\in (0,1/2]$ such that if $\varepsilon$, $A$ are such that $0<\varepsilon<A<+\infty$, there exists $C>0$ such that if $t\in [\varepsilon,A]$, $T\geq 1$, then
\begin{align}
\left|{\rm Tr}_{\Gamma,s}\left[Nh'(\widetilde{D}_{t,T})\right]-\widetilde{\chi}'^{-}(F)\right|\leq {C\over{T^\delta}}.
\end{align}
\end{thm}

\begin{thm}\label{t4.7}
The following identity holds in $\Omega^{\bullet}(S)/d\Omega^{\bullet}(S)$,
\begin{multline}
\lim_{T\to +\infty}\left\{\int_{1}^{+\infty}\left({\rm Tr}_{\Gamma,s}\left[Nh'\left(\widetilde{D}_{t,T}\right)\right]-\chi'(F)\right){{dt}\over{2t}}\right.\\
-\widetilde{h}^*_{L^2}\left(\nabla^{H^{\bullet}_{(2)}(\widetilde{X},\widetilde{F}|_{\widetilde{X}})},g^{H^{\bullet}_{(2)}(\widetilde{X},\widetilde{F}|_{\widetilde{X}})}_{L_2,0},g^{H^{\bullet}_{(2)}(\widetilde{X},\widetilde{F}|_{\widetilde{X}})}_{L_2,T}\right)-{\rm Tr}_s[f]T\\
\left.-{1\over 4}\left(\widetilde{\chi}'^{+}(F)-\widetilde{\chi}'^{-}(F)\right)\log(T)\right\}\\
=\int_{0}^{1}\left({\rm Tr}_{\Gamma,s}\left[N^{C^{\bullet}({W}^u,\widetilde{F})}h'\left(B_t^{C^{\bullet}({W}^u,\widetilde{F})}\right)\right]-\widetilde{\chi}'^{-}(F)\right){{dt}\over{2t}}\\
+\int_{1}^{+\infty}\left({\rm Tr}_{\Gamma,s}\left[N^{C^{\bullet}({W}^u,\widetilde{F})}h'\left(B_t^{C^{\bullet}({W}^u,\widetilde{F})}\right)\right]-{\chi}'(\widetilde{F})\right){{dt}\over{2t}}\\
-\widetilde{h}^*_{L^2}\left(\nabla^{H^{\bullet}_{(2)}(\widetilde{X},\widetilde{F}|_{\widetilde{X}})},g^{H^{\bullet}_{(2)}(\widetilde{X},\widetilde{F}|_{\widetilde{X}})}_{L_2,0},g^{H^{\bullet}_{(2)}(\widetilde{X},\widetilde{F}|_{\widetilde{X}})}_{C^{\bullet}({W}^u,\widetilde{F})}\right)+{1\over 4}\left(\widetilde{\chi}'^{-}(F)-\widetilde{\chi}'^{+}(F)\right)\log(\pi),
\end{multline}
where $*$ means that the factors $2i\pi$ are omitted (cf. \cite[Chapter 9]{BG}).
\end{thm}

The proofs of Theorems \ref{t3} and \ref{t4.7} will be given in next sections.

Let $\alpha>0$ be a positive constant which will be chosen later. Let $f:\mathbb{\bf R}\to [0,1]$ be a smooth even function  such that
\begin{align}
f(t)=1\ \ {\rm if}\ \ |t|\leq {\alpha\over 2},\ \ f(t)=0\ \ {\rm if}\ \ |t|\geq \alpha.
\end{align}
Set 
\begin{align}
g(t)=1-f(t).
\end{align}
\begin{definition}
For $t\in (0,1]$, $a\in\mathbb{\bf C}$, set
\begin{align}\nonumber
F_t(a)=\int_{-\infty}^{+\infty}\exp(iu\sqrt{2}a)\exp\left(-{u^2\over 2}\right)f(ut){du\over{\sqrt{2\pi}}},
\end{align}
\begin{align}
G_t(a)=\int_{-\infty}^{+\infty}\exp(iu\sqrt{2}a)\exp\left(-{u^2\over 2}\right)g(ut){du\over{\sqrt{2\pi}}}.\end{align}
\end{definition}
Then 
\begin{align}\label{a16}
\exp(-a^2)=F_t(a)+G_t(a).
\end{align}
The functions $F_t(a)$, $G_t(a)$ are even holomorphic functions. Therefore there exist holomorphic functions $\widetilde{F}_t(a)$, $\widetilde{G}_t(a)$ such that
\begin{align}\nonumber
F_t (a)=\widetilde{F}_t(a^2),
\end{align}
\begin{align}\label{a17}
G_t(a)=\widetilde{G}_t(a^2).
\end{align}
From (\ref{a16}), (\ref{a17}), we get
\begin{align}\label{a18}
\exp(-a)=\widetilde{F}_t(a)+\widetilde{G}_t(a).
\end{align}
The restrictions of $F_t$, $G_t$ to ${\bf R}$ lie in the Schwartz space $\mathcal{S}({\bf R})$. Therefore the restrictions of $\widetilde{F}_t$, $\widetilde{G}_t$ to ${\bf R}$ also lie in $\mathcal{S}({\bf R})$.

From (\ref{a18}), we deduce that 
\begin{align}\label{su4.18}
\exp\left(-\widetilde{C}^2_{t,T}\right)=\widetilde{F}_{\sqrt{t}}\left(\widetilde{C}^2_{t,T}\right)+\widetilde{G}_{\sqrt{t}}\left(\widetilde{C}^2_{t,T}\right).
\end{align}

We define $C_{t,T}$ and $D_{t,T}$ for $X\to M\xrightarrow{\pi} S$ in analogue with $\widetilde{C}_{t,T}$, $\widetilde{D}_{t,T}$ for $\widetilde{X}\to \widetilde{M}\xrightarrow{\widetilde{\pi}} S$ as in (\ref{edp2.15}) with respect to $g^F_T$.

We have the following analogue theorem of \cite[Theorem 11.3]{B} and \cite[Theorem 5.3]{Ma}.
\begin{theorem}\label{det4.5}
There exist $c>0$, $C>0$ such that for $t\in (0,1]$, $T\geq 0$, then
\begin{align}\label{de4.19}
\left|{\rm Tr}_{\Gamma,s}\left[N\widetilde{G}_{\sqrt{t}}\left(\widetilde{C}^2_{t,T}\right)\right]\right|\leq c\exp\left(-{C\over{t}}\right),\ \  \left|{\rm Tr}_{s}\left[N\widetilde{G}_{\sqrt{t}}\left({C}^2_{t,T}\right)\right]\right|\leq c\exp\left(-{C\over{t}}\right).
\end{align}
\end{theorem}
\begin{proof}
Set
\begin{align}
H_t(a)=\int_{-\infty}^{+\infty}\exp(iu\sqrt{2}a)\exp\left(-{u^2\over {2t^2}}\right)g(u){{du}\over{t\sqrt{2\pi}}}.
\end{align}
Then 
\begin{align}\label{a22}
G_{\sqrt{t}}(a)=H_{\sqrt{t}}\left({a\over {\sqrt{t}}}\right).
\end{align}
By \cite[eq. (13.23)]{BLe}, we find that for any $c\in {\bf R}_+$, $m\in {\bf N}$, there exist $c_m>0$, $C_m>0$ such that
\begin{align}\label{a23}
\sup_{a\in {\bf C}, |{\rm Im}(a)|\leq c}|a|^m |H_{\sqrt{t}} (a)|\leq c_m \exp\left(-{C_m\over{t}}\right).
\end{align}
Again there is a holomorphic function $\widetilde{H}_t(a)$ such that
\begin{align}\label{a24}
H_{\sqrt{t}}(a)=\widetilde{H}_{\sqrt{t}}(a^2)
\end{align}
and so  by (\ref{a22}), (\ref{a24})
\begin{align}\label{a25}
\widetilde{G}_{\sqrt{t}}(a)=\widetilde{H}_{\sqrt{t}}\left({a\over t}\right). 
\end{align}

Let $\Delta'$ be the contour in ${\bf C}$ defined by
\begin{multline}
\Delta'=\{x+iy| +\infty>x\geq -1,y=1\}\cup \{x+iy|x=-1, 1\geq y\geq -1\}\\
\cup \{x+iy|-1\leq x<+\infty,y=-1\}. 
\end{multline}
From (\ref{a23}), we deduce that
\begin{align}\label{ab27}
\sup_{a\in \Delta'}|a|^m \left|\widetilde{H}_{\sqrt{t}}(a)\right|\leq c\exp \left(-{C\over t}\right).
\end{align}
Let $\widetilde{H}_{t,p}(a)$ be a holomorphic function such that
\begin{align}
\left\{\begin{matrix}\lim_{a\to +\infty}\widetilde{H}_{\sqrt{t},p}(a)=0,\\{{\widetilde{H}^{p-1}_{\sqrt{t},p}(a)}\over{(p-1)!}}=\widetilde{H}_{\sqrt{t},p}(a).\end{matrix}\right.
\end{align}
By (\ref{a23}), we see that for any $m\in{\bf N}$,
\begin{align}
\sup_{a\in\Delta'}|a|^m \left|\widetilde{H}_{\sqrt{t},p}(a)\right|\leq c\exp\left(-{C\over{t}}\right).
\end{align}

By (\ref{a25}) and (\ref{ad5.63}), we have
\begin{align}\label{de4.29}
{\rm Tr}_{\Gamma,s}\left[N\widetilde{G}_{\sqrt{t}}\left(\widetilde{C}^2_{t,T}\right)\right]=\psi^{-1}_t{\rm Tr}_{\Gamma,s}\left[N\widetilde{H}_{\sqrt{t}}\left(\widetilde{A}^2_T\right)\right],
\end{align}
where $\widetilde{A}_T$ is defined by (\ref{ad5.13}).

Clearly
\begin{align}
\widetilde{H}_{\sqrt{t}} \left(\widetilde{A}^2_T\right)={1\over{2i\pi}}\int_{\Delta'}{{\widetilde{H}_{\sqrt{t}}(\lambda)}\over{\lambda-\widetilde{A}^2_T}}d\lambda.
\end{align}
Equivalently
\begin{align}\label{de4.31}
\widetilde{H}_{\sqrt{t}} \left(\widetilde{A}^2_T\right)={1\over{2i\pi}}\int_{\Delta'}{{\widetilde{H}_{\sqrt{t},p}(\lambda)}\over{\left(\lambda-\widetilde{A}^2_T\right)^p}}d\lambda.
\end{align}
Using (\ref{ab27}), (\ref{de4.31}) and proceeding as in \cite[Chapter 9]{B}, we find that for $t\in (0,1]$, $T\geq 1$,
\begin{align}\label{de4.32}
\left|{\rm Tr}_{\Gamma,s}\left[N\widetilde{H}_{\sqrt{t}}\left(\widetilde{A}^2_T\right)\right]\right|\leq c\exp\left(-{C\over{t}}\right). 
\end{align}
Using (\ref{de4.29})-(\ref{de4.32}), we get the first formula in (\ref{de4.19}). Similarly, one can get the second inequality in (\ref{de4.19}). 
\end{proof}

\begin{proposition}
The following formula holds, 
\begin{align}\label{4.7}
\lim_{T\to +\infty}\int_{0}^{1}\left({\rm Tr}_{s}\left[Nh'(D_{t,T})\right]-{\rm Tr}_{\Gamma,s}\left[Nh'(\widetilde{D}_{t,T})\right]\right){{dt}\over{2t}}=0.
\end{align}
\end{proposition}

\begin{proof}

Take $\alpha>0$ small enough so that for any $x\in \widetilde{M}$, $\widetilde{F}_t(\widetilde{D}^2_{t,T})(x,x)$ only depends on the behavior of $\widetilde{D}_{t,T}$ in a sufficiently small neighborhood of $x\in\widetilde{M}$. Then one has
\begin{align}\label{de4.34}
{\rm Tr}_{s}\left[N\widetilde{F}_{\sqrt{t}}\left({D}^2_{t,T}\right)\right]=
{\rm Tr}_{\Gamma,s}\left[N\widetilde{F}_{\sqrt{t}}\left(\widetilde{D}^2_{t,T}\right)\right].
\end{align}
Then by (\ref{de4.34}) and Theorem \ref{det4.5}, there exist $c>0$, $C>0$ such that for $t\in (0,1]$, $T\geq 0$,
\begin{align}\label{a35}
\left|{\rm Tr}_{s}\left[Nh'(D_{t,T})\right]-{\rm Tr}_{\Gamma,s}\left[Nh'(\widetilde{D}_{t,T})\right]\right|\leq c\exp\left(-{C\over{t}}\right).
\end{align}

By (\ref{a35}), for any $\varepsilon'>0$, there exists $\varepsilon>0$ such that
\begin{align}\label{a36}
\left|\int_0^{\varepsilon}\left({\rm Tr}_{s}\left[Nh'(D_{t,T})\right]-{\rm Tr}_{\Gamma,s}\left[Nh'(\widetilde{D}_{t,T})\right]\right){{dt}\over{2t}}\right|<\varepsilon' .
\end{align}
By Theorem \ref{t3} and \cite[Theorem 9.7]{BG}, we have
\begin{align}\label{a37}
\left|\int_{\varepsilon}^{1}\left({\rm Tr}_{s}\left[Nh'(D_{t,T})\right]-{\rm Tr}_{\Gamma,s}\left[Nh'(\widetilde{D}_{t,T})\right]\right){{dt}\over{2t}}\right|\leq {C\over {T^\delta}}{1\over \varepsilon}.
\end{align}
Then by (\ref{a36}) and (\ref{a37}), we get (\ref{4.7}).
\end{proof}

\begin{rem}
Using the methods in Section 6, as the estimates (\ref{ap6.31}) and (\ref{ap6.33}), one can give another proof of (\ref{4.7}).
\end{rem}

By \cite[(9.39) and (9.40)]{BG}, we have
\begin{multline}\label{4.8}
\varphi\left(-\int_0^1\left({\rm Tr}_{s}\left[Nh'(D_{t,T})\right]-{1\over 2}n\chi(F)h'(0)\right)\right.{{dt}\over{2t}}\\
\left.-\int_{1}^{+\infty}\left({\rm Tr}_{s}\left[Nh'(D_{t,T})\right]-\chi'(F)h'(0)\right){{dt}\over{2t}}\right)\\
=\mathcal{T}_{h}\left(T^H {M},g^{T{X}},\nabla^{{F}},g_T^{{F}}\right)\\
-\left({n\over 4}\chi(F)-{1\over 2}\chi'(F)\right)\left[\int_0^1 \left(h'\left(i\sqrt{t}/2\right)-h'(0)\right){{dt}\over t}+\int_1^{+\infty}h'\left(i\sqrt{t}/2\right){{dt}\over t}\right].
\end{multline}
For the $L^2$-analytic torsion form, by (\ref{ad2.16}) and (\ref{ad2.22}), we have 
\begin{multline}
\mathcal{T}_{L^2,h}\left(T^H {\widetilde{M}},g^{T\widetilde{X}},\nabla^{\widetilde{F}},g_T^{\widetilde{F}}\right)=-\varphi\left(\int_{0}^{1} \left[{\rm Tr}_{\Gamma,s}\left[{1\over 2}Nh'(\widetilde{D}_{t,T})\right]-{n\over 4}\chi(F)h'(0)\right]{{dt}\over t}\right)\\
-\int_{0}^{1}\left[{n\over 4}\chi(F)h'(0)-{1\over 2}\chi'(\widetilde{F})h'(0)-\left({n\over 4}\chi(F)-{1\over 2}{\chi}'(\widetilde{F})\right)h'\left({i\sqrt{t}\over 2}\right)\right]{{dt}\over t}\\
-\varphi\left(\int_{1}^{+\infty}\left[{\rm Tr}_{\Gamma,s}\left[{1\over 2}Nh'(\widetilde{D}_{t,T})\right]-{1\over 2}\chi'(\widetilde{F})h'(0)\right]{dt\over t}\right)\\+\int_{1}^{+\infty}\left({n\over 4}\chi(F)-{1\over 2}\chi'(\widetilde{F})\right)h'\left({i\sqrt{t}\over 2}\right){dt\over t}.
\end{multline}

Then we also have
\begin{multline}\label{4.9}
\varphi\left(-\int_0^1\left({\rm Tr}_{s}\left[Nh'(\widetilde{D}_{t,T})\right]-{1\over 2}n\chi(F)h'(0)\right)\right.{{dt}\over{2t}}\\
\left.-\int_{1}^{+\infty}\left({\rm Tr}_{s}\left[Nh'(\widetilde{D}_{t,T})\right]-\chi'(\widetilde{F})h'(0)\right){{dt}\over{2t}}\right)\\
=\mathcal{T}_{L^2,h}\left(T^H {\widetilde{M}},g^{T\widetilde{X}},\nabla^{\widetilde{F}},g_T^{\widetilde{F}}\right)\\
-\left({n\over 4}\chi(F)-{1\over 2}\chi'(\widetilde{F})\right)\left[\int_0^1 \left(h'\left(i\sqrt{t}/2\right)-h'(0)\right){{dt}\over t}+\int_1^{+\infty}h'\left(i\sqrt{t}/2\right){{dt}\over t}\right].
\end{multline}

By \cite[(9.71)]{BG},
\begin{align}
\int_0^1 \left(h'\left(i\sqrt{t}/2\right)-h'(0)\right){{dt}\over t}+\int_1^{+\infty}h'\left(i\sqrt{t}/2\right){{dt}\over t}=\Gamma'(1)+2(\log(2)-1).
\end{align}
Then by (\ref{4.8}) and (\ref{4.9}), we have
\begin{multline}\label{4.11}
\mathcal{T}_{L^2,h}\left(T^H {\widetilde{M}},g^{T\widetilde{X}},\nabla^{\widetilde{F}},g_T^{\widetilde{F}}\right)-\mathcal{T}_{h}\left(T^H {M},g^{T{X}},\nabla^{{F}},g_T^{{F}}\right)=\\
\varphi\left(-\int_0^1\left({\rm Tr}_{s}\left[Nh'(\widetilde{D}_{t,T})\right]-{1\over 2}n\chi(F)h'(0)\right)\right.{{dt}\over{2t}}\\
\left.-\int_{1}^{+\infty}\left({\rm Tr}_{s}\left[Nh'(\widetilde{D}_{t,T})\right]-\chi'(\widetilde{F})h'(0)\right){{dt}\over{2t}}\right)\\
-\varphi\left(-\int_0^1\left({\rm Tr}_{s}\left[Nh'(D_{t,T})\right]-{1\over 2}n\chi(F)h'(0)\right)\right.{{dt}\over{2t}}\\
\left.-\int_{1}^{+\infty}\left({\rm Tr}_{s}\left[Nh'(D_{t,T})\right]-\chi'(F)h'(0)\right){{dt}\over{2t}}\right)\\-{1\over 2}\left(\chi'(\widetilde{F}-\chi'(F))\right)\left(\Gamma'(1)+2(\log(2)-1)\right).
\end{multline}

By (\ref{4.7}), (\ref{4.11}), \cite[Theorem 9.8]{BG} and Theorem \ref{t4.7}, we have 
\begin{multline}
\lim_{T\to +\infty}\left(\mathcal{T}_{L^2,h}\left(T^H {\widetilde{M}},g^{T\widetilde{X}},\nabla^{\widetilde{F}},g_T^{\widetilde{F}}\right)-\mathcal{T}_{h}\left(T^H {M},g^{T{X}},\nabla^{{F}},g_T^{{F}}\right)\right.\\
+\widetilde{h}_{L^2}\left(\nabla^{H^{\bullet}_{(2)}(\widetilde{X},\widetilde{F}|_{\widetilde{X}})},g^{H^{\bullet}_{(2)}(\widetilde{X},\widetilde{F}|_{\widetilde{X}})}_{L_2,0},g^{H^{\bullet}_{(2)}(\widetilde{X},\widetilde{F}|_{\widetilde{X}})}_{L_2,T}\right)\\-\left.\widetilde{h}\left(\nabla^{H^{\bullet}({X},{F}|_{{X}})},g^{H^{\bullet}({X},{F}|_{{X}})}_{L_2,0},g^{H^{\bullet}({X},{F}|_{{X}})}_{L_2,T}\right)\right)\\
+\widetilde{h}\left(\nabla^{H^{\bullet}({X},{F}|_{{X}})},g^{H^{\bullet}({X},{F}|_{{X}})}_{L_2,0},g^{H^{\bullet}({X},{F}|_{{X}})}_{C^{\bullet}({W}^u,{F})}\right)\\
-\widetilde{h}_{L^2}\left(\nabla^{H^{\bullet}_{(2)}(\widetilde{X},\widetilde{F}|_{\widetilde{X}})},g^{H^{\bullet}_{(2)}(\widetilde{X},\widetilde{F}|_{\widetilde{X}})}_{L_2,0},g^{H^{\bullet}_{(2)}(\widetilde{X},\widetilde{F}|_{\widetilde{X}})}_{C^{\bullet}({W}^u,\widetilde{F})}\right)\\
=\lim_{T\to +\infty}\left(-\varphi \int_{1}^{+\infty}\left({\rm Tr}_{s}\left[Nh'(\widetilde{D}_{t,T})\right]-\chi'(\widetilde{F})\right){{dt}\over{2t}}\right.\\+\varphi \int_{1}^{+\infty}\left({\rm Tr}_{s}\left[Nh'({D}_{t,T})\right]-\chi'({F})\right){{dt}\over{2t}}\\
+\widetilde{h}_{L^2}\left(\nabla^{H^{\bullet}_{(2)}(\widetilde{X},\widetilde{F}|_{\widetilde{X}})},g^{H^{\bullet}_{(2)}(\widetilde{X},\widetilde{F}|_{\widetilde{X}})}_{L_2,0},g^{H^{\bullet}_{(2)}(\widetilde{X},\widetilde{F}|_{\widetilde{X}})}_{L_2,T}\right)\\-\left.\widetilde{h}\left(\nabla^{H^{\bullet}({X},{F}|_{{X}})},g^{H^{\bullet}({X},{F}|_{{X}})}_{L_2,0},g^{H^{\bullet}({X},{F}|_{{X}})}_{L_2,T}\right)
\right)\\
+\widetilde{h}\left(\nabla^{H^{\bullet}({X},{F}|_{{X}})},g^{H^{\bullet}({X},{F}|_{{X}})}_{L_2,0},g^{H^{\bullet}({X},{F}|_{{X}})}_{C^{\bullet}({W}^u,{F})}\right)\\
-\widetilde{h}_{L^2}\left(\nabla^{H^{\bullet}_{(2)}(\widetilde{X},\widetilde{F}|_{\widetilde{X}})},g^{H^{\bullet}_{(2)}(\widetilde{X},\widetilde{F}|_{\widetilde{X}})}_{L_2,0},g^{H^{\bullet}_{(2)}(\widetilde{X},\widetilde{F}|_{\widetilde{X}})}_{C^{\bullet}({W}^u,\widetilde{F})}\right)\\
-{1\over 2}(\chi'(\widetilde{F})-\chi'(F))(\Gamma'(1)+2(\log(2)-1))\\
=
-\varphi\int_{0}^{1}\left({\rm Tr}_{\Gamma,s}\left[N^{C^{\bullet}({W}^u,\widetilde{F})}h'\left(B_t^{C^{\bullet}({W}^u,\widetilde{F})}\right)\right]-\widetilde{\chi}'^{-}(F)\right){{dt}\over{2t}}\\
-\varphi\int_{1}^{+\infty}\left({\rm Tr}_{\Gamma,s}\left[N^{C^{\bullet}({W}^u,\widetilde{F})}h'\left(B_t^{C^{\bullet}({W}^u,\widetilde{F})}\right)\right]-{\chi}'(\widetilde{F})\right){{dt}\over{2t}}\\
+\varphi\int_{0}^{1}\left({\rm Tr}_{s}\left[N^{C^{\bullet}({W}^u,{F})}h'\left(B_t^{C^{\bullet}({W}^u,{F})}\right)\right]-\widetilde{\chi}'^{-}(F)\right){{dt}\over{2t}}\\
+\varphi\int_{1}^{+\infty}\left({\rm Tr}_{s}\left[N^{C^{\bullet}({W}^u,{F})}h'\left(B_t^{C^{\bullet}({W}^u,{F})}\right)\right]-{\chi}'({F})\right){{dt}\over{2t}}\\
-{1\over 2}(\chi'(\widetilde{F})-\chi'(F))(\Gamma'(1)+2(\log(2)-1))\\
=T_{L^2,h}\left(\widetilde{A}^{C^{\bullet}({W}^u,\widetilde{F})'},g^{C^{\bullet}({W}^u,\widetilde{F})}\right)-T_h\left(A^{C^{\bullet}({W}^u,{F})'},g^{C^{\bullet}({W}^u,{F})}\right).
\end{multline}
So the formula (\ref{4.6}) holds, then we get Theorem \ref{main}.

\end{proof}

\section{A proof of Theorem \ref{t4.7}}
\setcounter{equation}{0}

In this section, we will give a proof of Theorem \ref{t4.7}. In \cite{BG}, there are two proofs of \cite[Theorem 9.8]{BG}. In this paper, we will use the method in \cite[Chapter 10]{BG}. Almost the $L^2$ analogue results in \cite[Chapter 10]{BG} hold, except that it should be modified for the convergence as $t\to +\infty$, such as \cite[Theorem 10.37]{BG}. Also for the finite dimension of $F^{[0,1]}_T$ and etc., here it should be of finite $\Gamma$-dimension.

\subsection
{The harmonic oscillator near $\widetilde{\bf {B}}$
}\label{s5.1}
This subsection is the $L^2$-set up of \cite[Section 10.1]{BG}.

Recall that $\widetilde{\bf B}$ is the set of fiberwise critical points of $\widetilde{f}$. Then we have a $\mathbb{Z}_2$-graded vector bundle $T\widetilde{X}|_{\widetilde{\bf B}}$ (cf. \ref{s3.2}). The metric $g^{T\widetilde{X}}$ induces a metric $g^{T\widetilde{X}|_{\widetilde{\bf B}}}=g^{T\widetilde{X}^s|_{\widetilde{\bf B}}}\oplus g^{T\widetilde{X}^u|_{\widetilde{\bf B}}}$. Recall that $T\widetilde{X}^s|_{\widetilde{\bf B}}$ is the even part of $T\widetilde{X}|_{\widetilde{\bf B}}$ and $T\widetilde{X}^u|_{\widetilde{\bf B}}$ is the odd part of $T\widetilde{X}|_{\widetilde{\bf B}}$. Let $\pi': T\widetilde{X}|_{\widetilde{\bf B}}\to \widetilde{\bf B}$ denote the $\mathbb{Z}_2$-graded real vector bundle.

\begin{definition} (\cite [Definition 4.1]{BG})
For $x\in \widetilde{\bf B}$, let $\widetilde{I}_x$ (resp. $\widetilde{I}^0_x$) be the vector space of smooth (resp. square integrable) sections of $\pi'^* \Lambda ((T\widetilde{X}|_{\widetilde{\bf B}})^*)$ along the fiber $T\widetilde{X}_{x}$. 
\end{definition}

Let $d^{T\widetilde{X}|_{\widetilde{\bf B}}}$ be the de Rham operator acting along the fibers $\widetilde{I}$.

The connection $\nabla^{T\widetilde{X}|_{\widetilde{\bf B}}}$ induces a horizontal vector bundle $T^H T\widetilde{X}|_{\widetilde{\bf B}}$ on the total space $T\widetilde{X}|_{\widetilde{\bf B}}$. Then one verifies easily that with the notation in (\ref{ad1}), if $U, V\in T\widetilde{\bf B}$, $Z\in T\widetilde{X}|_{\widetilde{\bf B}}$,
$$T^H_Z(U,V)=R^{T\widetilde{X}|_{\widetilde{\bf B}}}(U,V)Z.$$

Also, with the notation in \cite[Definition 3.2]{BL} and in (\ref{nn1}), if $U\in T\widetilde{\bf B}$, if $s$ is a smooth section of $\widetilde{I}$ on $\widetilde{\bf B}$,
$$\nabla^{\widetilde{I}}_{U}s=\nabla_{U}^{\Lambda((T\widetilde{X}|_{\widetilde{\bf B}})^*)}s.$$

Let $q:T\widetilde{X}|_{\widetilde{\bf B}}\to {\bf R}$ be the smooth function, such that if $Z=(Z_+, Z_-)\in T\widetilde{X}|_{\widetilde{\bf B}}=T\widetilde{X}^s|_{\widetilde{\bf B}}\oplus T\widetilde{X}^u|_{\widetilde{\bf B}}$,
$$q(Z)={1\over 2}\left(\left|Z_+\right|^2-\left|Z_-\right|^2\right),$$
where $Z_+\in T\widetilde{X}^s$, $Z_-\in T\widetilde{X}^u$. Then $q$ is a fiberwise Morse function, whose only critical point $0$ has index ${\rm dim}T\widetilde{X}^u|_{\widetilde{\bf B}}$.

Clearly, $\widetilde{F}|_{\widetilde{\bf B}}$ is equipped with a flat connection $\nabla^{\widetilde{F}|_{\widetilde{\bf B}}}$. We denote by $\nabla^{\widetilde{I}\widehat{\otimes}\widetilde{F}|_{\widetilde{\bf B}}}$ the connection on $\widetilde{I}\widehat{\otimes}\widetilde{F}|_{\widetilde{\bf B}}$ induced by $\nabla^{\widetilde{I}}$ and $\nabla^{\widetilde{F}|_{\widetilde{\bf B}}}$. The metric $g^{\widetilde{F}}$ induces a metric $g^{\widetilde{F}|_{\widetilde{\bf B}}}$ on $\widetilde{F}|_{\widetilde{\bf B}}$. Let $\mathcal{F}$ be the restriction of $\widetilde{f}$ to $\widetilde{\bf B}$. We choose $\varepsilon>0$ as in \cite[Section 10.1]{BG}. Then by \cite[(9.3)]{BG}, using the notation in \cite[(4.6)]{BG}, if $x\in \widetilde{\bf B}$, $Z\in (T\widetilde{X}|_{\widetilde{\bf B}})_{x}$, $|Z|\leq \varepsilon$, set
\begin{align}\label{s5.1}
\widetilde{f}(Z)=\mathcal{F}(x)+q(Z).
\end{align}

For $T\in\mathbb{\bf R}$, the metric $g^{\widetilde{F}}_{T}$ induces the metric $g^{\widetilde{F}|_{\widetilde{\bf B}}}_{T}$, which is given by 
\begin{align}
g^{\widetilde{F}|_{\widetilde{\bf B}}}_{T}=e^{-2T\mathcal{F}}g^{\widetilde{F}|_{\widetilde{\bf  B}}}.
\end{align}

The exterior differentiation operator $\mathcal{C}^{\widetilde{I}\widehat{\otimes}\widetilde{F}|'_{\widetilde{\bf B}}}$, acting on $\Omega^{\bullet}(\widetilde{M}|_{\widetilde{\bf B}},\widetilde{F}|_{\widetilde{\bf B}})$, has degree $1$ and satisfies $(\mathcal{C}^{\widetilde{I}\widehat{\otimes}\widetilde{F}|'_{\widetilde{\bf B}}})^2=0$.
Thus $\mathcal{C}^{\widetilde{I}\widehat{\otimes}\widetilde{F}|'_{\widetilde{\bf B}}}$ defines a flat superconnection of total degree $1$ on $\widetilde{I}\widehat{\otimes}\widetilde{F}|_{\widetilde{\bf B}}$. As in \cite[(4.5)]{BG}, we have the identity
\begin{align}
\mathcal{C}^{\widetilde{I}\widehat{\otimes}\widetilde{F}|'_{\widetilde{\bf B}}}=d^{T\widetilde{X}|_{\widetilde{\bf B}}}+\nabla^{\widetilde{I}\widehat{\otimes}\widetilde{F}|_{\widetilde{\bf B}}}+i_{R^{T\widetilde{X}|_{\widetilde{B}}}Z}.
\end{align}
Given $T\in\mathbb{\bf R}$, let $\mathcal{C}_T^{\widetilde{I}\widehat{\otimes}\widetilde{F}|''_{\widetilde{\bf B}}}$ be the adjoint flat superconnection with respect to the metrics $g^{T\widetilde{X}|_{\widetilde{\bf B}}}$, $g_{T}^{\widetilde{F}|_{\widetilde{\bf B}}}$. Then we have (cf. \cite[(10.5)]{BG})
\begin{multline}
\mathcal{C}_T^{\widetilde{I}\widehat{\otimes}\widetilde{F}|''_{\widetilde{\bf B}}}
=d^{T\widetilde{X}|_{\widetilde{\bf B}}*}+2Ti_{Z_{+}-Z_{-}}+\nabla^{\widetilde{I}\widehat{\otimes}\widetilde{F}|_{\widetilde{\bf B}}}+\omega\left(\widetilde{F}|_{\widetilde{\bf B}},g^{\widetilde{F}|_{\widetilde{\bf B}}}\right)-2Td\mathcal{F}-R^{T\widetilde{X}|_{\widetilde{\bf B}}}Z\wedge,
\end{multline}
where 
$$\omega\left(\widetilde{F}|_{\widetilde{\bf B}},g^{\widetilde{F}|_{\widetilde{\bf B}}}\right)=\left(\nabla^{\widetilde{F}|_{\widetilde{\bf B}}}\right)^*-\nabla^{\widetilde{F}|_{\widetilde{\bf B}}}.$$

Set
\begin{align}\label{ad5.5}
\mathcal{C}^{\widetilde{I}\widehat{\otimes}\widetilde{F}|_{\widetilde{\bf B}}}_{T}={1\over 2}\left(\mathcal{C}_T^{\widetilde{I}\widehat{\otimes}\widetilde{F}|''_{\widetilde{\bf B}}}+\mathcal{C}^{\widetilde{I}\widehat{\otimes}\widetilde{F}|'_{\widetilde{\bf B}}}\right),\ \ \mathcal{D}_{T}^{\widetilde{I}\widehat{\otimes}\widetilde{F}|_{\widetilde{\bf B}}}={1\over 2}\left(\mathcal{C}_T^{\widetilde{I}\widehat{\otimes}\widetilde{F}|''_{\widetilde{\bf B}}}-\mathcal{C}^{\widetilde{I}\widehat{\otimes}\widetilde{F}|'_{\widetilde{\bf B}}}\right).
\end{align}

Set 
\begin{align}\nonumber
\overline{\mathcal{C}}^{\widetilde{I}\widehat{\otimes}\widetilde{F}|'_{\widetilde{\bf B}}}=e^{-T\widetilde{f}}
\mathcal{C}^{\widetilde{I}\widehat{\otimes}\widetilde{F}|'_{\widetilde{\bf B}}}e^{T\widetilde{f}},\ \ \overline{\mathcal{C}}^{\widetilde{I}\widehat{\otimes}\widetilde{F}|''_{\widetilde{\bf B}}}=e^{T\widetilde{f}}
\mathcal{C}_{0}^{\widetilde{I}\widehat{\otimes}\widetilde{F}|''_{\widetilde{\bf B}}}e^{-T\widetilde{f}},
\end{align}
\begin{align}\label{ad5.6}
\overline{\mathcal{C}}^{\widetilde{I}\widehat{\otimes}\widetilde{F}|_{\widetilde{\bf B}}}_{T}={1\over 2}\left(\overline{\mathcal{C}}_T^{\widetilde{I}\widehat{\otimes}\widetilde{F}|''_{\widetilde{\bf B}}}+\overline{\mathcal{C}}^{\widetilde{I}\widehat{\otimes}\widetilde{F}|'_{\widetilde{\bf B}}}\right),\ \ \overline{\mathcal{D}}_{T}^{\widetilde{I}\widehat{\otimes}\widetilde{F}|_{\widetilde{\bf B}}}={1\over 2}\left(\overline{\mathcal{C}}_T^{\widetilde{I}\widehat{\otimes}\widetilde{F}|''_{\widetilde{\bf B}}}-\overline{\mathcal{C}}^{\widetilde{I}\widehat{\otimes}\widetilde{F}|'_{\widetilde{\bf B}}}\right).
\end{align}
By (\ref{s5.1}) and (\ref{ad5.6}), we get
\begin{multline}\nonumber
\overline{\mathcal{C}}^{\widetilde{I}\widehat{\otimes}\widetilde{F}|_{\widetilde{\bf B}}}_{T}={1\over 2}\left(d^{T\widetilde{X}|_{\widetilde{\bf B}}*}+d^{T\widetilde{X}|_{\widetilde{\bf B}}}\right)+{T\over 2}(i_{Z_+-Z_-}+(Z_+-Z_-)\wedge)+\nabla^{\widetilde{I}\widehat{\otimes}\widetilde{F}|_{\widetilde{\bf B}}}\\
+{1\over 2}\omega\left(\widetilde{F}|_{\widetilde{\bf B}},g^{\widetilde{F}|_{\widetilde{B}}}\right)+{1\over 2}\left(i_{R^{T\widetilde{X}|_{\widetilde{\bf B}}}Z}-R^{T\widetilde{X}|_{\widetilde{\bf B}}}Z\wedge\right),
\end{multline}
\begin{multline}\label{ad5.7}
\overline{\mathcal{D}}_{T}^{\widetilde{I}\widehat{\otimes}\widetilde{F}|_{\widetilde{\bf B}}}={1\over 2}\left(d^{T\widetilde{X}|_{\widetilde{\bf B}}*}-d^{T\widetilde{X}|_{\widetilde{\bf B}}}\right)+{T\over 2}(i_{Z_+-Z_-}-(Z_+-Z_-)\wedge)\\
+{1\over 2}\omega\left(\widetilde{F}|_{\widetilde{\bf B}},g^{\widetilde{F}|_{\widetilde{B}}}\right)-Td{\mathcal{F}}-{1\over 2}\left(i_{R^{T\widetilde{X}|_{\widetilde{\bf B}}}Z}+R^{T\widetilde{X}|_{\widetilde{\bf B}}}Z\wedge\right).
\end{multline}

Let $e_1,\dots, e_{n_{+}}$ be an orthonormal basis of $T\widetilde{X}^s|_{\widetilde{\bf B}}$, let $e_{n_{+}+1},\dots,e_n$ be an orthonormal basis of $T\widetilde{X}^u|_{\widetilde{\bf B}}$. 
Set (cf. \cite[3.33]{BG})
$$c(e_i)=e_i\wedge -i_{e_i},\ \ \widehat{c}(e_i)=e_i\wedge +i_{e_i},$$
where we identify $T\widetilde{X}|_{\widetilde{\bf B}}$ with its dual by the metric.

Then we have
\begin{multline}\label{27}
\overline{\mathcal{C}}_T^{\widetilde{I}\widehat{\otimes}\widetilde{F}|_{\widetilde{\bf B}},2}=-{1\over 4}\left(\nabla_{e_i}+\left\langle R^{T\widetilde{X}|_{\widetilde{\bf B}}}Z,e_i\right\rangle\right)+{1\over 4}\left\langle e_i, R^{T\widetilde{X}|_{\widetilde{\bf B}}}e_j\right\rangle\widehat{c}(e_i)\widehat{c}(e_j)\\
-{1\over 4}\omega^2\left(\widetilde{F}|_{\widetilde{\bf B}},g^{\widetilde{F}|_{\widetilde{\bf B}}}\right)+{{T^2}\over 4}|Z|^2+{T\over 4}\left(\sum_{1\leq i\leq n_{+}}c(e_i)\widehat{c}(e_i)-\sum_{n_{+}+1\leq i\leq n_{+}+n_{-}}c(e_i)\widehat{c}(e_i)\right).
\end{multline}
In particular, by (\ref{27}), we get
\begin{multline}
\overline{\mathcal{C}}_T^{\widetilde{I}\widehat{\otimes}\widetilde{F}|_{\widetilde{\bf B}},2,(0)}=-{1\over 4}\Delta^{T\widetilde{X}|_{\widetilde{\bf B}}}+{{T^2}\over 4}|Z|^2
+{T\over 4}\left(\sum_{1\leq i\leq n_{+}}c(e_i)\widehat{c}(e_i)-\sum_{n_{+}+1\leq i\leq n_{+}+n_{-}}c(e_i)\widehat{c}(e_i)\right).
\end{multline}

Clearly,
\begin{align}
{\rm Sp}\left(\overline{\mathcal{C}}_T^{\widetilde{I}\widehat{\otimes}\widetilde{F}|_{\widetilde{\bf B}},2}\right)={\rm Sp}\left(\overline{\mathcal{C}}_T^{\widetilde{I}\widehat{\otimes}\widetilde{F}|_{\widetilde{\bf B}},2,(0)}\right).
\end{align}
The operator $\overline{\mathcal{C}}_T^{\widetilde{I}\widehat{\otimes}\widetilde{F}|_{\widetilde{\bf B}},2,(0)}$ is a harmonic oscillator. Take $\rho$ non zero in $\Lambda^{\rm max}(T^*\widetilde{X}^u|_{\widetilde{\bf B}})$. Let $\bar{\mathfrak{f}}_{T}$ be the one dimensional vector space spanned by $\exp(-T|Z|^2/2)\rho$. Then we have
\begin{align}\label{ad5.10}
{\rm ker}\overline{\mathcal{C}}_T^{\widetilde{I}\widehat{\otimes}\widetilde{F}|_{\widetilde{\bf B}},2,(0)}=\bar{\mathfrak{f}}_T\otimes \widetilde{F}|_{\widetilde{\bf B}}.
\end{align}

Let $\bar{\mathfrak{p}}_T$ be the orthogonal projection operator from $\widetilde{I}$ on $\bar{\mathfrak{f}}_T$. Then $\bar{\mathfrak{p}}_T$ extends to an endomorphism of $\widetilde{I}\widehat{\otimes}\widetilde{F}|_{\widetilde{\bf B}}$. Note that when acting on $\widetilde{I}\widehat{\otimes}\widetilde{F}|_{\widetilde{\bf B}}$, $\bar{\mathfrak{p}}_T$ is of the form $\bar{\mathfrak{p}}_T\otimes 1$.

For $T\in\mathbb{\bf R}^*_{+}$, using \cite[(4.12)]{BG}, we get
\begin{align}
{\rm Sp}\left(\overline{\mathcal{C}}_T^{\widetilde{I}\widehat{\otimes}\widetilde{F}|_{\widetilde{\bf B}},2,(0)}\right)={T\over 2}\mathbb{\bf N}.
\end{align}

\subsection{The eigenbundles associated to small eigenvalues}
This subsection is the $L^2$-set up of \cite[Section 10.2]{BG}

For $T\geq 0$, recall that $d_T^{\widetilde{X},*}$ is the adjoint of $d^{\widetilde{X}}$ with respect to the metrics $g^{T\widetilde{X}}$, $g_{T}^{\widetilde{F}}$. Let $\nabla_T^{\Omega_{(2)}^{\bullet}(\widetilde{X},\widetilde{F}|_{\widetilde{X}}),*}$ be the corresponding adjoint connection to $\nabla^{\Omega_{(2)}^{\bullet}(\widetilde{X},\widetilde{F}|_{\widetilde{X}})}$ and let $\widetilde{A}''_T$ be the superconnection on $\Omega_{(2)}^{\bullet}(\widetilde{X},\widetilde{F}|_{\widetilde{X}})$ adjoint to $\widetilde{A}'$ with respect to $g^{T\widetilde{X}}$, $g_T^{\widetilde{F}}$. Then
\begin{align}\nonumber
\widetilde{A}'=d^{\widetilde{X}}+\nabla^{\Omega_{(2)}^{\bullet}(\widetilde{X},\widetilde{F}|_{\widetilde{X}})}+i_{T^H}
\end{align}
and
\begin{align}\label{ad5.12}
\widetilde{A}_T''=d_T^{\widetilde{X},*}+\nabla_T^{\Omega_{(2)}^{\bullet}(\widetilde{X},\widetilde{F}|_{\widetilde{X}}),*}-T^H\wedge.
\end{align}

Set 
\begin{align}\label{ad5.13}
\widetilde{A}_T={1\over 2}\left(\widetilde{A}''_T+\widetilde{A}'\right),\ \ \ \widetilde{B}_T={1\over 2}\left(\widetilde{A}''_T-\widetilde{A}'\right).
\end{align}
Clearly,
\begin{align}
\widetilde{A}^2_T={1\over 4}\left[\widetilde{A}''_T,\widetilde{A}'\right].
\end{align}
In the sequel, we will also use the notations,
\begin{align}\label{ad5.15}
{\overline{A}}_T=e^{-T\widetilde{f}}\widetilde{A}_T e^{T\widetilde{f}},\ \ \ {\overline{B}}_T=e^{-T\widetilde{f}}\widetilde{B}_T e^{T\widetilde{f}}.
\end{align}

If $H\in\Lambda^{\bullet}(T^*S)\widehat{\otimes}{\rm End}(\Omega_{(2)}^{\bullet}(\widetilde{X},\widetilde{F}|_{\widetilde{X}}))$, let $H^{(0)}$ be the component of $H$ in ${\rm End}(\Omega_{(2)}^{\bullet}(\widetilde{X},\widetilde{F}|_{\widetilde{X}}))$. Clearly,
\begin{align}
{\rm Sp}\left(\widetilde{A}^2_T\right)={\rm Sp}\left({\overline{A}}^2_T\right)={\rm Sp}\left(\widetilde{A}^{2,(0)}_{T}\right)={\rm Sp}\left({\overline{A}}^{2,(0)}_{T}\right).
\end{align}

\begin{definition}
If $s\in S$, let $\widetilde{F}^{[0,1]}_{T,s}$ (resp. ${\overline{F}}^{[0,1]}_{T,s}$) be the eigenspaces of $\widetilde{A}^{2,(0)}_{T}$ (resp. ${\overline{A}}^{2,(0)}_{T}$) associated to eigenvalues $\lambda\in [0,1]$, let $\widetilde{P}^{[0,1]}_{T,s}$ (resp. ${\overline{P}}^{[0,1]}_{T,s}$) be the orthogonal projection operator form $(\Omega_{(2)}^{\bullet}(\widetilde{X}_s,\widetilde{F}|_{\widetilde{X}_s}),g_T^{\Omega_{(2)}^{\bullet}(\widetilde{X}_s,\widetilde{F}|_{\widetilde{X}_s})})$ on $\widetilde{F}^{[0,1]}_{T,s}$ (resp. from $(\Omega_{(2)}^{\bullet}(\widetilde{X}_s,\widetilde{F}|_{\widetilde{X}_s}),g^{\Omega_{(2)}^{\bullet}(\widetilde{X}_s,\widetilde{F}|_{\widetilde{X}_s})})$ on ${\overline{F}}^{[0,1]}_{T,s}$).
\end{definition}

Clearly, we have the obvious orthogonal splittings,
\begin{align}
\widetilde{F}^{[0,1]}_{T,s}=\bigoplus_{i=0}^{{\rm dim}\widetilde{X}}\widetilde{F}^{[0,1],i}_{T,s},\ \ {\overline{F}}^{[0,1]}_{T,s}=\bigoplus_{i=0}^{{\rm dim}\widetilde{X}}{\overline{F}}^{[0,1],i}_{T,s}.
\end{align}
Also,
\begin{align}
{\overline{P}}^{[0,1]}_{T,s}=e^{-T\widetilde{f}}\widetilde{P}^{[0,1]}_{T,s}e^{T\widetilde{f}}.
\end{align}
Let $M^i$ be the number of elements in $B^i$. Equivalently, $M^i$ is the number of critical points of $f$ in a given fiber $X$ whose index is equal to $i$.
\begin{theorem}\label{r2th5.3}
There exists $T_0\geq 0$ such that for $T\geq T_0$,
\begin{align}\label{6.19}
{\rm Sp}\left(\widetilde{A}^{2,(0)}_{T}\right)\subset \left[0,{1\over 4}\right]\cup [4,\infty),
\ {\rm dim}_{\Gamma}\left(\widetilde{F}^{[0,1],i}_{T,s}\right)=M^i,\ \ 1\leq i\leq {\rm dim}\widetilde{X}.
\end{align}
\end{theorem}

\begin{proof}
For a given $s\in S$, this result was established in \cite[Proposition 5.2 and Theorem 5.5]{BFKM}. Since $S$ is compact, a trivial uniformity argument shows that we can find $T_0\in {\bf R}_{+}$ such that (\ref{6.19}) holds for any $s\in S$, $T\geq T_0$.
\end{proof}

By the above, it follows that the $\widetilde{F}^{[0,1]}_{T,s}$, ${\overline{F}}^{[0,1]}_{T,s}$ are the fibers of smooth $\mathbb{\bf Z}$-graded vector bundles $\widetilde{F}^{[0,1]}_{T}$, ${\overline{F}}^{[0,1]}_{T}$ on $S$, which are subbundles of $\Omega_{(2)}^{\bullet}(\widetilde{X},\widetilde{F}|_{\widetilde{X}})$.

Clearly,
\begin{align}
h'(x)=(1+2x^2)e^{x^2}.
\end{align}
Put
\begin{align}\label{ad5.21}
r(x)=(1-2x)e^{-x}.
\end{align}

Let $\delta$ be the unit circle in $\mathbb{\bf C}$. Let $\Gamma=\Gamma_+\cup \Gamma_{-}$ be the contour defined by
\begin{align}\nonumber
\Gamma_{+}=\left\{z=x+iy|x\geq 2, y=\pm 1\right\}\cup \left\{z=x+iy|x=2, -1\leq y\leq 1\right\},
\end{align}
\begin{align}
\Gamma_{-}=\left\{z=x+iy|x\leq -2, y=\pm 1\right\}\cup \left\{z=x+iy|x=-2, -1\leq y\leq 1\right\}.
\end{align}

\begin{definition}
For $t\in \mathbb{\bf R}^*_{+}$, $T\geq T_0$, put
\begin{align}
\widetilde{K}_{t,T}=\psi^{-1}_t{1\over 2i\pi}\int_{\Gamma}{{r(t\lambda)}\over{\lambda-\widetilde{A}^2_T}}d\lambda\psi_t,\ 
\widetilde{L}_{t,T}={1\over{2i\pi}}\int_{\sqrt{t}\delta}{{h'(\lambda)}\over{\lambda-{\widetilde{D}}_{t,T}}}d\lambda.
\end{align}
\end{definition}

Then as \cite[Proposition 10.4]{BG}, we have
\begin{proposition}
The following identity holds,
\begin{align}\label{6.24}
h'(\widetilde{D}_{t,T})=\widetilde{K}_{t,T}+\widetilde{L}_{t,T}.
\end{align}
\end{proposition}

By (\ref{6.24}),
\begin{align}
{\rm Tr}_{\Gamma,s}\left[Nh'\left(\widetilde{D}_{t,T}\right)\right]={\rm Tr}_{\Gamma,s}\left[N\widetilde{K}_{t,T}\right]+{\rm Tr}_{\Gamma,s}\left[N\widetilde{L}_{t,T}\right].
\end{align}
By the same proof of \cite[Theorem 10.5]{BG}, we have the following analogue theorem in $L^2$-case.
\begin{theorem}\label{t5.5}
There exist $C>0$, $c>0$, $\delta\in (0,1]$ such that for $t\geq 1$, $T\geq T_0$,
\begin{align}\label{a27}
\left|{\rm Tr}_{\Gamma,s}\left[N \widetilde{K}_{t,T}\right]\right|\leq {{Ce^{-ct}}\over{T^{\delta}}}.
\end{align}
\end{theorem}

\begin{proof}
Set 
\begin{align}\label{ad5.27}
\widetilde{M}_{t,T,a}=\psi^{-1}_t {1\over{2i\pi}} \int_{\Gamma}{{\exp(-ta\lambda)}\over{\lambda-\widetilde{A}^2_T}}d\lambda \psi_t.
\end{align}
Then by (\ref{ad5.21}),
\begin{align}\label{ad29}
\widetilde{K}_{t,T}=\left(1+2{\partial\over{\partial a}}\right)\widetilde{M}_{t,T,a}|_{a=1}.
\end{align}

Put 
\begin{align}\label{ad5.29}
\overline{M}_{t,T,a}=e^{-T\widetilde{f}}\widetilde{M}_{t,T,a}e^{T\widetilde{f}}.
\end{align}
By (\ref{ad5.15}) and (\ref{ad5.27}),
\begin{align}
\overline{M}_{t,T,a}=\psi^{-1}_t {1\over{2i\pi}} \int_{\Gamma}{{\exp(-ta\lambda)}\over{\lambda-\overline{A}^2_T}}d\lambda \psi_t.
\end{align}
By (\ref{ad5.29}),
\begin{align}\label{ad32}
{\rm Tr}_{\Gamma,s}\left[N\widetilde{M}_{t,T,a}\right]={\rm Tr}_{\Gamma,s}\left[N\overline{M}_{t,T,a}\right].
\end{align}

Clearly if $d^{\widetilde{X},*}$ is the adjoint of $d^{\widetilde{X}}$ with respect to $g^{\Omega^{\bullet}_{(2)}(\widetilde{X},\widetilde{F}|_{\widetilde{X}})}$, then
\begin{align}\label{ad5.32}
d^{\widetilde{X},*}_{T}=d^{\widetilde{X},*}+2Ti_{\nabla \widetilde{f}}.
\end{align}
Set $D^{\widetilde{X}}=d^{\widetilde{X}}+d^{\widetilde{X},*}$. Recall that $\nabla^{\Omega^{\bullet}_{(2)}(\widetilde{X},\widetilde{F}|_{\widetilde{X}}),*}=\nabla^{\Omega^{\bullet}_{(2)}(\widetilde{X},\widetilde{F}|_{\widetilde{X}}),*}_{0}$ is the connection adjoint to $\nabla^{\Omega^{\bullet}_{(2)}(\widetilde{X},\widetilde{F}|_{\widetilde{X}})}$ with respect to $g^{\Omega^{\bullet}_{(2)}(\widetilde{X},\widetilde{F}|_{\widetilde{X}})}$. Let $(d\widetilde{f})^H$ be the horizontal component of $d\widetilde{f}$. Then
\begin{align}\label{ad5.33}
\nabla^{\Omega^{\bullet}_{(2)}(\widetilde{X},\widetilde{F}|_{\widetilde{X}}),*}_{T}=\nabla^{\Omega^{\bullet}_{(2)}(\widetilde{X},\widetilde{F}|_{\widetilde{X}}),*}-2T\left(d\widetilde{f}\right)^H.
\end{align}
From (\ref{ad5.12}), (\ref{ad5.13}), (\ref{ad5.32}) and (\ref{ad5.33}), we get
\begin{align}
\widetilde{A}_T={1\over 2}\left(D^{\widetilde{X}}+2Ti_{\nabla\widetilde{f}}\right)+\nabla^{\Omega^{\bullet}_{(2)}(\widetilde{X},\widetilde{F}|_{\widetilde{X}}),u}-T\left(d\widetilde{f}\right)^H-{1\over 2}c\left(T^H\right).
\end{align}
If $\nabla^{\Omega^{\bullet}_{(2)}\left(\widetilde{X},\widetilde{F}|_{\widetilde{X}}\right),u}=\nabla^{\Omega^{\bullet}_{(2)}\left(\widetilde{X},\widetilde{F}|_{\widetilde{X}}\right),u}_{0}$, by (\ref{ad5.15}), we obtain
\begin{align}\label{ad5.35}
\overline{A}_{T}={1\over 2}\left(D^{\widetilde{X}}+T\widehat{c}(\nabla \widetilde{f})\right)+\nabla^{\Omega^{\bullet}_{(2)}\left(\widetilde{X},\widetilde{F}|_{\widetilde{X}}\right),u}-{1\over 2}c\left(T^H\right).
\end{align}
The term $T(d\widetilde{f})^H$ disappeared in (\ref{ad5.35}).

Observe that fiberwise, $\widetilde{\bf B}$ is the zero set of $\nabla \widetilde{f}$. Also $\widehat{c}(\nabla\widetilde{f})$ anticommutes with the principal $c(i\xi)$ of $D^{\widetilde{X}}$. By the simplifying assumptions in \cite[Section 9.1]{BG}, near $\widetilde{\bf B}$,
\begin{align}\label{ad37}
\overline{A}_{T}=\overline{\mathcal{C}}^{\widetilde{I}\widehat{\otimes}\widetilde{F}|_{\widetilde{\bf B}}}_{T}.
\end{align}
By (\ref{ad5.10}), the kernel of $\overline{\mathcal{C}}^{\widetilde{I}\widehat{\otimes}\widetilde{F}|_{\widetilde{\bf B}},(0)}_{T}$ can be identified with $\widetilde{F}|_{\widetilde{\bf B}}\otimes o^u|_{\widetilde{\bf B}}$. Using (\ref{ad5.7}) and the fact, we get the easy formula,
\begin{align}
\bar{\mathfrak{p}}_T\overline{A}_T\bar{\mathfrak{p}}_T=\nabla^{\widetilde{F}|_{\widetilde{\bf B}},u}. 
\end{align}

Observe that
\begin{align}
{\rm Sp}\left(\nabla^{\widetilde{F}|_{\widetilde{\bf B}},u,2}\right)=\{0\}. 
\end{align}
Therefore,
\begin{align}\label{ad40}
{1\over{2i\pi}}\int_{\Gamma}{{\exp(-ta\lambda)}\over{\lambda-\nabla^{\widetilde{F}|_{\widetilde{\bf B}},u,2}}}d\lambda=0.
\end{align}

By (\ref{ad29}), (\ref{ad32}), (\ref{ad5.35}) and (\ref{ad37})-(\ref{ad40}), as in the proof of \cite[Theorem 10.5]{BG}, one need a $L^2$-version of \cite[(9.149)]{B}. We show the $L^2$-case in the Appendix $B$. Then the proof of the theorem is completed.

\end{proof}

\subsection{The projectors $\widetilde{\mathbb{\bf P}}^{[0,1]}_{T}$}
This subsection is the $L^2$-set up of \cite[Sections 10.3 and 10.4]{BG}.

Observe that $\Lambda^{\bullet}(T^*S)$ acts on $\Lambda^{\bullet}(T^*S)\widehat{\otimes }\Omega_{(2)}^{\bullet}(\widetilde{X},\widetilde{F}|_{\widetilde{X}})$. Let $\delta\subset \mathbb{\bf C}$ be the circle of centre $0$ and radius $1/4$. Recall that for $T\geq T_0$,
\begin{align}
{\rm Sp}\left(\widetilde{A}^2_{T}\right)\cap \delta=\emptyset.
\end{align}

\begin{definition}
For $T\geq T_0$, put
\begin{align}\label{d42}
\widetilde{\mathbb{\bf P}}^{[0,1]}_{T}={1\over {2i\pi}}\int_{\delta}{{d\lambda}\over{\lambda-\widetilde{A}^2_{T}}}.
\end{align}
\end{definition}

Clearly, $\widetilde{\mathbb{\bf P}}^{[0,1]}_{T}\in\Lambda^{\bullet}(T^*S)\widehat{\otimes}{\rm End}(\Omega_{(2)}^{\bullet}(\widetilde{X},\widetilde{F}|_{\widetilde{X}}))$. Then we write,
\begin{align}
\widetilde{\bf P}^{[0,1]}_{T}=\sum_{j=1}^{{\rm dim}S}\widetilde{\bf P}^{[0,1],(j)}_{T},
\end{align}
where $\widetilde{\bf P}_{T}^{[0,1],j}\in\Lambda^{j}(T^*S)\widehat{\otimes}{\rm End}(\Omega^{\bullet}_{(2)}(\widetilde{X},\widetilde{F}|_{\widetilde{X}}))$.

 Set $\widetilde{\mathbb{\bf F}}^{[0,1]}_{T}={\rm Im}(\widetilde{\mathbb{\bf P}}^{[0,1]}_{T})$, then $\widetilde{\mathbb{\bf F}}^{[0,1]}_{T}$ is a ${\bf Z}_2$-graded subbundle of $\Lambda^{\bullet}(T^*S)\widehat{\otimes}\Omega^{\bullet}_{(2)}(\widetilde{X},\widetilde{F}|_{\widetilde{X}})$. Also we have
\begin{align}\label{6.30}
\widetilde{\bf P}^{[0,1],(0)}_{T}=\widetilde{P}^{[0,1]}_{T}.
\end{align}

In the sequel, the operator $*$ acts on $\Lambda^{\bullet}(T^*S)\widehat{\otimes}{\rm End}(\Omega^{\bullet}_{(2)}(\widetilde{X},\widetilde{F}|_{\widetilde{X}}))$ as in \cite[(1.8)]{BG}, with respect to the metric $g_T^{\Omega_{(2)}^{\bullet}(\widetilde{X},\widetilde{F}|_{\widetilde{X}})}$. We will often say that if $k$ is such that $k^*=k$, then it is self-adjoint.

If $k\in \Lambda^{\bullet}(T^*S)\widehat{\otimes}{\rm End}(\Omega^{\bullet}_{(2)}(\widetilde{X},\widetilde{F}|_{\widetilde{X}}))$, we can write $k$ in the form, 
\begin{align}\label{ad5.31}
k=\sum_{j=0}^{{\rm dim}S}k^{(j)},\ \ k^{(j)}\in\Lambda^{j}(T^*S)\widehat{\otimes}{\rm End}(\Omega^{\bullet}_{(2)}(\widetilde{X},\widetilde{F}|_{\widetilde{X}})).
\end{align}
Observe that $\Lambda^{\bullet}(T^*S)\widehat{\otimes}{\rm End}(\Omega^{\bullet}_{(2)}(\widetilde{X},\widetilde{F}|_{\widetilde{X}}))$ is a $\mathbb{Z}$-graded bundle of algebras. Namely, if $k$ is taken as in (\ref{ad5.31}), ${\rm deg}k=p$ if for any $j$,
\begin{align}\nonumber
k^{(j)}\in \Lambda^{(j)}(T^*S)\widehat{\otimes}{\rm Hom}\left(\Omega^{\bullet}_{(2)}(\widetilde{X},\widetilde{F}|_{\widetilde{X}}),\Omega^{\bullet+p-j}_{(2)}(\widetilde{X},\widetilde{F}|_{\widetilde{X}})\right).
\end{align}
Moreover $\Lambda^{\bullet}(T^*S)\widehat{\otimes}{\rm End}(\Omega^{\bullet}_{(2)}(\widetilde{X},\widetilde{F}|_{\widetilde{X}}))$ inherits a filtration $F$ from the filtration of $\Lambda^{\bullet}(T^*S)$. We will say that ${\rm deg}(k)\geq 0$ if it is the sum of elements of non-negative degree. Also we will write that 
\begin{align}\nonumber
{\rm deg}(k)\leq 2F(k)
\end{align}
if for any $j$, ${\rm deg}(k^j)\leq 2j$.

\begin{theorem}\label{thm5.8}
For $T\geq T_0$, $\widetilde{\bf P}^{[0,1]}_{T}$ is an even projection operator acting on $\Lambda^{\bullet}(T^*S)\widehat{\otimes}\Omega^{\bullet}_{(2)}(\widetilde{X},\widetilde{F}|_{\widetilde{X}})$, which commutes with the action of $\Lambda^{\bullet}(T^*S)$ and with $\widetilde{A}'$, $\widetilde{A}''_{T}$, and is such that
\begin{align}\label{b46}
\widetilde{\bf P}^{[0,1],(0)}_{T}=\widetilde{P}^{[0,1]}_{T}.
\end{align}
Also
\begin{align}\label{c47}
\widetilde{\bf P}^{[0,1],*}_{T}=\widetilde{\bf P}^{[0,1]}_{T},\ \ 0\leq {\rm deg}\left(\widetilde{\bf P}^{[0,1]}_{T}\right)\leq 2F\left(\widetilde{\bf P}^{[0,1]}_{T}\right).
\end{align}
The linear map $\alpha\in\Lambda(T^*S)\widehat{\otimes}\widetilde{F}_{T}^{[0,1]}\to \widetilde{\bf P}^{[0,1]}_{T}\alpha\in \widetilde{\bf F}^{[0,1]}_{T}$ is an isomorphism of $\mathbb{Z}_2$-graded filtered vector bundles.

\end{theorem}

\begin{proof}
By definition, $\widetilde{\bf P}^{0,1}_{T}$ commutes with $\widetilde{A}^2_T$. Since $\widetilde{A}'$, $\widetilde{A}''_T$ and the elements of $\Lambda^{\bullet}(T^*S)$ commute with $\widetilde{A}^2_T$, they also commute with $\widetilde{\bf P}^{[0,1]}_{T}$. We write $\widetilde{A}^2_T$ in the form,
\begin{align}
\widetilde{A}^2_T=\widetilde{A}^{2,(0)}_{T}+\widetilde{A}^{2,(>0)}_{T}.
\end{align}
Then if $\lambda\in\delta$,
\begin{align}\label{b47}
\left(\lambda-\widetilde{A}^2_T\right)^{-1}=\sum_{i=0}^{{\rm dim}X}\left(\lambda-\widetilde{A}_{T}^{2,(0)}\right)^{-1}\widetilde{A}^{2,(>0)}_{T}\cdots \widetilde{A}^{2,(>0)}_{T}\left(\lambda-\widetilde{A}_{T}^{2,(0)}\right)^{-1},
\end{align}
so that $\widetilde{A}^{2,(>0)}_{T}$ appears $i$ times in the right-hand side of (\ref{b47}). The term corresponding to $i=0$ is obviously equal to the projection operator $\widetilde{P}^{[0,1]}_{T}$, i.e. (\ref{b46}) holds. Also since $\widetilde{A}^{2,*}_{T}=\widetilde{A}^2_{T}$, the first identity in (\ref{c47}) also holds. Also $\widetilde{A}^{2,(0)}$ is of degree $0$. Using \cite[Proposition 10.7]{BG}, (\ref{d42}) and (\ref{b47}), we get the second identity in (\ref{c47}).

Obviously, $\widetilde{\bf F}^{[0,1]}_{T}$ inherits a filtration from the filtration of $\Lambda^{\bullet}(T^*S)\widehat{\otimes}\Omega^{\bullet}_{(2)}(\widetilde{X},\widetilde{F}|_{\widetilde{X}})$. If $\beta\in\widetilde{\bf F}^{[0,1],\geq j}_{T}$, then
\begin{align}
\beta=\beta^{(j)}+\beta^{(j+1)}+\cdots.
\end{align}
Since $\widetilde{\bf P}^{[0,1]}_{T}\beta=\beta$, using (\ref{b46}), we get 
\begin{align}
\widetilde{P}^{[0,1]}_{T}\beta^{(j)}=\beta^{(j)},
\end{align}
so that $\beta^{(j)}\in \Lambda^{j}(T^*S)\widehat{\otimes}\widetilde{F}^{[0,1]}_{T}$. This way, we defined an injective linear map ${\rm Gr}^j (\widetilde{\bf F}^{[0,1]}_{T})\to \Lambda^j (T^*S)\widehat{\otimes}\widetilde{F}^{[0,1]}_{T}$. An obvious inverse for this map is just $\alpha\in\Lambda^j (T^*S)\widehat{\otimes}\widetilde{F}^{[0,1]}_{T}\to \left[\widetilde{\bf P}^{[0,1]}_{T}\alpha\right]\in {\rm Gr}^j (\widetilde{\bf F}^{[0,1]}_{T})$. Therefore $\alpha\in\Lambda^{\bullet}(T^*S)\widehat{\otimes}\widetilde{F}^{[0,1]}_{T}\to \widetilde{\bf P}^{[0,1]}_{T}\alpha\in\widetilde{\bf F}^{[0,1]}_{T}$ is an isomorphism of ${\bf Z}_2$-graded filtered vector bundles. 
\end{proof}

\subsection{The projectors $\widetilde{\mathbb{\bf P}}^{[0,1]}_{t,T}$}
This subsection is the $L^2$-set up of \cite[Section 10.5]{BG}.

Since $T\widetilde{M}=T^{H}\widetilde{M}\oplus T\widetilde{X}$, we have a smooth identification $T\widetilde{M}\cong \pi^*TS\oplus T\widetilde{X}$. Therefore we have the identification,
\begin{align}
\Lambda(T^*\widetilde{M})\cong \pi^*\Lambda^{\bullet}(T^*S)\widehat{\otimes}\Lambda(T^*\widetilde{X}).
\end{align}
Let $N^{\Lambda^{\bullet}(T^*\widetilde{M})}$ be the operator defining the $\mathbb{\bf  Z}$-grading of $\Lambda(T^*\widetilde{M})$. Then $N^{\Lambda^{\bullet}(T^*\widetilde{M})}$ acts naturally on the space $\Lambda^{\bullet}(T^*S)\widehat{\otimes}\Omega_{(2)}^{\bullet}(\widetilde{X},\widetilde{F}|_{\widetilde{X}})$.

For $t>0$, $T\geq 0$, let $\widetilde{A}''_{t,T}$ be the adjoint superconnection to $\widetilde{A}'$ with respect to the metrics $g^{T\widetilde{X}}/t$, $g^{\widetilde{F}}$. Then we have
\begin{align}\label{ad5.53}
t^{-N^{\Lambda^{\bullet}(T^*\widetilde{M})/2}}\widetilde{A}'t^{N^{\Lambda^{\bullet}(T^*\widetilde{M})/2}}={1\over{\sqrt{t}}}\widetilde{A}',\ \ t^{-N^{\Lambda^{\bullet}(T^*\widetilde{M})/2}}\widetilde{A}''_Tt^{N^{\Lambda^{\bullet}(T^*\widetilde{M})/2}}={1\over{\sqrt{t}}}\widetilde{A}''_{t,T}.\end{align}
We also have
\begin{align}
{\rm Sp}\left(\widetilde{A}^2_{t,T}\right)\subset \left[0,{t\over 4}\right]\cup [2t,+\infty).
\end{align}

\begin{definition}
For $T\geq T_0$, $t>0$, put
\begin{align}
\widetilde{\mathbb{\bf P}}^{[0,1]}_{t,T}={{1\over{2i\pi}}\int_{t\delta}{{d\lambda}\over{\lambda-\widetilde{A}^2_{t,T}}}}.
\end{align}
\end{definition}

Set $\widetilde{\bf F}^{[0,1]}_{t,T}={\rm Im}\left(\widetilde{\bf P}^{[0,1]}_{t,T}\right)$, then $\widetilde{\bf P}^{[0,1]}_{t,T}$ is an even projection operator and $\widetilde{\bf F}^{[0,1]}_{t,T}$ is a finite $\Gamma$-dimensional subbundle of $\Lambda^{\bullet}(T^*S)\widehat{\otimes}\Omega^{\bullet}_{(2)}(\widetilde{X},\widetilde{F}|_{\widetilde{X}})$.

\begin{proposition}
The following identity holds,
\begin{align}\nonumber
\widetilde{\bf P}^{[0,1]}_{t,T}=t^{-N^{\Lambda^{\bullet}(T^*\widetilde{M})}/2}\widetilde{\bf P}^{[0,1]}_{T}
t^{-N^{\Lambda^{\bullet}(T^*\widetilde{M})}/2},
\end{align}
\begin{align}\label{ad5.56}
\widetilde{\bf F}^{[0,1]}_{t,T}=t^{-N^{\Lambda^{\bullet}(T^*\widetilde{M})}/2}\widetilde{\bf F}^{[0,1]}_{T}
.
\end{align}
\end{proposition}
\begin{proof}
This is a consequence of (\ref{ad5.53}).
\end{proof}

\begin{definition}
For $T\geq T_0$, put
\begin{align}
\widetilde{\mathbb{\bf P}}^{[0,1]}_{\infty,T}={1\over{2i\pi}}\int_{\delta}{{d\lambda}\over{\lambda-{1\over 4}\left[\widetilde{A}',d_T^{\widetilde{X},*}\right]}}.
\end{align}
\end{definition}
Then by the same argument of \cite[Proposition 10.14]{BG}, we have the following $L^2$-analogue of it.
\begin{proposition}
Given $\alpha\in\Lambda^{\bullet}(T^*S)\widehat{\otimes}\Omega_{(2)}^{\bullet}(\widetilde{X},\widetilde{F}|_{\widetilde{X}})$, as $t\to +\infty,$
\begin{align}
\widetilde{\mathbb{\bf P}}^{[0,1]}_{t,T}\alpha=\widetilde{\mathbb{\bf P}}^{[0,1]}_{\infty,T}\alpha+\mathcal{O}(1/t)\alpha.
\end{align}
\end{proposition}

For $t>0$, $T\geq 0$,
\begin{align}
\widetilde{C}_{t,T}=t^{N/2}\widetilde{A}_{t,T}t^{-N/2}.
\end{align}

\begin{definition}
Put
\begin{align}
{\widehat{\mathbb{\bf P}}}^{[0,1]}_{t,T}={1\over{2i\pi}}\int_{t\delta}{{d\lambda}\over{\lambda-\widetilde{C}^2_{t,T}}}.
\end{align}
\end{definition}

Set ${\widehat{\mathbb{\bf F}}}^{[0,1]}_{t,T}={\rm Im}\left({\widehat{\mathbb{\bf P}}}^{[0,1]}_{t,T}\right)$. Then ${\widehat{\mathbb{\bf P}}}^{[0,1]}_{t,T}$ is a projection operator and ${\widehat{\mathbb{\bf F}}}^{[0,1]}_{t,T}$ is a finite $\Gamma$-dimensional vector bundle. 

\begin{proposition}
The following identities hold,
\begin{align}\nonumber
\widehat{\bf P}^{[0,1]}_{t,T}=\psi_{t}^{-1}\widetilde{\bf P}^{[0,1]}_{T}\psi_t,
\end{align}
\begin{align}\nonumber
\widehat{\bf P}^{[0,1]}_{t,T}=t^{N/2}\widetilde{\bf P}^{[0,1]}_{t,T}t^{-N/2},
\end{align}
\begin{align}\label{ad5.61}
\widehat{\bf F}^{[0,1]}_{t,T}=\psi_{t}^{-1}\widetilde{\bf F}^{[0,1]}_{T}=t^{N/2}\widetilde{\bf F}^{[0,1]}_{t,T}. 
\end{align}
Also, as $t\to +\infty$,
\begin{align}\label{ad5.62}
\widehat{\bf P}^{[0,1]}_{t,T}=P_{T}^{[0,1]}s+\mathcal{O}\left(1/{\sqrt{t}}\right).
\end{align}
\end{proposition}
\begin{proof}
By (\ref{ad2.18}),
\begin{align}\label{ad5.63}
\widetilde{C}_{t,T}=\psi^{-1}_{t}\sqrt{t}\widetilde{A}_{T}\psi_t.
\end{align}
From (\ref{d42}), (\ref{ad5.56}) and (\ref{ad5.63}), we get (\ref{ad5.61}). By (\ref{b46}) and (\ref{ad5.61}), we get (\ref{ad5.62}). 
\end{proof}

\subsection{The maps $\widetilde{\mathbb{\bf P}}^{\infty}_{T}$}
This subsection is the $L^2$-set up of \cite[Section 10.6]{BG}.

\begin{definition}
Let $\widetilde{P}^{\infty}$ be the map
\begin{align}\nonumber
\alpha\in \Omega_{(2)}^{\bullet}(\widetilde{X},\widetilde{F}|_{\widetilde{X}})\to \widetilde{P}^{\infty}\alpha=\sum_{{x}\in\widetilde{B}} \widetilde{W}^u ({x})^* \int_{\overline{\widetilde{W}^u({x})}}\alpha\in C^{\bullet}(W^u,\widetilde{F}).
\end{align}
\end{definition}

Similarly as \cite[Definitions 5.2, 5.7]{BG}, we have $\widetilde{P}^{\infty}:\Omega_{(2)}^{\bullet}(\widetilde{X},\widetilde{F}|_{\widetilde{X}})\to C^{\bullet}(W^u,\widetilde{F})$ and $\widetilde{\mathbb{\bf P}}^{\infty}:\Omega_{(2)}^{\bullet}(\widetilde{M},\widetilde{F})\to \Lambda^{\bullet}(T^*S)\widehat{\otimes}C^{\bullet}(W^u,\widetilde{F})$. Then $\widetilde{P}^{\infty}$ and $\widetilde{\mathbb{\bf P}}^{\infty}$ are chain maps which preserves the $\mathbb{\bf Z}$-grading. Also $\widetilde{\mathbb{\bf P}}^{\infty}$ preserves the filtrations associated to $\Lambda^{\bullet}(T^*S)$.

\begin{definition}
For $T\geq T_0$, let $\widetilde{\mathbb{\bf P}}^{\infty}_{T}:\widetilde{\mathbb{\bf F}}^{[0,1]}_{T}\to \Lambda^{\bullet}(T^*S)\widehat{\otimes}C^{\bullet}(W^u,\widetilde{F})$, $\widetilde{P}^{\infty}_{T}:\widetilde{F}^{[0,1]}_T\to C^{\bullet}(W^u,\widetilde{F})$ be the restrictions of $\widetilde{\mathbb{\bf P}}^{\infty}$, $\widetilde{P}^{\infty}$ to $\widetilde{\mathbb{\bf F}}^{[0,1]}_T$, $\widetilde{F}^{[0,1]}_{T}$.
\end{definition}

Observe that $\widetilde{\mathbb{\bf P}}^{[0,1]}_{T}-\widetilde{P}^{[0,1]}_{T}$ contains only terms of positive degree in the Grassmann variables in $\Lambda^{\bullet}(T^*S)$. 

\begin{proposition}
The following $L^2$-analogue of \cite[(10.78)]{BG} holds,
\begin{align}\label{6.39}
\left(\widetilde{\bf P}^{\infty}_{T}\right)^{-1}=\widetilde{\bf P}^{[0,1]}_{T}\left(\widetilde{P}^{\infty}_{T}\right)^{-1}\left(\widetilde{\bf P}^{\infty}\widetilde{\bf P}^{[0,1]}_{T}\left(\widetilde{P}^{\infty}_{T}\right)^{-1}\right)^{-1},
\ \left(\widetilde{\bf P}^{\infty}_{T}\right)^{-1,(0)}=\left(\widetilde{P}^{\infty}_{T}\right)^{-1}.
\end{align}
\end{proposition}
\begin{proof}
First by (\ref{b46}), we get
\begin{align}\label{d61}
\left(\widetilde{\bf P}^{\infty}\widetilde{\bf P}^{[0,1]}_{T}\left(\widetilde{P}^{\infty}_{T}\right)^{-1}\right)^{(0)}=1.
\end{align}
By (\ref{d61}), $\widetilde{\bf P}^{\infty}\widetilde{\bf P}^{[0,1]}_{T}(P^{\infty}_{T})^{-1}$ is one to one. Recall that $(\widetilde{P}^{\infty}_{T})^{-1}$ and $\widetilde{\bf P}^{\infty}$ preserve the total degree. By Theorem \ref{thm5.8}, $\widetilde{\bf P}^{\infty}_{T}$ increases the total degree. Therefore $\widetilde{\bf P}^{\infty}\widetilde{\bf P}^{[0,1]}_{T}(\widetilde{P}^{\infty}_{T})^{-1}$ increases the total degree. Using the invertibility of $\widetilde{\bf P}^{\infty}\widetilde{\bf P}^{[0,1]}_{T}(\widetilde{P}^{\infty}_{T})^{-1}$, the first equation in (\ref{6.39}) follows. From (\ref{b46}) and (\ref{d61}), we get the second equation in (\ref{6.39}). 
\end{proof}

\subsection{The generalized metric $\mathfrak{g}_T^{C^{\bullet}(W^u,\widetilde{F})}$}
This subsection is the $L^2$-set up of \cite[Section 10.7]{BG}.

The map $(\widetilde{\mathbb{\bf P}}^{\infty}_{T})^{-1}$ identifies $\Lambda^{\bullet}(T^*S)\widehat{\otimes}C^{\bullet}(W^u, \widetilde{F})$ and $\widetilde{\mathbb{\bf F}}^{[0,1]}_{T}\subset \Lambda^{\bullet}(T^*S)\widehat{\otimes}\Omega_{(2)}^{\bullet}(\widetilde{X},\widetilde{F}|_{\widetilde{X}})$. In the sequel, we will consider $(\widetilde{\mathbb{\bf P}}^{\infty}_{T})^{-1}$ as a map from $\Lambda^{\bullet}(T^*S)\widehat{\otimes}C^{\bullet}(W^u,\widetilde{F})$ into $\Lambda^{\bullet}(T^*S)\widehat{\otimes}\Omega_{(2)}^{\bullet}(\widetilde{X},\widetilde{F}|_{\widetilde{X}})$.

Let $(\widetilde{\mathbb{\bf P}}^{\infty}_{T})^{-1,*}$ be the adjoint of $(\widetilde{\mathbb{\bf P}}^{\infty}_{T})^{-1}$ with respect to the metrics $g^{C^{\bullet}(W^u,\widetilde{F})}$, $g_T^{\Omega_{(2)}^{\bullet}(\widetilde{X},\widetilde{F}|_{\widetilde{X}})}$. Then $(\widetilde{\mathbb{\bf P}}^{\infty}_{T})^{-1,*}$ maps $\Lambda^{\bullet}(T^*S)\widehat{\otimes}\Omega_{(2)}^{\bullet}(\widetilde{X},\widetilde{F}|_{\widetilde{X}})$ into $\Lambda^{\bullet}(T^*S)\widehat{\otimes}C^{\bullet}(W^u,\widetilde{F})$.

\begin{definition}
For $T\geq T'_0$, put
\begin{align}
\mathfrak{g}_{T}^{C^{\bullet}(W^u,\widetilde{F})}=(\widetilde{\mathbb{\bf P}}^{\infty}_{T})^{-1,*}(\widetilde{\mathbb{\bf P}}^{\infty}_{T})^{-1},\ g_{T}^{C^{\bullet}(W^u,\widetilde{F})}=(\widetilde{P}^{\infty}_T)^{-1,*}(\widetilde{P}^{\infty}_{T})^{-1}.
\end{align}
\end{definition}

Observe that $\mathfrak{g}_{T}^{C^{\bullet}(W^u,\widetilde{F})}$ is a generalized metric on $C^{\bullet}(W^u,\widetilde{F})$ in the sense of \cite[Section 2.9]{BG}. Also $g_T^{C^{\bullet}(W^u,\widetilde{F})}$ is a standard metric on $C^{\bullet}(W^u,\widetilde{F})$, which is such that the $C^i(W^u,\widetilde{F})'s$ are orthogonal in $C^{\bullet}(W^u,\widetilde{F})$ with respect to $g_T^{C^{\bullet}(W^u,\widetilde{F})}$.

The following $L^2$-analogue of \cite[Theorem 10.21]{BG} is obviously holds.
\begin{theorem}
For $T\geq T'_0$,
\begin{align}
\mathfrak{g}_{T}^{C^{\bullet}(W^u,\widetilde{F}),(0)}=g_{T}^{C^{\bullet}(W^u,\widetilde{F})}.
\end{align}
Moreover,
\begin{align}\label{6.42}
\widetilde{\mathbb{\bf P}}^{\infty}_{T}\widetilde{A}'(\widetilde{\mathbb{\bf P}}^{\infty}_{T})^{-1}=\widetilde{A}^{C^{\bullet}(W^u,\widetilde{F})'},
\ \widetilde{\mathbb{\bf P}}^{\infty}_{T}\widetilde{A}''_T(\widetilde{\mathbb{\bf P}}^{\infty}_{T})^{-1}=\left(\mathfrak{g}_{T}^{C^{\bullet}(W^u,\widetilde{F})}\right)^{-1}\widetilde{A}^{C^{\bullet}(W^u,\widetilde{F})''}\mathfrak{g}_{T}^{C^{\bullet}(W^u,\widetilde{F})}.
\end{align}
\end{theorem}

\subsection{The maps $\widetilde{\mathbb{\bf P}}^{\infty}_{t,T}$ and the generalized metrics $\mathfrak{g}_{t,T}^{C^{\bullet}(W^u,\widetilde{F})}$}
This subsection is the $L^2$-set up of \cite[Section 10.8]{BG}.

Let $N^{\Lambda^{\bullet}(T^*S)\widehat{\otimes}C^{\bullet}(W^u,\widetilde{F})}$ be the number operator of $\Lambda^{\bullet}(T^*S)\widehat{\otimes}C^{\bullet}(W^u,\widetilde{F})$.

\begin{definition}
Given $T\geq T'_0$, $t>0$, let $\widetilde{\mathbb{\bf P}}^{\infty}_{t,T}:\widetilde{\mathbb{\bf F}}^{[0,1]}_{t,T}\to \Lambda^{\bullet}(T^*S)\widehat{\otimes}C^{\bullet}(W^u,\widetilde{F})$ be the restriction of $\widetilde{\mathbb{\bf P}}^{\infty}$ to $\widetilde{\mathbb{\bf F}}^{[0,1]}_{t,T}$.
\end{definition}

Recall that $\widetilde{\mathbb{\bf P}}^{\infty}_{\infty,T}:\widetilde{\mathbb{\bf F}}^{[0,1]}_{\infty,T}\to \Lambda^{\bullet}(T^*S)\widehat{\otimes}(W^u,\widetilde{F})$ is an isomorphism. Also we have the following $L^2$-analogue of \cite[Proposition 10.24]{BG}.
\begin{proposition}\label{p6.15}
The map $\widetilde{\mathbb{\bf P}}^{\infty}_{t,T}$ is invertible. Moreover,
\begin{align}\label{6.43}
\left(\widetilde{\mathbb{\bf P}}^{\infty}_{t,T}\right)^{-1}=t^{-N^{\Lambda^{\bullet}(T^*\widetilde{M})}/2}\left(\widetilde{\mathbb{\bf P}}^{\infty}_{T}\right)^{-1}t^{N^{\bullet}(T^*S)\widehat{\otimes}C^(W^u,\widetilde{F})/2}.
\end{align}
As $t\to +\infty$,
\begin{align}
\left(\widetilde{\mathbb{\bf P}}^{\infty}_{t,T}\right)^{-1}=\left(\widetilde{\mathbb{\bf P}}^{\infty}_{\infty,T}\right)^{-1}+\mathcal{O}(1/t).
\end{align}
\end{proposition}

Let $g_{t,T}^{\Omega_{(2)}^{\bullet}(\widetilde{X},\widetilde{F}|_{\widetilde{X}})}$ be the metric on $\Omega_{(2)}^{\bullet}(\widetilde{X},\widetilde{F}|_{\widetilde{X}})$ which is associated to the metrics $g^{T\widetilde{X}}/t$, $g^{\widetilde{F}}_{T}$ on $T\widetilde{X}$, $\widetilde{F}$. Let $(\widetilde{\mathbb{\bf P}}^{\infty}_{t,T})^{-1,*}$ be the adjoint of $(\widetilde{\mathbb{\bf P}}^{\infty}_{t,T})^{-1}$ with respect to the metrics $g^{\Omega_{(2)}^{\bullet}(\widetilde{X},\widetilde{F}|_{\widetilde{X}})}_{t,T}$, $g^{C^{\bullet}(W^u,\widetilde{F})}$.

Let $(\widetilde{P}^{\infty}_{T})^{-1,*}$ be the adjoint of $(\widetilde{P}^{\infty}_{T})^{-1}$ with respect to the metrics $g_{T}^{\Omega_{(2)}^{\bullet}(\widetilde{X},\widetilde{F}|_{\widetilde{X}})}$, $g^{C^{\bullet}(W^u,\widetilde{F})}$.

\begin{proposition}
The following identity holds,
\begin{align}\label{6.47}
t^{N/2}\left(\widetilde{\mathbb{\bf P}}_{t,T}\right)^{-1}t^{-N^{C^{\bullet}(W^u,\widetilde{F})}/2}=\psi^{-1}_{t}\left(\widetilde{\mathbb{\bf P}}^{\infty}_{T}\right)\psi_t.
\end{align}
\end{proposition}

\begin{proof}
This follows from Proposition \ref{p6.15}.
\end{proof}

Let $(\widetilde{\bf P}^{\infty}_{t,T})^{-1,*}_{0}$ be the adjoint of $(\widetilde{\bf P}^{\infty}_{t,T})^{-1}$ with respect to the metrics $g_{T}^{\Omega_{(2)}^{\bullet}(\widetilde{X},\widetilde{F}|_{\widetilde{X}})}$, $g^{C^{\bullet}(W^u,\widetilde{F})}$.

\begin{proposition}\label{p6.18}
There is a smooth section $\widetilde{J}$ of
$$\left(\Lambda^{\bullet}(T^*S)\widehat{\otimes}{\rm Hom}\left(C^{\bullet}(W^u,\widetilde{F}),\Omega_{(2)}^{\bullet}\left(\widetilde{X},\widetilde{F}|_{\widetilde{X}}\right)\right)\right)^{\rm even},$$
such that as $t\to +\infty$,
\begin{align}\nonumber
t^{N/2}\left(\widetilde{\mathbb{\bf P}}^{\infty}_{t,T}\right)^{-1}t^{-N^{C^{\bullet}(W^u,\widetilde{F})}/2}=\left(\widetilde{P}^{\infty}_{T}\right)^{-1}+{\widetilde{J}\over{\sqrt{t}}}+\mathcal{O}(1/t),
\end{align}
\begin{align}\label{6.48}
t^{-N^{C^{\bullet}(W^u,\widetilde{F})}/2}\left(\widetilde{\mathbb{\bf P}}^{\infty}_{t,T}\right)^{-1,*}_{0}t^{N/2}=\left(\widetilde{P}^{\infty}_{T}\right)^{-1,*}+{{\widetilde{J}_0^*}\over{t}}+\mathcal{O}(1/t).
\end{align}
\end{proposition}

\begin{proof}
Using (\ref{6.30}), (\ref{6.39}) and (\ref{6.47}), we get the first identity in (\ref{6.48}). By taking adjoints, we obtain the second identity. 
\end{proof}

\begin{definition}
Put 
\begin{align}
\mathfrak{g}_{t,T}^{C^{\bullet}(W^u,\widetilde{F})}=\left(\widetilde{\mathbb{\bf P}}^{\infty}_{t,T}\right)^{-1,*}\left(\widetilde{\mathbb{\bf P}}^{\infty}_{t,T}\right)^{-1}.
\end{align}
\end{definition}

Then $\mathfrak{g}_{t,T}^{C^{\bullet}(W^u,\widetilde{F})}$ is a generalized metric on $C^{\bullet}(W^u,\widetilde{F})$.

\begin{theorem}\label{t6.20}
There is a smooth section $\widetilde{H}$ of
$$\left(\Lambda^{\bullet}(T^*S)\widehat{\otimes}{\rm End}\left(C^{\bullet}(W^u,\widetilde{F})\right)\right)^{\rm even},$$
such that as $t\to +\infty$,
\begin{align}\label{6.50}
t^{-N^{C^{\bullet}(W^u,\widetilde{F})}/2+n/2}\mathfrak{g}_{t,T}^{C^{\bullet}(W^u,\widetilde{F})}t^{-N^{C^{\bullet}(W^u,\widetilde{F})}/2}=g_{T}^{C^{\bullet}(W^u,\widetilde{F})}+{\widetilde{H}\over{\sqrt{t}}}+\mathcal{O}(1/t),
\end{align}
\begin{multline}\label{6.51}
t^{N^{C^{\bullet}(W^u,\widetilde{F})}/2}\left[\mathfrak{g}_{t,T}^{C^{\bullet}(W^u,\widetilde{F})}\right]^{-1}{\partial\over{\partial t}}\left[\mathfrak{g}_{t,T}^{C^{\bullet}(W^u,\widetilde{F})}\right]t^{-N^{C^{\bullet}(W^u,\widetilde{F})}/2}\\
=\left(N^{C^{\bullet}(W^u,\widetilde{F})}-{n\over 2}\right){1\over t}+\mathcal{O}\left(1/t^{3/2}\right).
\end{multline}
\end{theorem}

\begin{proof}
First by definitions of $\left(\widetilde{\bf P}^{\infty}_{t,T}\right)^{-1,*}$ and $\left(\widetilde{\bf P}^{\infty}_{t,T}\right)_{0}^{-1,*}$, we have
\begin{align}
\left(\widetilde{\bf P}^{\infty}_{t,T}\right)^{-1,*}=\left(\widetilde{\bf P}^{\infty}_{t,T}\right)_{0}^{-1,*}t^{N-n/2},
\end{align}
then
\begin{multline}\label{6.53}
t^{-N^{C^{\bullet}(W^u,\widetilde{F})}/2}\mathfrak{g}_{t,T}^{C^{\bullet}(W^u,\widetilde{F})}t^{-N^{C^{\bullet}(W^u,\widetilde{F})}/2}\\
=\left(t^{-N^{C^{\bullet}(W^u,\widetilde{F})}/2-n/2}\left(\widetilde{\bf P}^{\infty}_{t,T}\right)^{-1,*}_{0}t^{N/2}\right)\left(t^{N/2}\left(\widetilde{\bf P}^{\infty}_{t,T}\right)^{-1}t^{-N^{C^{\bullet}(W^u,\widetilde{F})}/2}\right).
\end{multline}
By Proposition \ref{p6.18} and by (\ref{6.53}), we get (\ref{6.50}).

Also
\begin{multline}
\mathfrak{g}_{t,T}^{C^{\bullet}(W^u,\widetilde{F})}=
t^{N^{C^{\bullet}(W^u,\widetilde{F})}/2-n/2}\left[t^{n/2}t^{-N^{C^{\bullet}(W^u,\widetilde{F})}/2}\mathfrak{g}_{t,T}^{C^{\bullet}(W^u,\widetilde{F})}t^{-N^{C^{\bullet}(W^u,\widetilde{F})}/2}\right]t^{N^{C^{\bullet}(W^u,\widetilde{F})}/2}.
\end{multline}
Therefore,
\begin{multline}\label{6.55}
{\partial\over{\partial t}}\mathfrak{g}_{t,T}^{C^{\bullet}(W^u,\widetilde{F})}={{N^{C^{\bullet}(W^u,\widetilde{F})}}\over{2t}}\mathfrak{g}_{t,T}^{C^{\bullet}(W^u,\widetilde{F})}+\mathfrak{g}_{t,T}^{C^{\bullet}(W^u,\widetilde{F})}{{N^{C^{\bullet}(W^u,\widetilde{F})}}\over{2t}}\\
+t^{N^{C^{\bullet}(W^u,\widetilde{F})}/2-n/2}{\partial\over{\partial t}}\left(t^{n/2}t^{-N^{C^{\bullet}(W^u,\widetilde{F})}/2}\mathfrak{g}_{t,T}^{C^{\bullet}(W^u,\widetilde{F})}t^{-N^{C^{\bullet}(W^u,\widetilde{F})}/2}\right)t^{N^{C^{\bullet}(W^u,\widetilde{F})}/2}.
\end{multline}
By (\ref{6.43}), (\ref{6.48}) and (\ref{6.53}), one verifies easily that
$$t^{n/2}t^{-N^{C^{\bullet}(W^u,\widetilde{F})}/2}\mathfrak{g}_{t,T}^{C^{\bullet}(W^u,\widetilde{F})}t^{-N^{C^{\bullet}(W^u,\widetilde{F})}/2}$$
is a polynomial in $1/{\sqrt{t}}$. Therefore, as $t\to +\infty$,
\begin{align}\label{6.56}
{\partial\over{\partial t}}\left(t^{n/2}t^{-N^{C^{\bullet}(W^u,\widetilde{F})}/2}\mathfrak{g}_{t,T}^{C^{\bullet}(W^u,\widetilde{F})}t^{-N^{C^{\bullet}(W^u,\widetilde{F})}/2}\right)=\mathcal{O}\left(1/t^{3/2}\right).
\end{align}
From (\ref{6.50}), (\ref{6.55}) and (\ref{6.56}), we get (\ref{6.51}).
\end{proof}

\begin{proposition}
The following identity holds,
\begin{align}\label{6.57}
t^{N^{C^{\bullet}(W^u,\widetilde{F})}/2}\widetilde{\mathbb{\bf P}}^{\infty}_{t,T}\widetilde{A}'\left(\widetilde{\mathbb{\bf P}}^{\infty}_{t,T}\right)^{-1}t^{-N^{C^{\bullet}(W^u,\widetilde{F})}/2}=\sqrt{t}\widetilde{\partial}+\nabla^{C^{\bullet}(W^u,\widetilde{F})}.
\end{align}
Also as $t\to +\infty$,
\begin{align}\label{6.58}
t^{N^{C^{\bullet}(W^u,\widetilde{F})}/2}\widetilde{\mathbb{\bf P}}^{\infty}_{t,T}\widetilde{A}''_{t,T}\left(\widetilde{\mathbb{\bf P}}^{\infty}_{t,T}\right)^{-1}t^{-N^{C^{\bullet}(W^u,\widetilde{F})}/2}=\widetilde{P}^{\infty}_{T}\sqrt{t}d^{\widetilde{X},*}\left(\widetilde{P}^{\infty}_{T}\right)^{-1}+\mathcal{O}(1).
\end{align}
\end{proposition}

\begin{proof}
Identity (\ref{6.57}) follows from (\ref{6.42}). Also by (\ref{6.42}), we have
\begin{multline}
\widetilde{\bf P}^{\infty}_{t,T}\widetilde{A}''_{t,T}\left(\widetilde{\bf P}^{\infty}_{t,T}\right)^{-1}=\left[t^{-N^{C^{\bullet}(W^u,\widetilde{F})}+n/2}\left(\widetilde{\bf P}^{\infty}_{t,T}\right)^{-1,*}\left(\widetilde{\bf P}^{\infty}_{t,T}\right)^{-1}\right]^{-1}\\
\left(t^{-N^{C^{\bullet}(W^u,\widetilde{F})}}\widetilde{A}^{C^{\bullet}(W^u,\widetilde{F})''}t^{C^{\bullet}(W^u,\widetilde{F})}\right)t^{-N^{C^{\bullet}(W^u,\widetilde{F})}+n/2}\left(\widetilde{\bf P}^{\infty}_{t,T}\right)^{-1,*}\left(\widetilde{\bf P}^{\infty}_{t,T}\right)^{-1}.
\end{multline}
Then we have
\begin{multline}\label{6.60}
t^{N^{C^{\bullet}(W^u,\widetilde{F})}/2}\widetilde{\bf P}^{\infty}_{t,T}\widetilde{A}''_{t,T}\left(\widetilde{\bf P}^{\infty}_{t,T}\right)^{-1}t^{-N^{C^{\bullet}(W^u,\widetilde{F})}/2}=t^{N^{C^{\bullet}(W^u,\widetilde{F})}/2}\left[\mathfrak{g}_{t,T}^{C^{\bullet}(W^u,\widetilde{F})}\right]^{-1}\\
t^{N^{C^{\bullet}(W^u,\widetilde{F})}/2}\left(t^{-N^{C^{\bullet}(W^u,\widetilde{F})}/2}\widetilde{A}^{C^{\bullet}(W^u,\widetilde{F})''}t^{N^{C^{\bullet}(W^u,\widetilde{F})}/2}\right)t^{-N^{C^{\bullet}(W^u,\widetilde{F})}/2}\\
\mathfrak{g}_{t,T}^{C^{\bullet}(W^u,\widetilde{F})}t^{-N^{C^{\bullet}(W^u,\widetilde{F})}/2}.
\end{multline}
Using Theorem \ref{t6.20} and (\ref{6.60}), we find that as $t\to +\infty$,
\begin{multline}
t^{N^{C^{\bullet}(W^u,\widetilde{F})}/2}\widetilde{\bf P}^{\infty}_{t,T}\widetilde{A}''_{t,T}\left(\widetilde{\bf P}^{\infty}_{t,T}\right)^{-1}t^{-N^{C^{\bullet}(W^u,\widetilde{F})}/2}\\=\left(g_{T}^{C^{\bullet}(W^u,\widetilde{F})}\right)^{-1}\sqrt{t}\widetilde{\partial}^* g_T^{C^{\bullet}(W^u,\widetilde{F})}+\mathcal{O}(1),
\end{multline}
which is equivalent to (\ref{6.58}). The proof is completed.
\end{proof}

\subsection{The superconnection forms for $\widetilde{\mathbb{\bf F}}^{[0,1]}_{T}$}
This subsection is the $L^2$-set up of \cite[Sections 10.10, 10.13 and 10.14]{BG}.

Observe that 
\begin{align}
\widetilde{\mathbb{\bf P}}^{[0,1]}_{t,T}={1\over{2i\pi}}\int_{\sqrt{t}\delta}{{d\lambda}\over{\lambda-\widetilde{B}_{t,T}}}.
\end{align}
Using the holomorphic functional calculus, we find that
\begin{align}\nonumber
h(\widetilde{B}_{t,T})\widetilde{\mathbb{\bf P}}^{[0,1]}_{t,T}={1\over{2i\pi}}\int_{\sqrt{t}\delta}{{h(\lambda)}\over{\lambda-\widetilde{B}_{t,T}}}d\lambda,
\end{align}
\begin{align}\label{bd5.87}
h'(\widetilde{B}_{t,T})\widetilde{\mathbb{\bf P}}^{[0,1]}_{t,T}={1\over{2i\pi}}\int_{\sqrt{t}\delta}{{h'(\lambda)}\over{\lambda-\widetilde{B}_{t,T}}}d\lambda.
\end{align}
Similarly, 
\begin{align}\nonumber
{\widehat{\mathbb{\bf P}}}^{[0,1]}_{t,T}={1\over{2i\pi}}\int_{\sqrt{t}\delta}{{d\lambda}\over{\lambda-\widetilde{D}_{t,T}}},
\end{align}
\begin{align}\nonumber
h\left(\widetilde{D}_{t,T}\right){\widehat{\mathbb{\bf P}}}^{[0,1]}_{t,T}={1\over{2i\pi}}\int_{\sqrt{t}\delta}{{h(\lambda)}\over{\lambda-\widetilde{D}_{t,T}}}d\lambda,
\end{align}
\begin{align}\nonumber
h'\left(\widetilde{D}_{t,T}\right){\widehat{\mathbb{\bf P}}}^{[0,1]}_{t,T}={1\over{2i\pi}}\int_{\sqrt{t}\delta}{{h'(\lambda)}\over{\lambda-\widetilde{D}_{t,T}}}d\lambda,
\end{align}

\begin{align}
\widetilde{L}_{t,T}=h'\left(\widetilde{D}_{t,T}\right){\widehat{\mathbb{\bf P}}}^{[0,1]}_{t,T}.
\end{align}

\begin{definition}
For $t\in\mathbb{\bf R}^*_{+}$, $T\geq T'_0$, put
\begin{align}\nonumber
\widetilde{a}_{t,T}=\sqrt{2i\pi}\varphi {\rm Tr}_{\Gamma,s}\left[h\left(\widetilde{B}_{t,T}\right)\widetilde{\mathbb{\bf P}}^{[0,1]}_{t,T}\right],
\end{align}
\begin{align}\label{ad5.89}
\widetilde{b}_{t,T}={1\over 2}\varphi{\rm Tr}_{\Gamma,s}\left[\left(N-{n\over 2}\right)h'\left(\widetilde{B}_{t,T}\right)\widetilde{\mathbb{\bf P}}^{[0,1]}_{t,T}\right].
\end{align}
\end{definition}

In (\ref{ad5.89}), we may replace $\widetilde{B}_{t,T}$ by $\widetilde{D}_{t,T}$ and $\widetilde{\bf P}^{[0,1]}_{[t,T]}$ by $\widehat{\bf P}^{[0,1]}_{t,T}$. Then $\widetilde{a}_{t,T}$, $\widetilde{b}_{t,T}$ are forms on $S$.

Let $\widetilde{I}^0$ be the vector bundle of $L_2$ sections of $\Lambda^{\bullet}(T^*\widetilde{X})\widehat{\otimes}\widetilde{F}$ along the fibers $\widetilde{X}$, and let $\|\ \|_0$ be the norm of $\widetilde{I}^0$ associated to the Hermitian product (\ref{ad2.5}). If $L\in \mathcal{L}(\widetilde{I}^0)$, for $p\geq 1$, put
\begin{align}\label{ad5.90}
\|L\|_{\Gamma,p}={\rm Tr}_{\Gamma}\left[\left(L^* L\right)^{p/2}\right]^{1/p}.
\end{align}
Then (\ref{ad5.90}) defines a norm on a vector subspace of $\mathcal{L}(\widetilde{I}^0)$. For $p=1$, we get the $\Gamma$-trace class operators.

As in \cite[(10.159)]{BG}, given $t'>0$, if $t$ is close enough to $t'$, instead of (\ref{bd5.87}), we can write 
\begin{align}\label{bd5.91}
h\left(\widetilde{B}_{t,T}\right)\widetilde{\bf P}^{[0,1]}_{t,T}={1\over{2i\pi}}\int_{\sqrt{t'}\delta}{{h(\lambda)}\over{\lambda-\widetilde{B}_{t,T}}}d\lambda,
\end{align}
the key point being that the contour of integration in (\ref{bd5.91}) does not depend on $t$.

For $T\geq T'_0$, there exists $d_1>0$ such that 
\begin{align}\nonumber
\left|{\rm Sp}\left(\widetilde{B}_{T}\right)\right|\subset \delta\cup (2d_1,+\infty).
\end{align}

Set 
\begin{multline}\nonumber
\Delta''_+=\left\{x+iy|x=-{{d_1}\over {2}},+\infty>y\geq d_1\right\}\cup \left\{x+iy|-{{d_1}\over 2}\leq x\leq {d_1\over 2},y=d_1\right\}\\
\cup \left\{x+iy|x={d_1\over 2},d_1\leq y<+\infty\right\},
\end{multline}
\begin{multline}\nonumber
\Delta''_-=\left\{x+iy|x={d_1\over 2}, -\infty<y\leq -d_1\right\}\cup \left\{x+iy|-{d_1\over 2}\leq x\leq {d_1\over 2}, y=-d_1\right\}\\
\cup \left\{x+iy|x=-{d_1\over 2}, -d_1\geq y>-\infty\right\},
\end{multline}
\begin{align}\nonumber
\Delta''=\Delta''_+\cup \Delta''_-.
\end{align}

Set
\begin{align}\label{ad92}
\widetilde{G}_{t,T}=\psi^{-1}_t {1\over{2i\pi}}\int_{\Delta''}{{h(\sqrt{t}\lambda)}\over{\lambda-\widetilde{B}_T}}d\lambda\psi_t.
\end{align}
Then 
\begin{align}\label{ad93}
\widetilde{a}_{t,T}=\sqrt{2i\pi}\varphi {\rm Tr}_{\Gamma,s}\left[h\left(\widetilde{B}_{t,T}\right)\right]-\sqrt{2i\pi}\varphi {\rm Tr}_{\Gamma,s}\left[\widetilde{G}_{t,T}\right].
\end{align}

\begin{proposition}\label{ap5.28}
For $T\geq T'_0$, there exist $C>0$, $c>0$ such that for $t\geq 1$,
\begin{align}\label{ad94}
\left\| \widetilde{G}_{t,T}\right\|_{\Gamma,1}\leq C e^{-ct}.
\end{align}
\end{proposition}
\begin{proof}
Take $p\in{\bf N}$, $p>{\rm dim}\widetilde{X}$. Let $h_p (\lambda)$ be the unique holomorphic function on ${\bf C}\setminus {\bf R}$ such that

(1) As $\lambda\to \pm i\infty$, $h_p(\lambda)\to 0$.

(2) The following identity holds,
\begin{align}
{{h_p^{(p-1)}(\lambda)}\over{(p-1)!}}=h(\lambda). 
\end{align}
Clearly, if $\lambda\in\Delta$,
\begin{align}\label{cd5.96}
|{\rm Re}(\lambda)|\leq {1\over 2}|{\rm Im}(\lambda)|. 
\end{align}
Using (\ref{cd5.96}), we find that there exist $C>0$, $C'>0$ such that if $\lambda\in \Delta''$,
\begin{align}\label{ad97}
\left|h_{p}\left(\sqrt{t}\lambda\right)\right|\leq C\exp\left(-C't|\lambda|^2\right).
\end{align}

Clearly,
\begin{align}\label{ad98}
{1\over{2i\pi}}\int_{\Delta''}{{k(\sqrt{t}\lambda)}\over{\lambda-\widetilde{B}_T}}d\lambda={1\over{2i\pi}}\int_{\Delta''}{{h_p (\sqrt{t}\lambda)}\over{\sqrt{t}^{p-1}(\lambda-\widetilde{B}_T)^p}}d\lambda.
\end{align}
If $\lambda\in\Delta''$, we have the expansion,
\begin{align}
\left(\lambda-\widetilde{B}_T\right)^{-1}=\left(\lambda-\widetilde{B}^{(0)}_T\right)^{-1}+\left(\lambda-\widetilde{B}^{(0)}_T\right)^{-1}\widetilde{B}^{(\geq 1)}_T \left(\lambda-\widetilde{B}^{(0)}_T\right)^{-1}
+\cdots
\end{align}
and the expansion only contains a finite number of terms. Then we find that there is $C>0$ such that if $\lambda\in\Delta''$,
\begin{align}\label{ad100}
\left\|\left(\lambda-\widetilde{B}_T\right)^{-1}\right\|_{\infty}\leq C.
\end{align}

Fix $\lambda_0\in\Delta''$. Since $p>{\rm dim}\widetilde{X}$, and $\widetilde{B}_T$ is a fiberwise elliptic of order $1$,
\begin{align}\label{ad101}
\left\|\left(\lambda_0-\widetilde{B}_T\right)^{-1}\right\|_{\Gamma,p}<+\infty.
\end{align}
If $\lambda\in\Delta''$,
\begin{align}\label{ad102}
\left(\lambda-\widetilde{B}_T\right)^{-1}=\left(\lambda_0-\widetilde{B}_T\right)^{-1}+(\lambda-\lambda_0)\left(\lambda_0-\widetilde{B}_T\right)^{-1}\left(\lambda-\widetilde{B}_T\right)^{-1}.
\end{align}
From (\ref{ad100})-(\ref{ad102}), we find that if $\lambda\in\Delta''$,
\begin{align}\label{ad103}
\left\|\left(\lambda-\widetilde{B}_T\right)^{-1}\right\|_{\Gamma,p}\leq C(1+|\lambda|)\left\|\left(\lambda_0-\widetilde{B}_T\right)^{-1}\right\|_{\Gamma,p}\leq C'(1+|\lambda|). 
\end{align}
Using (\ref{ad103}), we find that if $\lambda\in\Delta''$,
\begin{align}\label{ad104}
\left\|\left(\lambda-\widetilde{B}_T\right)^{-p}\right\|_{\Gamma,1}\leq C(1+|\lambda|)^p. 
\end{align}

From (\ref{ad92}), (\ref{ad97}), (\ref{ad98}) and (\ref{ad104}), we get (\ref{ad94}). The proof is completed. 

\end{proof}

\begin{theorem}
For $T\geq T'_0$, the form $\widetilde{a}_{t,T}$ is odd and closed, and the form $\widetilde{b}_{t,T}$ is even. Moreover,
\begin{align}\label{6.77}
{\partial\over{\partial t}}\widetilde{a}_{t,T}=d{\widetilde{b}_{t,T}\over t}.
\end{align}
Also as $t\to +\infty$, there exists $\gamma'>0$ such that
\begin{align}\nonumber
\widetilde{a}_{t,T}=h\left(\nabla^{H_{(2)}^{\bullet}(\widetilde{X},\widetilde{F}|_{\widetilde{X}})},g_{L_2,T}^{H_{(2)}^{\bullet}(\widetilde{X},\widetilde{F}|_{\widetilde{X}})}\right)+\mathcal{O}(t^{-\gamma'}),
\end{align}
\begin{align}\label{6.78}
\widetilde{b}_{t,T}={1\over 2}\chi'(\widetilde{F})-{n\over 4}\chi(F)+\mathcal{O}(t^{-\gamma'}).
\end{align}
\end{theorem}

\begin{proof}
By the same argument of \cite[Theorem 10.37]{BG}, we have that the form $\widetilde{a}_{t,T}$ is odd and closed, the form $\widetilde{b}_{t,T}$ is even and (\ref{6.77}) holds.

By (\ref{ad93}), Proposition \ref{ap5.28} and \cite{Schick;NonCptTorsionEst,SS}, we get the first formula in (\ref{6.78}). One can get the second formula in (\ref{6.78}) similarly. 
\end{proof}

Set 
\begin{align}
\widetilde{c}=\widetilde{a}_{t,T}+{{dt}\over{t}}\widetilde{b}_{t,T},
\end{align}
then $\widetilde{c}$ is closed.

\begin{definition}
For $T\geq T'_0$, put
\begin{align}
\widetilde{S}^{[0,1]}_{h}(T)=-\int_{1}^{+\infty}(\widetilde{b}_{t,T}-\widetilde{b}_{\infty,T}){{dt}\over{t}}.
\end{align}
\end{definition}

The following $L^2$-analogue of \cite[Proposition 10.41]{BG} holds.
\begin{proposition}\label{p5.25}
The following identity holds,
\begin{align}
\widetilde{S}^{[0,1]}_{h}(T)=-\int_{1}^{+\infty}\left(\varphi{\rm Tr}_{\Gamma,s}\left[N\widetilde{L}_{t,T}\right]-\chi'(\widetilde{F})\right){{dt}\over{2t}}.
\end{align}
\end{proposition}

Define $$U_{L^2,h}\left(\widetilde{A}^{C^{\bullet}(W^u,\widetilde{F})'},\mathfrak{g}_{t,T}^{C^{\bullet}(W^u,\widetilde{F})}\right)$$ as in \cite[Definition 2.49]{BG} using the $\Gamma$-trace.
Then for $T\geq T'_0$, in $\Omega^{\bullet}(S)/d\Omega^{\bullet}(S)$, we have
\begin{align}\label{5.71}
\widetilde{S}^{[0,1]}_{h}(T)=U_{L^2,h}\left(\widetilde{A}^{C^{\bullet}(W^u,\widetilde{F})'},\mathfrak{g}_{t,T}^{C^{\bullet}(W^u,\widetilde{F})}\right).
\end{align}

Let $\eta_{T}$ be the even form associated to $\widetilde{A}^{C^{\bullet}(W^u,\widetilde{F})'}$, $\mathfrak{g}_{T}^{C^{\bullet}(W^u,\widetilde{F})}$, which was as introduced in \cite[(2.143)]{BG} using $\Gamma$-trace.

\begin{definition}
For $u\in\mathbb{\bf R}^*_{+}$, let $\mathfrak{h}_{T,u}^{C^{\bullet}(W^u,\widetilde{F})}$ be the generalized metric on $C^{\bullet}(W^u,\widetilde{F})$,
\begin{align}\label{89}
\mathfrak{h}^{C^{\bullet}(W^u,\widetilde{F})}_{T,u}=u^{N^{C^{\bullet}(W^u,\widetilde{F})/2}}\mathfrak{g}^{C^{\bullet}(W^u,\widetilde{F})}_{T}u^{N^{C^{\bullet}(W^u,\widetilde{F})/2}}.
\end{align}
\end{definition}

By the same argument in \cite[Theorem 10.49]{BG}, the following identity holds in $\Omega^{\bullet}(S)/d\Omega^{\bullet}(S)$,
\begin{multline}\label{5.73}
-U_{L^2,h}\left(\widetilde{A}^{C^{\bullet}(W^u,\widetilde{F})'},\mathfrak{g}_{t,T}^{C^{\bullet}(W^u,\widetilde{F})}\right)-\widetilde{h}_{L^2}\left(\nabla^{H^{\bullet}_{(2)}(\widetilde{X},\widetilde{F}|_{\widetilde{X}})},g^{H^{\bullet}_{(2)}(\widetilde{X},\widetilde{F}|_{\widetilde{X}})}_{L_2,0},g^{H^{\bullet}_{(2)}(\widetilde{X},\widetilde{F}|_{\widetilde{X}})}_{L_2,T}\right)\\
=\varphi\left\{\int_0^1 \left({\rm Tr}_{\Gamma,s}\left[N^{C^{\bullet}(W^u,\widetilde{F})}h'\left(\widetilde{B}_t^{C^{\bullet}(W^u,\widetilde{F})}\right)\right]-\widetilde{\chi}'^{-}(F)\right){{dt}\over{t}}\right.\\
\left.+\int_1^{+\infty} \left({\rm Tr}_{\Gamma,s}\left[N^{C^{\bullet}(W^u,\widetilde{F})}h'\left(\widetilde{B}_t^{C^{\bullet}(W^u,\widetilde{F})}\right)\right]-{\chi}'(\widetilde{F})\right){{dt}\over{t}}\right\}\\
-\int_0^1 \left(\varphi{\rm Tr}_{\Gamma,s}\left[{1\over 2}\left(\mathfrak{h}^{C^{\bullet}(W^u,\widetilde{F})}_{T,t}\right)^{-1}{\partial\over{\partial t}}\mathfrak{h}^{C^{\bullet}(W^u,\widetilde{F})}_{T,t}h'\left(\widetilde{B}^{C^{\bullet}(W^u,\widetilde{F})}_{t,T}\right)\right]-{{\eta_T}\over{t}}\right)dt\\
-\widetilde{h}_{L^2}\left(\nabla^{H^{\bullet}_{(2)}(\widetilde{X},\widetilde{F}|_{\widetilde{X}})},g^{H^{\bullet}_{(2)}(\widetilde{X},\widetilde{F}|_{\widetilde{X}})}_{L_2,0},g^{H^{\bullet}_{(2)}(\widetilde{X},\widetilde{F}|_{\widetilde{X}})}_{C^{\bullet}(W^u,\widetilde{F})}\right)\\
+\widetilde{h}_{L^2}\left(\nabla^{C^{\bullet}(W^u,\widetilde{F})},\mathfrak{g}_T^{C^{\bullet}(W^u,\widetilde{F})},g^{C^{\bullet}(W^u,\widetilde{F})}\right).
\end{multline}

Let $\widetilde{A}^{C^{\bullet}(W^u,\widetilde{F})''}_{T,u}$ be the adjoint of $\widetilde{A}^{C^{\bullet}(W^u,F)'}$ with respect to $\mathfrak{h}^{C^{\bullet}(W^u,\widetilde{F})}_{T,u}$. Put
\begin{align}
\widetilde{B}^{C^{\bullet}(W^u,\widetilde{F})}_{T,u}={1\over 2}\left(\widetilde{A}^{C^{\bullet}(W^u,\widetilde{F})''}_{T,u}-\widetilde{A}^{C^{\bullet}(W^u,F)'}\right).
\end{align}

Let $N^{H_{(2)}^{\bullet}(\widetilde{X},\widetilde{F}|_{\widetilde{X}})}$ be the number operator of $H_{(2)}^{\bullet}(\widetilde{X},\widetilde{F}|_{\widetilde{X}})$.

\begin{theorem}\label{t12.1}
As $T\to +\infty$,
\begin{align}\nonumber
\int_{0}^{1}\left(\varphi{\rm Tr}_{\Gamma,s}\left[{1\over 2}\left(\mathfrak{h}_{T,t}^{C^{\bullet}(W^u,\widetilde{F})}\right)^{-1}{\partial\over{\partial t}}\mathfrak{h}_{T,t}^{C^{\bullet}(W^u,\widetilde{F})}h'\left(\widetilde{B}_{T,t}^{C^{\bullet}(W^u,\widetilde{F})}\right)\right]-{\eta_{T}\over{t}}\right)dt\to 0,
\end{align}
\begin{multline}\label{91}
\widetilde{h}_{L^2}\left(\nabla^{C^{\bullet}(W^u,\widetilde{F})},\mathfrak{g}_{T}^{C^{\bullet}(W^u,\widetilde{F})},g^{C^{\bullet}(W^u,\widetilde{F})}\right)-{\rm Tr}_{s}[f]T\\
-{1\over 4}\left(\widetilde{\chi}'^{+}(F)-\widetilde{\chi}'^{-}(F)\right)\log(T)\to {1\over 4}\left(\widetilde{\chi}'^{-}(F)-\widetilde{\chi}'^{+}(F)\right)\log(\pi).
\end{multline}
\end{theorem}

The rest part of this section is devoted to prove Theorem \ref{t12.1}.

Recall that the operator $\overline{\mathcal{C}}_{T}^{\widetilde{I}\widehat{\otimes}\widetilde{F}|_{\widetilde{\mathbb{\bf B}}},2}$ was constructed in Section \ref{s5.1}.

Given $n'\in {\bf N}$, let $Q^{n'}(T\widetilde{X}^u|_{\widetilde{\bf B}})$ be the algebra of invariant polynomials of degree $n'$ on $T\widetilde{X}^u|_{\widetilde{\bf B}}$. 

\begin{definition}
For $T\geq 1$, put
\begin{align}
{\overline{\mathbb{P}}}_{T}={1\over{2i\pi}}\int_{\delta}{{d\lambda}\over{\lambda-\overline{\mathcal{C}}^{\widetilde{I}\widehat{\otimes}\widetilde{F}_{\widetilde{\mathbb{\bf B}}},2}_{T}}}.
\end{align}
\end{definition}
Then $\overline{\mathbb{P}}_{T}$ is a projection acting on $\Lambda^{\bullet}(T^*S)\widehat{\otimes}\widetilde{F}|_{\widetilde{\bf B}}$, with finite $\Gamma$-dimensional range $\overline{\mathbb{F}}_T$.
Then by \cite[Theorem XII.5]{ReS}, for $k\in {\bf N}$ large enough, 
\begin{align}\label{fd5.116}
\overline{\mathbb{F}}_T={\rm ker}\overline{\mathcal{D}}_T^{\widetilde{I}\widehat{\otimes}\widetilde{F}|_{\widetilde{\bf B}},k}.
\end{align}

As \cite[(10.191)]{BG}, for $n'$ large enough, we also have,
\begin{align}\label{cd5.116}
\overline{\mathbb{F}}_T \subset \exp\left(-T|Z|^2/2\right)Q^{n'}\left(T\widetilde{X}^u|_{\widetilde{\bf B}}\right)\widehat{\otimes}\Lambda^\bullet (T^*S)\widehat{\otimes}\Lambda^\bullet (T^*\widetilde{X}|_{\widetilde{\bf B}})\widehat{\otimes}\widetilde{F}|_{\widetilde{\bf B}}. 
\end{align}

\subsection{The maps $\widetilde{{\bf J}}_{T}$ and ${\overline{\bf e}}_{T}$}
This subsection is the $L^2$-set up of \cite[Section 10.15]{BG}.

Put 
\begin{align}\label{cd5.116}
{\overline{\mathbb{\bf P}}}^{[0,1]}_{T}={1\over{2i\pi}}\int_{\delta}{{d\lambda}\over{\lambda-{\overline{A}}^2_{T}}}.
\end{align}
Then 
\begin{align}\label{94}
{\overline{\mathbb{\bf P}}}^{[0,1]}_{T}=e^{-T\widetilde{f}}\widetilde{\mathbb{\bf P}}^{[0,1]}_{T}e^{T\widetilde{f}}.
\end{align}
So ${\overline{\mathbb{\bf P}}}^{[0,1]}_{T}$ is a projector acting on $\Lambda^{\bullet}(T^*S)\widehat{\otimes}\Omega_{(2)}^{\bullet}(\widetilde{X},\widetilde{F}|_{\widetilde{X}})$.

Let $\gamma:\mathbb{\bf R}\to [0,1]$ be a smooth function such that
\begin{align}
\gamma(a)=1\ {\rm for}\ a<{1\over 2},
\ \ =0\ {\rm for}\ a>1.
\end{align}
Let $\varepsilon_0$ be chosen as in \cite[Section 9.1]{BG}. If $Z\in\mathbb{\bf R}^n$, set
\begin{align}
\mu(Z)=\gamma(|Z|/\varepsilon_0).
\end{align}
Then
\begin{align}
\mu(Z)=1\ {\rm if}\ |Z|\leq \varepsilon_0/2,
\ \ =0\ {\rm if}\ |Z|\geq \varepsilon_0.
\end{align}

If $T>0$, set
\begin{align}
\alpha_T=\int_{\mathbb{\bf R}^n}\mu^2(Z)\exp(-T|Z|^2)dZ.
\end{align}
Then there is $c>0$ such that as $T\to +\infty$, 
\begin{align}
\alpha_T=\left({\pi\over T}\right)^{n/2}+\mathcal{O}\left(e^{-cT}\right).
\end{align}

Take $x\in \widetilde{\bf B}$. Let $\rho_x\in\Lambda^{\rm max}(T^{u,*}\widetilde{X}|_{\widetilde{\bf B}})_{x}$ be of norm 1. Then $\rho_x$ is determined up to sign. It defines a section $o(x)\otimes \Lambda^{\rm max}(T^{u,*}\widetilde{X}|_{\widetilde{\bf B}})_{x}$.

\begin{definition}
Let $\widetilde{{\bf J}}_{T}:\Lambda^{\bullet}(T^*S)\widehat{\otimes}C^{\bullet}(W^u,\widetilde{F})\to \Lambda^{\bullet}(T^*S)\widehat{\otimes}\Omega_{(2)}^{\bullet}(\widetilde{X},\widetilde{F}|_{\widetilde{X}})$ be such that if $h\in o(x)\otimes \widetilde{F}_x$, then
\begin{align}
\widetilde{{\bf J}}_T h={{\mu(Z)}\over{\alpha^{1/2}_{T}}}{\overline{\mathbb{P}}}_{T}\left[\exp\left(-T|Z|^2/2\right)\rho_x\right]h.
\end{align} 
\end{definition}

The induced map $\widetilde{J}_T:C^{\bullet}(W^u,\widetilde{F})\to \Omega_{(2)}^{\bullet}(\widetilde{X},\widetilde{F}|_{\widetilde{X}})$ is given by
\begin{align}\label{ad5.124}
\widetilde{J}_T h={{\mu(Z)}\over{\alpha^{1/2}_T}}\exp\left(-T|Z|^2/2\right)\rho_x h.
\end{align}

\begin{definition}
Let ${\overline{\bf e}}_T:\Lambda^{\bullet}(T^*S)\widehat{\otimes}C^{\bullet}(W^u,\widetilde{F})\to {\overline{\bf F}}^{[0,1]}_{T}$ be given by 
\begin{align}\label{102}
{\overline{\bf e}}_T={\overline{\bf P}}^{[0,1]}_{T}\widetilde{{\bf J}}_T.
\end{align}
\end{definition}
The induced map ${\overline{e}}_T:C^{\bullet}(W^u,\widetilde{F})\to \Omega_{(2)}^{\bullet}(\widetilde{X},\widetilde{F}|_{\widetilde{X}})$ is given by
\begin{align}
{\overline{e}}_T={\overline{P}}^{[0,1]}_{T}\widetilde{J}_T.
\end{align}

In the sequel, we write that as $T\to +\infty$, a family of smooth sections on $\widetilde{M}$ is $\mathcal{O}(e^{-cT})$ if the sup norm of the derivatives is $\mathcal{O}(e^{-cT})$. By the same argument of \cite[Theorem 10.56]{BG}, we have the following $L^2$-extension.

\begin{theorem}\label{t5.31}
There is $c>0$ such that as $T\to +\infty$, for any $s\in C^{\bullet}(W^u,\widetilde{F})$,
\begin{align}\label{cd5.127}
\left({\overline{\bf e}}_{T}-\widetilde{{\bf J}}_{T}\right)s=\mathcal{O}\left(e^{-cT}\right)\ uniformly\ on\ \widetilde{M}.
\end{align}
\end{theorem}

\begin{proof}
By (\ref{cd5.116}) and holomorphic functional calculus, we know that for any $k\in{\bf N}^*$, 
\begin{align}\label{cd5.128}
\overline{\bf P}^{[0,1]}_{T}={1\over{2i\pi}}\int_{\delta}{{d\lambda}\over{\lambda-\overline{A}^{2k}_{T}}}.
\end{align}
In the sequel, we choose $k\in{\bf N}^*$ large enough so that (\ref{fd5.116}) holds.

Take $x\in \widetilde{B}$, $h\in o(x)\otimes \widetilde{F}_x$. If $\lambda\in\delta$,
\begin{align}
\left(\lambda-\overline{A}^{2k}_T\right){{\widetilde{\bf J}_T h}\over{\lambda}}-{\widetilde{\bf J}}_T h=-\overline{A}^{2k}_T{{\widetilde{\bf J}_T h}\over{\lambda}}, 
\end{align}
and so,
\begin{align}\label{cd5.130}
{{\widetilde{\bf J}_T h}\over{\lambda}}-\left(\lambda-\overline{A}^{2k}_T\right)^{-1}\widetilde{\bf J}_T h=-\left(\lambda-\overline{A}^{2k}_T\right)^{-1} \overline{A}^{2k}_T {{\widetilde{\bf J}_T h}\over{\lambda}}.
\end{align}
Now, by the fundamental assumptions in \cite[Section 9.1]{BG}, with the required identifications, on $\{x'\in\widetilde{X},d^{\widetilde{X}}(x,x')\leq \varepsilon\}$, the operator $\overline{A}^2_T$ coincides with $\overline{\mathcal{C}}^{\widetilde{I}\widehat{\otimes}\widetilde{F}|_{\widetilde{\bf B}},2}_{T}$. Since $\mu(Z)=1$ for $|Z|\leq \varepsilon/2$, using (\ref{fd5.116}), we get
\begin{align}\label{cd5.132}
\overline{A}^{2k}_T \widetilde{J}_T h(Z)=0\ \ {\rm for}\ |Z|\leq \varepsilon/2.
\end{align}
By (\ref{cd5.116}), we deduce from (\ref{cd5.132}) that there exists $c>0$ such that as $T\to +\infty$,
\begin{align}\label{cd5.133}
\left|\overline{A}^{2k}_T \widetilde{\bf J}_T h\right|=\mathcal{O}\left(e^{-cT}\right)|h|.
\end{align}

Since 
\begin{align}\nonumber
\widetilde{A}={1\over 2} D^{\widetilde{X}}+\nabla^{\Omega_{(2)}^\bullet(\widetilde{X},\widetilde{F}|_{\widetilde{X}})}-{1\over 2}c(T^H),
\end{align}
we see that
\begin{align}\label{cd5.134}
\overline{A}^{2k}_{T}=\left(D^{\widetilde{X}}\over 2\right)^{2k}+\widetilde{K}_T,
\end{align}
where $\widetilde{K}_T$ is a differential operator of order $2k-1$, whose coefficients depend polynomially on $T$, the polynomial being of degree $2k$.

If $q\in {\bf N}$, let $\|\ \|_q$ be the norm on the fiberwise $q^{\rm th}$ Sobolev space of sections of $\Lambda^\bullet (T^*\widetilde{X})\widehat{\otimes}\Lambda^\bullet (T^* S)$ defined in \cite{BFKM,SS}.

Since $D^{\widetilde{X},2k}$ is elliptic of degree $2k$, given $q\in {\bf N}$, by \cite[Lemma 1.4]{Shu} (see Lemma \ref{r2lA1}), there exists $C>0$ such that for $s\in \Lambda^\bullet (T^* S)\widehat{\otimes}\Omega^\bullet_{(2)}(\widetilde{X},\widetilde{F}|_{\widetilde{X}})$, 
\begin{align}\label{cd5.135}
\|s\|_{q+2k}\leq C\left(\left\|D^{\widetilde{X},2k}s\right\|_q+\|s\|_0\right).
\end{align}
By the considerations which follow (\ref{cd5.134}) and (\ref{cd5.135}), we see that given $q\in {\bf N}$, there exists $C>0$ such that for $\lambda\in\delta$, $T\geq 1$,
\begin{align}\label{cd5.136}
\|s\|_{q+2k}\leq C\left(\left\|\left(\lambda-\overline{A}^{2k}_T\right)s\right\|_q+T^{2k}\|s\|_{q+2k-1}\right).
\end{align}
Also given $q\in{\bf N}$, there exists $C>0$ such that for $A>0$, $s\in\Lambda^\bullet(T^*S)\widehat{\otimes}\Omega^\bullet_{(2)}(\widetilde{X},\widetilde{F}|_{\widetilde{X}})$,
\begin{align}\label{cd5.137}
\|s\|_{q+2k-1}\leq C\left({{\|s\|_{q+2k}}\over A}+A^{q+2k-1}\|s\|_0\right).
\end{align}
From (\ref{cd5.136}), (\ref{cd5.137}), we deduce that there exists $C>0$, $k'\in{\bf N}$ such that for $\lambda\in\delta$, $s\in \Lambda^\bullet (T^*S)\widehat{\otimes}\Omega^\bullet_{(2)}(\widetilde{X},\widetilde{F}|_{\widetilde{X}})$,
\begin{align}\label{cd5.138}
\|s\|_{q+2k}\leq C\left(\left\|\left(\lambda-\overline{A}^{2k}_T\right)\right\|_{q}+T^{2k'}\|s\|_0\right).
\end{align}

Also by (\ref{6.19}), for $T\geq T_0$, we know that ${\rm Sp}(\overline{A}_T^{2k,(0)})\cap \delta=\emptyset$. More precisely, since $\overline{A}^{(2k),(0)}_T$ is self-adjoint, by (\ref{6.19}), there exists $C'>0$ such that for $T\geq T_0$, $\lambda\in\delta$, $s\in\Omega^\bullet_{(2)}(\widetilde{X},\widetilde{F}|_{\widetilde{X}})$,
\begin{align}\label{cd5.139}
\left\|\left(\lambda-\overline{A}_T^{2k,(0)}\right)^{-1}s\right\|_0\leq C'\|s\|_0.
\end{align}
By (\ref{cd5.138}), (\ref{cd5.139}), we get
\begin{align}\label{cd5.140}
\left\|\left(\lambda-\overline{A}_T^{2k,(0)}\right)^{-1}s\right\|_{q+2k}\leq C' T^{2k'}\|s\|_q.
\end{align}
Moreover,
\begin{align}\label{cd5.141}
\left(\lambda-\overline{A}^{2k}_T\right)^{-1}=\left(\lambda-\overline{A}_T^{2k,(0)}\right)^{-1}+\left(\lambda-\overline{A}_T^{2k,(0)}\right)^{-1}\overline{A}_T^{2k,(0)}\left(\lambda-\overline{A}_T^{2k,(0)}\right)^{-1}+\cdots 
\end{align}
and the expansion in (\ref{cd5.141}) contains a finite number of terms. Also by (\ref{ad5.35}), $\overline{A}_T^{2k,(0)}$ is a differential operator of order $2k-1$, which depends polynomially on $T$. 
By (\ref{cd5.140}), (\ref{cd5.141}), we find that there exists $C''>0$, $k''\in{\bf N}$ such that for $\lambda\in\delta$, $T\geq T_0$, $s\in\Omega^\bullet_{(2)}(\widetilde{X},\widetilde{F}|_{\widetilde{X}})$,
\begin{align}\label{cd5.142}
\left\|\left(\lambda-\overline{A}^{2k}_T\right)^{-1}s\right\|_{q+2k}\leq C'' T^{2k'}\|s\|_q.
\end{align}

From (\ref{cd5.133}), (\ref{cd5.142}), we deduce that there exists $c>0$ such that for $T\geq T_0$,
\begin{align}\label{cd5.143}
\left\|\left(\lambda-\overline{A}^{2k}_T\right)\overline{A}^{2k}_T \widetilde{\bf J}_Th\right\|_{q+2k}=\mathcal{O}\left(e^{-cT}\right)\|h\|.
\end{align}
Using (\ref{cd5.143}) and Sobolev's inequalities (cf. \cite{H}), we see that there exists $c>0$ such that for $\lambda\in\delta$, $T\geq T_0$, $x\in \widetilde{B}$, $h\in\widetilde{F}_x$,
\begin{align}\label{cd5.144}
\left|\left(\lambda-\overline{A}^{2k}_T\right)^{-1}\overline{A}^{2k}_T \widetilde{\bf J}_T h\right|=\mathcal{O}\left(e^{-cT}\right)\|h\|.
\end{align}
From (\ref{cd5.128}), (\ref{cd5.130}), (\ref{cd5.144}), we get (\ref{cd5.127}). The proof is completed.

\end{proof}

\begin{definition}
For $T\geq T_0$, let $\widetilde{{\bf e}}_{T}:\Lambda^{\bullet}(T^*S)\widehat{\otimes}C^{\bullet}(W^u,\widetilde{F})\to \widetilde{{\bf F}}^{[0,1]}_{T}$ be the linear map,
\begin{align}
\widetilde{{\bf e}}_{T}=e^{T\widetilde{f}}{\overline{\bf e}}_{T}.
\end{align}
\end{definition}
By (\ref{94}), (\ref{102}), we get
\begin{align}
\widetilde{{\bf e}}_{T}=\widetilde{{\bf P}}^{[0,1]}_{T}e^{T\widetilde{f}}\widetilde{{\bf J}}_{T}.
\end{align}
Then $\widetilde{{\bf e}}_{T}$ commutes with $\Lambda^{\bullet}(T^*S)$. The induced map $\widetilde{e}_T: C^{\bullet}(W^u,\widetilde{F})\to \Omega_{(2)}^{\bullet}(\widetilde{X},\widetilde{F}|_{\widetilde{X}})$ is given by,
\begin{align}
\widetilde{e}_T=\widetilde{P}_T^{[0,1]}e^{T\widetilde{f}}\widetilde{J}_T.
\end{align}

In the sequel, we consider $\widetilde{{\bf e}}_{T}$ as a linear map from $\Lambda^{\bullet}(T^*S)\widehat{\otimes} C^{\bullet}(W^u,\widetilde{F})$ into $\Lambda^{\bullet}(T^*S)\widehat{\otimes}\Omega_{(2)}^{\bullet}(\widetilde{X},\widetilde{F}|_{\widetilde{X}})$. Let $\widetilde{{\bf e}}^*_T: \Lambda^{\bullet}(T^*S)\widehat{\otimes}\Omega_{(2)}^{\bullet}(\widetilde{X},\widetilde{F}|_{\widetilde{X}})\to \Lambda^{\bullet}(T^*S)\widehat{\otimes}C^{\bullet}(W^u,\widetilde{F})$ be the adjoint of $\widetilde{{\bf e}}_T$ with respect to the metrics $g^{C^{\bullet}(W^u,\widetilde{F})}$, $g_{T}^{\Omega_{(2)}^{\bullet}(\widetilde{X},\widetilde{F}|_{\widetilde{X}})}$.

Recall that $C^{\bullet}(W^u,{\bf R})=\oplus _{x\in\widetilde{B}}o(x)$. We will denote by $\mathcal{O}_{D}(1/T)$ an element of $\Lambda^{\bullet}(T^*S)\widehat{\otimes}{\rm End}(C^{\bullet}(W^u,{\bf R}))$ which commutes with $\Lambda^{\bullet}(T^*S)$, which is of positive degree in $\Lambda^{\bullet}(T^*S)$, which preserves the $\Lambda^{\bullet}(T^*S)\widehat{\otimes}o(x)$, which is $\mathcal{O}(1/T)$ as $T\to +\infty$. Of course $\mathcal{O}_{D}(1/T)$ then acts on $C^{\bullet}(W^u,\widetilde{F})$, and preserves the $o(x)\otimes \widetilde{F}_{x}$, $x\in \widetilde{B}$. Also the various $\mathcal{O}_{D}(1/T)$ commute with each other.

Using Theorem \ref{t5.31}, the following $L^2$-extension of \cite[Proposition 10.58]{BG} can be proved using the same argument of \cite[Proposition 10.58]{BG}.

\begin{proposition}\label{14.5}
As $T\to +\infty$,
\begin{align}
\widetilde{{\bf e}}_{T}^* \widetilde{{\bf e}}_{T}=1+\mathcal{O}_{D}(1/T)+\mathcal{O}\left(e^{-cT}\right).
\end{align}
\end{proposition}

Clearly, $\widetilde{{\bf P}}^{\infty}_T\widetilde{{\bf e}}_T\in\Lambda^{\bullet}(T^*S)\widehat{\otimes}{\rm End}(C^{\bullet}(W^u,\widetilde{F}))$. Let $\mathcal{F}:C^*(W^u,\widetilde{F})\to C^*(W^u,\widetilde{F})$ be acting on $\widetilde{F}_{x}\otimes o^u_{x}$, $x\in\widetilde{B}$, by multiplication by $\widetilde{f}(x)$. Also we still denote by $N:C^*(W^u,\widetilde{F})\to C^*(W^u,\widetilde{F})$ the operator acting on $C^i(W^u,\widetilde{F})$ by multiplication by $i$. The following $L^2$-extension of \cite[Theorem 10.59]{BG}
can be proved by the same argument of \cite[Theorem 10.59]{BG}.

\begin{theorem}\label{14.6}
There exists $c>0$ such that as $T\to +\infty$,
\begin{align}
\widetilde{{\bf P}}^{\infty}_T \widetilde{{\bf e}}_{T}=e^{T\mathcal{F}}\left({\pi\over T}\right)^{N^{C^{\bullet}(W^u,\widetilde{F})/2}-n/4}\left(1+\mathcal{O}_{D}(1/T)+\mathcal{O}\left(e^{-cT}\right)\right).
\end{align}
In particular for $T\geq T'_0$ large enough, $\widetilde{{\bf P}}^{\infty}_{T}{\bf e}_T\in(\Lambda^{\bullet}(T^*S)\widehat{\otimes}{\rm End}(C^{\bullet}(W^u,\widetilde{F})))$ is invertible. 
\end{theorem}

\subsection{A proof of the first part of Theorem \ref{t12.1}}
This subsection is the $L^2$-set up of \cite[Section 10.16]{BG}.

By Theorem \ref{14.6}, for $T\geq T'_0$ large enough, the map $\widetilde{{\bf P}}^{\infty}_T \widetilde{{\bf e}}_{T}$ is invertible. Therefore, for $T$ large enough,
\begin{align}\label{110}
\mathfrak{g}_{T}^{C^{\bullet}(W^u,\widetilde{F})}=\left(\widetilde{{\bf P}}^{\infty}_{T}\widetilde{{\bf e}}_{T}\right)^{*,-1}(\widetilde{{\bf e}}^*_T \widetilde{{\bf e}}_{T})\left(\widetilde{{\bf P}}^{\infty}_{T}\widetilde{{\bf e}}_{T}\right)^{-1}.
\end{align}
By (\ref{89}),
\begin{multline}\label{111}
e^{-T\mathcal{F}}t^{N^{C^{\bullet}(W^u,\widetilde{F})/2}}\widetilde{A}^{C^{\bullet}(W^u,\widetilde{F})'}t^{-N^{C^{\bullet}(W^u,\widetilde{F})/2}}e^{T\mathcal{F}}\\
=\left(\sqrt{t}e^{-T\mathcal{F}}\widetilde{\partial} e^{T\mathcal{F}}+\nabla^{C^{\bullet}(W^u,\widetilde{F})}+Td\mathcal{F}\right),
\end{multline}
\begin{multline}\nonumber
e^{-T\mathcal{F}}t^{N^{C^{\bullet}(W^u,\widetilde{F})/2}}\widetilde{A}_{T,t}^{C^{\bullet}(W^u,\widetilde{F})''}t^{-N^{C^{\bullet}(W^u,\widetilde{F})/2}}e^{T\mathcal{F}}\\
=\left(\mathfrak{g}_{T}^{C^{\bullet}(W^u,\widetilde{F})}e^{T\mathcal{F}}\right)^{-1}\left(\sqrt{t}\widetilde{\partial}^*+\nabla^{C^{\bullet}(W^u,\widetilde{F}),*}\right)\mathfrak{g}_{T}^{C^{\bullet}(W^u,\widetilde{F})}e^{T\mathcal{F}}.
\end{multline}

Set 
\begin{align}\label{112}
\widetilde{{\bf k}}_T=e^{-T\mathcal{F}}\widetilde{{\bf P}}^{\infty}_{T}\widetilde{{\bf e}}_{T}.
\end{align}
Clearly, by (\ref{110}), (\ref{112}), we get
\begin{align}\label{113}
e^{T\mathcal{F}}\mathfrak{g}_{T}^{C^{\bullet}(W^u,\widetilde{F})}e^{T\mathcal{F}}=(\widetilde{{\bf k}}_T)^{-1,*}\left(\widetilde{{\bf e}}^*_{T}\widetilde{{\bf e}}_{T}\right)(\widetilde{{\bf k}}_{T})^{-1}.
\end{align}

From (\ref{110})-(\ref{113}), we obtain,
\begin{multline}
e^{-T\mathcal{F}}t^{N^{C^{\bullet}(W^u,\widetilde{F})/2}}\widetilde{A}_{T,t}^{C^{\bullet}(W^u,\widetilde{F})''}t^{-N^{C^{\bullet}(W^u,\widetilde{F})/2}}e^{T\mathcal{F}}=\widetilde{{\bf k}}_T \left(\widetilde{{\bf e}}^*_T \widetilde{{\bf e}}_{T}\right)^{-1}\widetilde{{\bf k}}^*_{T}\\
\left(\sqrt{t}e^{T\mathcal{F}}\widetilde{\partial}^* e^{-T\mathcal{F}}+\nabla^{C^{\bullet}(W^u,\widetilde{F}),*}-Td\mathcal{F}\right)\left(\widetilde{{\bf k}}_T \left(\widetilde{{\bf e}}^*_T \widetilde{{\bf e}}_T\right)^{-1}\widetilde{{\bf k}}^*_T\right)^{-1}.
\end{multline}
Also by Proposition \ref{14.5} and Theorem \ref{14.6}, as $T\to +\infty$,
\begin{align}\label{115}
\widetilde{{\bf k}}_T (\widetilde{{\bf e}}^*_T \widetilde{{\bf e}}_T)^{-1}\widetilde{{\bf k}}^*_T=\left({\pi\over T}\right)^{N^{C^{\bullet}(W^u,\widetilde{F})-n/2}}\left(1+\mathcal{O}_{D}(1/T)\right)+\mathcal{O}\left(e^{-cT}\right).
\end{align}

By (\ref{115}), we deduce that
\begin{multline}\label{116}
\widetilde{{\bf k}}_T (\widetilde{{\bf e}}^*_T \widetilde{{\bf e}}_T)^{-1}\widetilde{{\bf k}}^*_T (-Td\mathcal{F})\left(\widetilde{{\bf k}}_T (\widetilde{{\bf e}}^*_T \widetilde{{\bf e}}_T)^{-1}\widetilde{{\bf k}}^*_T\right)^{-1}
=-Td\mathcal{F}+\mathcal{O}_{D}(1/T)+\mathcal{O}\left(e^{-cT}\right).
\end{multline}

Similar as \cite[(10.251)]{BG}, there exists $c>0$ such that as $T\to +\infty$,
\begin{align}\label{117}
e^{-T\mathcal{F}}\widetilde{\partial} e^{T\mathcal{F}}=\mathcal{O}\left(e^{-cT}\right),\ e^{T\mathcal{F}}\widetilde{\partial}^* e^{-T\mathcal{F}}=\mathcal{O}\left(e^{-cT}\right).
\end{align}

From (\ref{111})-(\ref{117}), we deduce that given $t\in (0,1]$, as $T\to +\infty$,
\begin{multline}\label{118}
e^{-T\mathcal{F}}t^{N^{C^{\bullet}(W^u,\widetilde{F})/2}}\widetilde{A}_{T,t}^{C^{\bullet}(W^u,\widetilde{F})}t^{-N^{C^{\bullet}(W^u,\widetilde{F})/2}}e^{T\mathcal{F}}=\nabla^{C^{\bullet}(W^u,\widetilde{F}),u}\\
+\mathcal{O}_{D}(1/T)+\left(1+\sqrt{t}\right)\mathcal{O}\left(e^{-cT}\right).
\end{multline}

By \cite[(2.143) and (2.150)]{BG},
\begin{multline}\label{119}
\int_{0}^{1}\left(\varphi{\rm Tr}_{\Gamma,s}\left[{1\over 2}\left(\mathfrak{h}_{T,t}^{C^{\bullet}(W^u,\widetilde{F})}\right)^{-1}{\partial\over{\partial t}}\mathfrak{h}_{T,t}^{C^{\bullet}(W^u,\widetilde{F})}h'\left(\widetilde{B}_{T,t}^{C^{\bullet}(W^u,\widetilde{F})}\right)\right]-{{\eta_T}\over{t}}\right)dt\\
=\int_{0}^{1}\varphi{\rm Tr}_{\Gamma,s}\left[\left(N^{C^{\bullet}(W^u,\widetilde{F})}+\left(\mathfrak{g}_{T}^{C^{\bullet}(W^u,\widetilde{F})}\right)^{-1}N^{C^{\bullet}(W^u,\widetilde{F})}\mathfrak{g}_{T}^{C^{\bullet}(W^u,\widetilde{F})}\right)\right.\\
\left. \left(h'\left(\widetilde{B}_{T,t}^{C^{\bullet}(W^u,\widetilde{F})}\right)-h'\left(\widetilde{B}_{T,0}^{C^{\bullet}(W^u,\widetilde{F})}\right)\right)\right]{{dt}\over{2t}}.
\end{multline}
Also by (\ref{113}),
\begin{multline}\label{120}
e^{-T\mathcal{F}}\left(N^{C^{\bullet}(W^u,\widetilde{F})}+\left(\mathfrak{g}_{T}^{C^{\bullet}(W^u,\widetilde{F})}\right)^{-1}N^{C^{\bullet}(W^u,\widetilde{F})}\mathfrak{g}_{T}^{C^{\bullet}(W^u,\widetilde{F})}\right)e^{T\mathcal{F}}\\
=N^{C^{\bullet}(W^u,\widetilde{F})}+\widetilde{{\bf k}}_T(\widetilde{{\bf e}}^*_T \widetilde{{\bf e}}_T)^{-1}\widetilde{{\bf k}}^*_T N^{C^{\bullet}(W^u,\widetilde{F})}\left(\widetilde{{\bf k}}_T (\widetilde{{\bf e}}^*_T \widetilde{{\bf e}}_T)^{-1}\widetilde{{\bf k}}^*_T\right)^{-1}.
\end{multline}
By (\ref{115}) and (\ref{120}), it is clear that
$$\widetilde{{\bf k}}_T(\widetilde{{\bf e}}^*_T \widetilde{{\bf e}}_T)^{-1}\widetilde{{\bf k}}^*_T N^{C^{\bullet}(W^u,\widetilde{F})}\left(\widetilde{{\bf k}}_T (\widetilde{{\bf e}}^*_T \widetilde{{\bf e}}_T)^{-1}\widetilde{{\bf k}}^*_T\right)^{-1}$$
remains uniformly bounded as $T\to +\infty$.

From (\ref{118}), and from the above boundedness result, we see that as $T\to +\infty$, the integrand in the right-hand side of (\ref{119}) tends to $0$. Moreover, by (\ref{118}), we can use dominated convergence in this integral, which tends to $0$ as $T\to +\infty$. We have thus established the first convergence result in Theorem \ref{t12.1}.

\subsection{A proof of the second part of Theorem \ref{t12.1}}
This subsection is the $L^2$-set up of \cite[Section 10.17]{BG}.

Clearly,
\begin{multline}\label{121}
\widetilde{h}_{L^2}\left(\nabla^{C^{\bullet}(W^u,\widetilde{F})},\mathfrak{g}_{T}^{C^{\bullet}(W^u,\widetilde{F})},g^{C^{\bullet}(W^u,\widetilde{F})}\right)\\
=\widetilde{h}_{L^2}\left(\nabla^{C^{\bullet}(W^u,\widetilde{F})},\mathfrak{g}_{T}^{C^{\bullet}(W^u,\widetilde{F})},e^{-2T\mathcal{F}}g^{C^{\bullet}(W^u,\widetilde{F})}\right)\\
+\widetilde{h}_{L^2}\left( \nabla^{C^{\bullet}(W^u,\widetilde{F})},e^{-2T\mathcal{F}}g^{C^{\bullet}(W^u,\widetilde{F})}, g^{C^{\bullet}(W^u,\widetilde{F})}\right).
\end{multline}
By (\ref{ad4.1}), we have 
\begin{align}
\widetilde{h}_{L^2}\left( \nabla^{C^{\bullet}(W^u,\widetilde{F})},e^{-2T\mathcal{F}}g^{C^{\bullet}(W^u,\widetilde{F})}, g^{C^{\bullet}(W^u,\widetilde{F})}\right)=T {\rm Tr}_{s}[f].
\end{align}
Also,
\begin{multline}\label{123}
\widetilde{h}_{L^2}\left(\nabla^{C^{\bullet}(W^u,\widetilde{F})},\mathfrak{g}_{T}^{C^{\bullet}(W^u,\widetilde{F})},e^{-2T\mathcal{F}}g^{C^{\bullet}(W^u,\widetilde{F})}\right)\\
=\widetilde{h}_{L^2}\left(\nabla^{C^{\bullet}(W^u,\widetilde{F})}+Td\mathcal{F},e^{T\mathcal{F}}\mathfrak{g}_{T}^{C^{\bullet}(W^u,\widetilde{F})}e^{T\mathcal{F}},g^{C^{\bullet}(W^u,\widetilde{F})}\right).
\end{multline}
Using (\ref{113}), (\ref{115}), (\ref{116}) and (\ref{123}), we find that as $T\to +\infty$,
\begin{multline}\label{124}
\widetilde{h}_{L^2}\left(\nabla^{C^{\bullet}(W^u,\widetilde{F})},\mathfrak{g}_{T}^{C^{\bullet}(W^u,\widetilde{F})},e^{-2T\mathcal{F}}g^{C^{\bullet}(W^u,\widetilde{F})}\right)\\+{1\over 4}\left(\widetilde{\chi}'^{-}(F)-\widetilde{\chi}'^{+}(F)\right)\log (T)
\to {1\over 4}\left(\widetilde{\chi}'^{-}(F)-\widetilde{\chi}'^{+}(F)\right)\log (\pi).
\end{multline}

From (\ref{121})-(\ref{124}), we get the second equation in (\ref{91}). The proof of Theorem \ref{t12.1} is completed. Then by Theorem \ref{t5.5}, Proposition \ref{p5.25}, (\ref{5.71}), (\ref{5.73}) and Theorem \ref{t12.1}, we get Theorem \ref{t4.7}.

\section{A proof of Theorem \ref{t3}}
\setcounter{equation}{0}

By \cite[Section 5.5]{BG}, we may assume that $f$ is fiberwise nice, also we may assume that the conditions in \cite[Section 9.1]{BG} hold. Then the conditions in \cite[Section 11.1]{BG} hold.

Let $d_{T}^{\widetilde{X},*}$ be the adjoint of $d^{\widetilde{X}}$ with respect to the metrics $g^{T\widetilde{X}}$, $g^{\widetilde{F}}_{T}$. Let $\nabla_{T}^{\Omega_{(2)}^{\bullet}(\widetilde{X},\widetilde{F}|_{\widetilde{X}}),*}$ be the corresponding adjoint connection to $\nabla^{\Omega_{(2)}^{\bullet}(\widetilde{X},\widetilde{F}|_{\widetilde{X}})}$, and let $\widetilde{A}''_T$ be the adjoint superconnection on $\Omega^{\bullet}_{(2)}(\widetilde{X},\widetilde{F}|_{\widetilde{X}})$ to $\widetilde{A}'$ with respect $g^{T\widetilde{X}}$, $g^{\widetilde{F}}_{T}$.

Set 
\begin{align}
\widetilde{A}_T={1\over 2}\left(\widetilde{A}''_T+\widetilde{A}'\right),\ \ \widetilde{B}_{T}={1\over 2}\left(\widetilde{A}''_T-\widetilde{A}'\right).
\end{align}

Take $c_1\in (0,1]$. Put
\begin{align}
U_T=\left\{\lambda\in\mathbb{C},{1\over 4}|\lambda|\leq c_1\sqrt{T},|\lambda|\geq {1\over 8}\right\}.
\end{align}
Then we have the following $L^2$-analogue of \cite[Theorems 11.16 and 11.17]{BG}. The proofs are the same as \cite[Theorems 11.16 and 11.17]{BG}.

\begin{thm}\label{t1}
For $c_1\in (0,1]$ small enough, for any integer $p\geq {\rm dim}X+2$, there exists $C>0$ such that for $T\geq 1$ and $\lambda\in U_T$, 
\begin{multline}\label{ad6.3}
\left|{\rm Tr}_{\Gamma,s}\left[N\left(\lambda-\widetilde{B}_T\right)^{-p}\right]-{\rm Tr}_{\Gamma,s}\left[N^{C^{\bullet}({W}^u,\widetilde{F})}\left(\lambda-\widetilde{B}_0^{C^{\bullet}({W}^u,\widetilde{F})}\right)^{-p}\right]\right|\\
\leq {C\over{\sqrt{T}}}(1+|\lambda|)^{p+1}.
\end{multline}
\end{thm}

\begin{proof}

Let $\{f_\alpha\}$ be an orthonormal basis of $TS$ and $\{f^\alpha\}$ be its dual. Let $\{e_i\}$ be an orthonormal basis of $TX$ and $\widetilde{e}_i$ be its lift to $T\widetilde{X}$. Let ${^1}\nabla^{\Lambda^\bullet(T^*S)\widehat{\otimes}\Lambda^\bullet (T^*\widetilde{X})}$ be the connection along the fibers $\widetilde{X}$ on $\Lambda^\bullet (T^*S)\widehat{\otimes}\Lambda^\bullet (T^*\widetilde{X})$ defined by (cf. \cite[(3.41)]{BG})
\begin{align}\nonumber
{^1}\nabla^{\Lambda^\bullet(T^*S)\widehat{\otimes}\Lambda^\bullet (T^*\widetilde{X})}=\nabla^{\Lambda^\bullet(T^*S)\widehat{\otimes}\Lambda^\bullet (T^*\widetilde{X})}+{1\over 2}\left\langle \widetilde{S}\widetilde{e}_i, {f}_\alpha^H\right\rangle \sqrt{2}c(\widetilde{e}_i)f^{\alpha}+{1\over 2}\left\langle \widetilde{S} {f}_{\alpha}^H, {f}_\beta^H\right\rangle f^\alpha f^\beta. 
\end{align}
Put (cf. \cite[(3.44)]{BG})
\begin{align}\nonumber
{^1}\nabla_{t}^{\Lambda^\bullet(T^S)\widehat{\otimes}\Lambda^\bullet (T^*\widetilde{X})\widehat{\otimes}\widetilde{F},u}=\psi_t ^{-1}{^1}\nabla^{\Lambda^\bullet(T^S)\widehat{\otimes}\Lambda^\bullet (T^*\widetilde{X})\widehat{\otimes}\widetilde{F},u}\psi_t.
\end{align}

Then as \cite[(11.53)]{BG}, we have
\begin{multline}\label{ad6.4}
\widetilde{B}_T=-{1\over 2}\widehat{c}(e_i)\left(^1 \nabla_{1/2,e_i}^{\Lambda^\bullet (T^\bullet S)\widehat{\otimes}\Lambda^\bullet (T^\bullet \widetilde{X})\widehat{\otimes}\widetilde{F},u}-T\left\langle \nabla\widetilde{f},e_i\right\rangle\right)-{T\over 2}c\left(\nabla \widetilde{f}\right)\\
+{1\over 4}c(e_i)\omega\left(\nabla^{\widetilde{F}},g^{\widetilde{F}}\right)(e_i)+{1\over 2}f^\alpha \omega\left(\nabla^{\widetilde{F}},g^{\widetilde{F}}\right)\left(f^H_\alpha\right).
\end{multline}
From (\ref{ad6.4}), we get
\begin{align}\label{ad6.5}
\overline{B}_{T}=-{1\over 2}\widehat{c}(e_i)^1 \nabla_{1/2,e_i}^{\Lambda^\bullet (T^\bullet S)\widehat{\otimes}\Lambda^\bullet (T^\bullet \widetilde{X})\widehat{\otimes}\widetilde{F},u}-{T\over 2}c\left(\nabla \widetilde{f}\right)+{1\over 2}f^{\alpha}\omega\left(\nabla^{\widetilde{F}},g^{\widetilde{F}}\right)\left(f^H_\alpha\right),
\end{align}
By (\ref{ad6.5}), we find that the $\overline{B}_T^{(>0)}$ does not depend on $T$. We will write $\overline{B}^{(>0)}$ instead of $\overline{B}_T^{(>0)}$.

Clearly, in (\ref{ad6.3}), we can replace $\widetilde{B}_T$ by $\overline{B}_T$. By the simplifying assumptions made in \cite[Section 9.1]{BG}, on the support of $\mu$,
\begin{align}\label{ad6.6}
\overline{B}_T=\overline{\mathcal{D}}_T^{\widetilde{I}\widehat{\otimes}\widetilde{F}|_{\widetilde{\bf B}}}.
\end{align}
As in (\ref{ad5.10}),
\begin{align}\label{ad6.7}
{\rm ker}\overline{\mathcal{D}}_T^{\widetilde{I}\widehat{\otimes}\widetilde{F}|_{\widetilde{\bf B}}}=\overline{\mathfrak{f}}_T\otimes \widetilde{F}|_{\widetilde{\bf B}}.
\end{align}

Moreover,
\begin{align}\label{ad6.8}
\left(\lambda-\overline{B}_T\right)^{-1}=\left(\lambda-\overline{B}_T^{(0)}\right)^{-1}+\left(\lambda-\overline{B}_T^{(0)}\right)^{-1} \left(\overline{B}^{>0}\right)\left(\lambda-\overline{B}_T^{(0)}\right)^{-1}+\cdots,
\end{align}
and the expansion in (\ref{ad6.8}) only contains a finite number of terms. Set
\begin{align}
\overline{P}^{(1,+\infty)}_T=1-\overline{P}^{[0,1]}_T.
\end{align}
If $\lambda\in U_T$, put
\begin{align}\nonumber
L_{T,1}=\overline{P}^{[0,1]}_{T}\left(\lambda-\overline{B}_T^{(0)}\right)^{-1}\overline{P}^{[0,1]}_T,\ \ L_{T,2}=\overline{P}^{[0,1]}_{T}\left(\lambda-\overline{B}_T^{(0)}\right)^{-1}\overline{P}^{(1,+\infty)}_T,
\end{align}
\begin{align}
L_{T,3}=\overline{P}^{(1,+\infty)}_{T}\left(\lambda-\overline{B}_T^{(0)}\right)^{-1}\overline{P}^{[0,1]}_T,\ \ L_{T,4}=\overline{P}^{(1,+\infty)}_{T}\left(\lambda-\overline{B}_T^{(0)}\right)^{-1}\overline{P}^{(1,+\infty)}_T.
\end{align}

We still use the notation in (\ref{ad5.90}) to define the norms $\|\ \|_{\Gamma,p}$. By proceeding as in Theorem \ref{Bt7}, we can define $m_T(\lambda)\in {\rm End}(C^\bullet (\widetilde{W}^u,\widetilde{F}))$ such that for $T\geq 0$ large enough,
\begin{align}
L_{T,1}=\left(m_T(\lambda)\lambda\right)^{-1},
\end{align}
and moreover by (\ref{B44}), if $c_1>0$ is small enough, if $\lambda\in U_T$,
\begin{align}
\left\|m_T^{-1}(\lambda)-1\right\|_{\infty}\leq {C\over{\sqrt{T}}}(1+|\lambda|).
\end{align}
By using (\ref{ad6.6}), (\ref{ad6.7}) and by proceeding as in the proof of Theorem \ref{Bt9}, we find that for $2\leq j\leq 4$,
\begin{align}\label{ad6.13}
\|L_{T,j}\|_{\Gamma,p-1}\leq C,\ \ \|L_{T,j}\|_{\infty}\leq {C\over \sqrt{T}}.
\end{align}
From (\ref{ad6.8})-(\ref{ad6.13}) we find that to establish (\ref{ad6.3}), in (\ref{ad6.8}), we may as well replace $(\lambda-\overline{B}_T^{(0)})^{-1}$ by $\overline{P}^{[0,1]}_{T}/\lambda$.

Let $\overline{p}_T$ be the orthogonal projection operator form $\Omega^{\bullet}_{(2)}(\widetilde{X},\widetilde{F}|_{\widetilde{X}})$ on ${\rm Im}(\widetilde{J}_T)\subset \Omega^{\bullet}_{(2)}(\widetilde{X},\widetilde{F}|_{\widetilde{X}})$. We claim that 
\begin{align}\label{ad6.14}
\left\|\left(\overline{P}^{[0,1]}_T-1\right)\overline{p}_T\right\|_{\infty}=\mathcal{O}\left(e^{-cT}\right).
\end{align}
In fact, since $\widetilde{J}_T:C^\bullet(\widetilde{W}^u,\widetilde{F})\to \Omega^{\bullet}_{(2)}(\widetilde{X},\widetilde{F}|_{\widetilde{X}})$ is an isometric embedding, (\ref{ad6.14}) follows from that $\overline{e}_T-\widetilde{J}_T=\mathcal{O}(e^{-cT})$ as $T\to +\infty$.

Recall that $\|\ \|_{\Gamma.2}$ is the $\Gamma$-Hilbert-Schmidt norm. Then
\begin{align}\label{ad6.15}
\left\|\overline{P}^{[0,1]}_T-\overline{p}_T\right\|^2_{\Gamma,2}=\left\|\overline{P}^{[0,1]}_T\right\|^2_{\Gamma,2}+\left\|\overline{p}_T\right\|^2_{\Gamma,2}-2{\rm Re}{\rm Tr}_{\Gamma}\left[\overline{P}^{[0,1]}_{T}\overline{p}_T\right].
\end{align}
Since the $\Gamma$-dimension of the images of $\overline{P}^{[0,1]}_T$ and $\overline{p}_T$ are both equal to ${\rm dim}(F|_{\bf B})$, using (\ref{ad6.14}), (\ref{ad6.15}), we get
\begin{align}\label{ad6.16}
\left\|\overline{P}^{[0,1]}_T-\overline{p}_T\right\|^2_{\Gamma,2}=\mathcal{O}\left(e^{-cT}\right).
\end{align}
Let $Q_T$ be the orthogonal projection operator from $\Omega^\bullet_{(2)}(\widetilde{X},\widetilde{F}|_{\widetilde{X}})$ on $\widetilde{F}^{[0,1]}_T+{\rm Im}(\widetilde{J}_T)$. Then 
\begin{multline}\label{ad6.17}
\left\|\overline{P}^{[0,1]}_T-\overline{p}_T\right\|_{\Gamma,1}=\left\|\left(\overline{P}^{[0,1]}_T-\overline{p}_T\right)Q_T\right\|_{\Gamma,1}\\
\leq \left\|\overline{P}^{[0,1]}_T-\overline{p}_T\right\|_{\Gamma,2}\left\|Q_T\right\|_{\Gamma,1}\leq 2{\rm dim}\left(F|_{\bf B}\right)\left\|\overline{P}^{[0,1]}_T-\overline{p}_T\right\|_{\Gamma,2}. 
\end{multline}
From (\ref{ad6.16}), (\ref{ad6.17}), we get
\begin{align}
\left\|\overline{P}^{[0,1]}_T-\overline{p}_T\right\|_{\Gamma,1}=\mathcal{O}\left(e^{-cT}\right).
\end{align}
By the above, it follows that to establish (\ref{ad6.3}), we may as well replace in (\ref{ad6.8}) $(\lambda-\overline{B}_T^{(0)})^{-1}$ by $\overline{p}_T/\lambda$.

Let $r:\Lambda(T^*\widetilde{X}|_{\widetilde{\bf B}})\to \Lambda^{\rm max}(T^*\widetilde{X}^u|_{\widetilde{\bf B}})$ the obvious orthogonal projection operator. Then using (\ref{ad5.10}), (\ref{ad5.124}), we find easily that if $s\in\Omega^{\bullet}_{(2)}(\widetilde{X},\widetilde{F}|_{\widetilde{X}})$,
\begin{align}\label{ed6.19}
\overline{p}_T s(Z)
={\mu(Z)\over \alpha_T}\exp\left(-T|Z|^2/2\right)r\int_{T\widetilde{X}|_{\widetilde{\bf B}}}\mu(Z')\exp\left(-T|Z'|^2/2\right)s(Z')dv_{T\widetilde{X}}(Z').
\end{align}

Using (\ref{ad6.6}), (\ref{ed6.19}), we get
\begin{align}\label{ed6.20}
\overline{p}_T \overline{B}_T \overline{p}_T=\overline{p}_T \overline{\mathcal{D}}^{\widetilde{I}\widehat{\otimes}\widetilde{F}|_{\widetilde{\bf B}}}_T\overline{p}_T. 
\end{align}

By (\ref{ad5.124}), since here $d\widetilde{\mathcal{F}}=0$, we get
\begin{align}
\overline{\mathcal{D}}^{\widetilde{I}\widehat{\otimes}\widetilde{F}|_{\widetilde{\bf B}},(>0)}_T={1\over 2}\left(\omega\left(\widetilde{F}|_{\widetilde{\bf B}},g^{\widetilde{F}|_{\widetilde{\bf B}}}\right)-\widehat{c}\left(R^{T\widetilde{X}|_{\widetilde{\bf B}}}Z\right)\right).
\end{align}
Since $c(R^{T\widetilde{X}|_{\widetilde{\bf B}}}Z)$ is an odd operator, we get
\begin{align}\label{ed6.22}
rc\left(R^{T\widetilde{X}|_{\widetilde{\bf B}}}\right)r=0.
\end{align}
By (\ref{ed6.19}), (\ref{ed6.20}), (\ref{ed6.22}), we obtain,
\begin{align}
\overline{p}\overline{B}^{(>0)}\overline{p}_T={1\over 2}\omega\left(\widetilde{F}|_{\widetilde{\bf B}},g^{\widetilde{F}|_{\widetilde{\bf B}}}\right).
\end{align}

Using the above, we get (\ref{ad6.3}). The proof is completed.

\end{proof}

\begin{thm}\label{t2}
For $c_1\in (0,1]$ small enough, given any integer $p\geq {\rm dim}X+1$, there exists $C>0$ such that for $T\geq 1$ large enough, and $\lambda\in U_T$,
\begin{align}\label{df6.24}
\left\|\left(\lambda-\widetilde{B}_T\right)^{-1}\right\|_{\Gamma,p}\leq C(1+|\lambda|)^{p}.
\end{align}
\end{thm}

\begin{proof}
Again, we can replace in (\ref{df6.24}) $\widetilde{B}_T$ by $\overline{B}_T$. First we claim that (\ref{df6.24}) holds for $\overline{B}_T^{(0)}$. Using (\ref{ad6.6}), the proof is the same as the proofs of Theorems \ref{Bt7} and \ref{Bt9}. To get (\ref{df6.24}), we use (\ref{ad6.8}) and the fact that $\overline{B}^{(>0)}$ is of order $0$. The proof is completed. 
\end{proof}

Now we prove Theorem \ref{t3}.

Take $p\in\mathbb{\bf N}$. Let $k_p(\lambda)$ be the unique holomorphic function on $\mathbb{\bf C}\setminus\mathbb{\bf R}$ such that

--As $\lambda\to \pm i\infty$, $k_p(\lambda)\to 0$.

--The following identity holds,
\begin{align}
{{k_p^{(p-1)}(\lambda)}\over{(p-1)!}}=h'(\lambda).
\end{align}
Clearly, if $\lambda\in \Delta$,
\begin{align}\label{1}
|{\rm Re}(\lambda)|\leq {1\over 2}|{\rm Im}(\lambda)|.
\end{align}
Using (\ref{1}), we find that there exist $C>0$, $C'>0$ such that if $\lambda\in\Delta$,
\begin{align}\label{2}
\left|k_p(\sqrt{t}\lambda)\right|\leq C\exp\left(-C't|\lambda|^2\right).
\end{align}

By the argument of \cite[Theorem 11.7]{BG}, we find that for $T\geq 0$ large enough,
\begin{align}\label{ap6.28}
\left|{\rm Sp}\left(\widetilde{B}^{(0)}_{T}\right)\right|\subset \left(T\over \pi\right)^{1/2} e^{-T}\left[0,d\right]\cup [1,+\infty).
\end{align}

Let $\Delta=\Delta_{+}\cup \Delta_{-}$ be a contour defined by 
\begin{align}\nonumber
\Delta_{+}=\left\{z=x+iy|x=\pm{1\over 4}, {1\over 2}\leq y<+\infty\right\}\cup \left\{z=x+iy|y={1\over 2}, -{1\over 4}\leq x\leq {1\over 4}\right\},
\end{align}
\begin{align}
\Delta_{-}=\left\{z=x+iy|x=\pm{1\over 4}, -\infty< y\leq -{1\over 2}\right\}\cup\left\{z=x+iy|y=-{1\over 2}, -{1\over 4}\leq x\leq {1\over 4}\right\}.
\end{align}
Let $\delta$ be the unit circle in ${\bf C}$. Then for $T\geq 0$ large enough, for $t>0$, we have (cf. \cite[(12.47)]{BG})
\begin{align}\label{3}
h'\left(\widetilde{D}_{t,T}\right)=\psi_{t}^{-1}{1\over{2i\pi}}\int_{\Delta}{{h'(\sqrt{t}\lambda)}\over{\lambda-\widetilde{B}_{T}}}d\lambda \psi_t + {1\over{2i\pi}}\int_{\delta/4}{{h'(\lambda)}\over{\lambda-\widetilde{D}_{t,T}}}d\lambda .
\end{align}
Using Theorems \ref{t1} and \ref{t2}, (\ref{2}) and (\ref{3}), there exists $C>0$ such that for $t\in [\varepsilon,A]$, $T\geq 1$, we have
\begin{align}\label{ap6.31}
\left|\psi_{t}^{-1}{1\over{2i\pi}}\int_{\Delta}{{h'(\sqrt{t}\lambda)}\over{\lambda-\widetilde{B}_{T}}}d\lambda \psi_t \right|\leq {C\over {\sqrt{T}}}.
\end{align}
For $T>0$ large enough, proceedings as in the proof of \cite[Theorem 11.19]{BG}, one has
\begin{align}\label{ap6.32}
\lim_{t\to 0^+}{1\over{2i\pi}}\int_{\delta/4}{{h'(\lambda)}\over{\lambda-\widetilde{D}_{t,T}}}d\lambda=\widetilde{\chi}'^{-}(F).
\end{align}
Then by (\ref{ap6.28}) and (\ref{ap6.32}), for $T>0$ large enough, $t\in [\epsilon, A]$, one has
\begin{align}\label{ap6.33}
\left|\int_{\delta/4}{{h'(\lambda)}\over{\lambda-\widetilde{D}_{t,T}}}d\lambda-\widetilde{\chi}'^{-}(F)\right|\leq {C\over{\sqrt{T}}}.
\end{align}
Then by (\ref{ap6.31}) and (\ref{ap6.33}), we get the $L^2$ case of \cite[(12.48)]{BG}, i.e., given $\varepsilon$, $A$ with $0\leq \varepsilon\leq A\leq +\infty$, there exists $C>0$ such that if $t\in [\varepsilon,A]$, $T\geq 1$,
\begin{align}
\left|{\rm Tr}_{\Gamma,s}\left[Nh'\left(\widetilde{D}_{t,T}\right)\right]-\widetilde{\chi}'^{-}(F)\right|\leq {C\over{\sqrt{T}}}.
\end{align}
Then we get Theorem \ref{t3}.

In the general case, as in \cite[(12.49)]{BG}, (\ref{ad6.6}) is replaced by 
\begin{align}
\overline{B}^{2}_T=\overline{\mathcal{D}}_T^{\widetilde{I}\widehat{\otimes}\widetilde{F}|_{\widetilde{\bf B}},2}.
\end{align}    
Let $\Gamma''$ be a contour defined by
\begin{align}\nonumber
\Gamma''=\left\{z=x+iy|x=-{1\over 4}, -{1\over 2}\leq y\leq {1\over 2}\right\}\cup \left\{z=x+iy|y=\pm {1\over 2},-\infty<x\leq -{1\over 4}\right\}.
\end{align}
Then (\ref{ap6.31}) is replaced by
\begin{align}\label{ap6.36}
\psi^{-1}_t {1\over{2\pi i}}\int_{\Gamma''}{{(1+2t\lambda)e^{t\lambda}}\over{\lambda-\widetilde{B}^2_T}}d\lambda\psi_t.
\end{align} 
As in \cite{BG}, proceedings as in Appendix $B$, we get the corresponding estimate for (\ref{ap6.36}). The proof is completed. 

\appendix

\section {$L^2$ analogue of some results of \cite{BLe} in the current case}

\setcounter{equation}{0}

In this appendix, we generalize some results in \cite{BLe} to the current case.

For $h\in o(x)\otimes \widetilde{F}|_{x}$, recall that $\widetilde{J}_T(h)$ has been defined in (\ref{ad5.124}). Let 
$$\widetilde{E}_T=\widetilde{J}_T\left(L^2 \left(\widetilde{F}|_{\widetilde{B}}\right)\right)$$
be the image of $\widetilde{J}_T$. Since $\widetilde{J}_T$ is an isometry, $\widetilde{E}_T\subset \Omega_{(2)}^{\bullet}(\widetilde{X},\widetilde{F}|_{\widetilde{X}})$ is closed.

Let $\widetilde{E}^{\bot}_T$ denote the orthogonal complement of $\widetilde{E}_T$ in $\Omega_{(2)}^{\bullet}(\widetilde{X},\widetilde{F}|_{\widetilde{X}})$, that is
$$\Omega_{(2)}^{\bullet}(\widetilde{X},\widetilde{F}|_{\widetilde{X}})=\widetilde{E}_T\oplus \widetilde{E}^{\bot}_T.$$

Let $\widetilde{P}_T$ (resp. $\widetilde{P}^{\bot}_{T}$) denote the orthogonal projection from $\Omega_{(2)}^{\bullet}(\widetilde{X},\widetilde{F}|_{\widetilde{X}})$ onto $\widetilde{E}_T$ (resp. $\widetilde{E}^{\bot}_{T}$).

Recall that we have the decompositions 
$$\widetilde{A}'=d^{\widetilde{X}}+\nabla^{\Omega_{(2)}^{\bullet}(\widetilde{X},\widetilde{F}|_{\widetilde{X}})}+i_{T^H},\ \widetilde{A}''=d^{\widetilde{X},*}+\nabla^{\Omega^{\bullet}(\widetilde{X},\widetilde{F}|_{\widetilde{X}}),*}-T^H\wedge.$$
For any $T\geq 0$, following \cite{W}, set
$$d^{\widetilde{X}}_{T}=e^{-T\widetilde{f}}d^{\widetilde{X}}e^{T\widetilde{f}},\ \delta^{\widetilde{X},*}_T=e^{T\widetilde{f}}d^{\widetilde{X},*}e^{-T\widetilde{f}}.$$
Then $\delta^{\widetilde{X},*}_T$ is the formal adjoint of $d^{\widetilde{X}}_{T}$ with respect to the usual inner product on $\Omega_{(2)}^{\bullet}(\widetilde{X},\widetilde{F}|_{\widetilde{X}})$.

Set
$$\widetilde{D}'_T=d^{\widetilde{X}}_{T}+\delta^{\widetilde{X},*}_T,\ \widetilde{D}'^2_T=\left(d^{\widetilde{X}}_{T}+\delta^{\widetilde{X},*}_T\right)^2=d^{\widetilde{X}}_{T}\delta^{\widetilde{X},*}_T+\delta^{\widetilde{X},*}_Td^{\widetilde{X}}_{T}.$$
Then $\widetilde{D}'^2_T$ preserves the $\mathbb{Z}$-grading of $\Omega_{(2)}^{\bullet}(\widetilde{X},\widetilde{F}|_{\widetilde{X}})$.

Following Bismut-Lebeau \cite[Section 9]{BLe}, we define
\begin{align}\nonumber
\widetilde{D}'_{T,1}=\widetilde{P}_T\widetilde{D}'_{T}\widetilde{P}_T,\ \ \widetilde{D}'_{T,2}=\widetilde{P}_{T}\widetilde{D}'_T \widetilde{P}^{\bot}_T,
\end{align}
\begin{align}
\widetilde{D}'_{T,3}=\widetilde{P}^{\bot}_T\widetilde{D}'_T \widetilde{P}_T,\ \ \widetilde{D}'_{T,4}=\widetilde{P}^{\bot}_{T}\widetilde{D}'_T\widetilde{P}^{\bot}_T. 
\end{align}
We then write the operator $\widetilde{D}'_T$ in matrix form
\begin{align}
\widetilde{D}'_T=\left(\begin{matrix}\widetilde{D}'_{T,1}&\widetilde{D}'_{T,2}\\ \widetilde{D}'_{T,3}&\widetilde{D}'_{T,4}\end{matrix}\right).
\end{align}

We recall the following elliptic estimates needed in this paper. 
\begin{lemma}(\cite[Lemma 1.4]{Shu})\label{r2lA1}
Let $A:C^\infty(\widetilde{X},\widetilde{F}|_{\widetilde{X}})\to C^\infty(\widetilde{X},\widetilde{F}|_{\widetilde{X}})$ be a $C^\infty$-bounded uniformly elliptic differential operator of order $m$. For any $i,j\geq 0$, there exists a constant $C$ such that for any $s\in C_0^{\infty}(\widetilde{X},\widetilde{F}|_{\widetilde{X}})$,
$$\|s\|_{i+m}\leq C\left(\|As\|_i +\|s\|_j\right).$$
\end{lemma}

Let $\mathbb{H}^{1}(\widetilde{X},\widetilde{F}|_{\widetilde{X}})$ denote the first Sobolev space with respect to a (fixed, $\Gamma$-invariant) first Sobolev norm on $\Omega^*(\widetilde{X},\widetilde{F}|_{\widetilde{X}})$. Using Lemma \ref{r2lA1}, by proceeding as in \cite[Section 9 b)]{BLe}, one deduces that

(i) The following identity holds (cf. \cite[(6.15)]{Zhang}):
\begin{align}\label{r2A3}
\widetilde{D}'_{T,1}=0.
\end{align}

(ii) There exists $C>0$ such that for any $T\geq 1$, any $s\in \widetilde{E}^{\bot}_T\cap \mathbb{H}^1 (\widetilde{X},\widetilde{F}|_{\widetilde{X}})$, $s'\in \widetilde{E}_T\cap \mathbb{H}^{1}(\widetilde{X},\widetilde{F}|_{\widetilde{X}})$, then
\begin{align}\label{r2A4}
\left\|\widetilde{D}'_{T,2}s\right\|_0 \leq C \left({{\|s\|_1}\over{\sqrt{T}}}+\|s\|_0\right)\ {\rm and}\ \left\|\widetilde{D}'_{T,3}\right\|_0 \leq C\left({{\|s'\|_1}\over{\sqrt{T}}}+\|s'\|_0\right). 
\end{align}

(iii) There exist $T_0>0$ and $c>0$ such that for any $s\in \widetilde{E}^{\bot}_T\cap \mathbb{H}^1(\widetilde{X},\widetilde{F}|_{\widetilde{X}})$, then
\begin{align}\label{r2A5}
\left\|\widetilde{D}'_{T,4}s\right\|_0 \geq c\left(\|s\|_1 +\sqrt{T}\|s\|_0\right). 
\end{align}

(iv) For any $T\geq T_0$, $\lambda\in \bf{C}$, $|\lambda|\leq {{c\sqrt{T}}\over{2}}$, $s\in \widetilde{E}_T^{\bot}$, then
\begin{align}\label{r2A6}
\left\|\left(\lambda-\widetilde{D}'_{T,4}\right)^{-1} s\right\|_0 \leq {C\over{\sqrt{T}}}\|s\|_0,\ \ \left\|\left(\lambda-\widetilde{D}'_{T,4}\right)^{-1}s\right\|_{1}\leq \|s\|_0.
\end{align}.

Following \cite[(9.113)]{BLe}, for $T\geq 1$, set
\begin{align}\label{r2A8}
U_{T}=\left\{\lambda\in\mathbb{C}: 1\leq |\lambda|\leq {{c\sqrt{T}}\over{4}}\right\}. 
\end{align}

Note that one can also extend \cite[Definitions 9.17 and 9.22]{BLe} to the current case, and they still have the similar properties.

\begin{definition}\label{Ade5}(Compare with \cite[Definition 9.17]{BLe})
If $H$, $H'$ are $\Gamma$-Hilbert modules, set
$$\mathcal{L}_{\Gamma,p}=\left\{A\in\mathcal{L}(H,H');\ {\rm Tr}_{\Gamma}\left[(A^*A)^{p/2}\right]<+\infty\right\}.$$
If $A\in \mathcal{L}_{\Gamma,p}(H,H')$, set
$$\|A\|_{\Gamma,p}=\left\{{\rm Tr}_{\Gamma}\left[(A^*A)^{p/2}\right]\right\}^{1/p}.$$
\end{definition}

From (\ref{r2A3}) to (\ref{r2A8}), one can show that there exists $T_0\geq 1$ such that for any $T\geq T_0$, $\lambda\in U_T$, $\lambda-\widetilde{D}'_T$ is invertible.

\begin{proposition}\label{r2PA2}
If $p\geq {\rm dim}X+1$, there exists $C>0$ such that for $T\geq T_0$, $\lambda\in\bf{C}$, $\lambda\in {{c\sqrt{T}}\over{2}}$, then
\begin{align}\nonumber
\left\| \left(\lambda-\widetilde{D}'_{T,4}\right)^{-1}\right\|_{\infty}\leq {C\over{\sqrt{T}}},
\end{align}
\begin{align}\nonumber
\left\|\left(\lambda-\widetilde{D}'_{T,4}\right)^{-1}\right\|_{\Gamma,p}\leq C,
\end{align}
\begin{align}\label{r2Aa8}
\left\|\widetilde{D}'_{T,2}\left(\lambda-\widetilde{D}'_{T,4}\right)^{-1}\right\|_{\infty}\leq {C\over{\sqrt{T}}}.
\end{align}
\end{proposition}
\begin{proof}
The first line of (\ref{r2Aa8}) is from (\ref{r2A6}). Also
\begin{align}\nonumber
\left\|\left(\lambda-\widetilde{D}'_{T,4}\right)^{-1}\right\|_{\Gamma,p}\leq \left\|\left(D^{\widetilde{X}}+\sqrt{-1}\right)^{-1}\right\|_{\Gamma,p}\left\|\left(D^{\widetilde{X}}+\sqrt{-1}\right)\left(\lambda-\widetilde{D}'_{T,4}\right)^{-1}\right\|_{\infty}.
\end{align}
Since $D^{\widetilde{X}}$ is elliptic of order one, when $p\geq {\rm dim}X+1$, $\|(D^{\widetilde{X}}+\sqrt{-1})^{-1}\|_{\Gamma,p}<+\infty$. 
Also by (\ref{r2A6}), for $T\geq T_0$,
\begin{align}
\left\|\left(D^{\widetilde{X}}+\sqrt{-1}\right)\left(\lambda-\widetilde{D}'_{T,4}\right)^{-1}\right\|_{\infty}\leq C.
\end{align}
The second line in (\ref{r2Aa8}) follows. Using (\ref{r2A4}) and (\ref{r2A6}), we get the third line in (\ref{r2Aa8}).
\end{proof}

\begin{definition}
For $T\geq T_0$, $\lambda\in U_T$, let $M_T(\lambda)$ be the linear map from $\widetilde{E}_T\cap \mathbb{H}^1 (\widetilde{X},\widetilde{F}|_{\widetilde{X}})$ into $\widetilde{E}^{\bot}_T$
\begin{align}
M_T(\lambda)=\lambda-\widetilde{D}'_{T,2}\left(\lambda-\widetilde{D}'_{T,4}\right)^{-1}\widetilde{D}'_{T,3}. 
\end{align}
\end{definition}

\begin{theorem}\label{Bt7}
$M_T(\lambda)$ is invertible and for any integer $p\geq {\rm dim}X+1$, there exist $C>0$ such that for $T\geq T_0$, $\lambda\in U_T$,
\begin{align}\nonumber
\left\|M_T^{-1}(\lambda)\right\|_{\infty}\leq C,\ \left\|\widetilde{D}'_{T,3}M_T^{-1}(\lambda)\right\|_{\infty}\leq C,\ \left\| M^{-1}_T(\lambda)\right\|_{\Gamma,p}\leq C(1+|\lambda|),
\end{align}
\begin{align}\label{ada36}
\left\|\widetilde{J}_T^{-1}\left(M_T^{-1}(\lambda)\right)^p \widetilde{J}_T-\lambda^{-p}\right\|_{\Gamma,1}\leq {C\over{\sqrt{T}}}(1+|\lambda|)^{p+1}.
\end{align}
\end{theorem}

\begin{proof}

For $\lambda\in U_T$, set
\begin{align}\label{B38}
m_T(\lambda)=1-{\lambda}^{-1}\widetilde{D}'_{T,2}\left(\lambda-\widetilde{D}'_{T,4}\right)^{-1}\widetilde{D}'_{T,3}.
\end{align}
Clearly
\begin{align}\label{B39}
M_T(\lambda)=\lambda m_T(\lambda).
\end{align}

For $\lambda\in U_T$, by (\ref{r2A4}) and (\ref{r2A6}), there exits $C_1>0$ such that for $T\geq T_0$, 
\begin{align}\nonumber
\left\|\widetilde{D}'_{T,2}\left(\lambda-\widetilde{D}'_{T,4}\right)^{-1}\right\|_{\infty}\leq {C_1\over{\sqrt{T}}}.
\end{align}
By \cite[(6.16)]{Zhang}, there exists $C_2>0$ such that  
\begin{align}\label{r2A13}
\left\|\lambda^{-1}{\widetilde{D}'_{T,3}}\right\|_{\infty}\leq C_2.
\end{align}
Then there exists $C>0$ such that
\begin{align}\label{B43}
\left\|{\lambda}^{-1}\widetilde{D}'_{T,2}\left(\lambda-\widetilde{D}'_{T,4}\right)^{-1}\widetilde{D}'_{T,3}\right\|_{\infty}\leq {C\over\sqrt{T}}.
\end{align}

From (\ref{B38})-(\ref{B43}), it is clear that if $1/\sqrt{T}$ and $|\lambda|/\sqrt{T}$ are small enough, the operator $m_T (\lambda)$ is invertible, and moreover for $T\geq 1$
\begin{align}\label{B44}
\left\|m_T^{-1}(\lambda)-1\right\|_{\infty}\leq {C\over{\sqrt{T}}}(1+|\lambda|).
\end{align}
In particular
\begin{align}\label{B45}
\left\|m_T^{-1}(\lambda)\right\|_{\infty}\leq C.
\end{align}
By (\ref{B39}), we get
\begin{align}\label{B46}
M_T^{-1}(\lambda)=\lambda^{-1} m_T^{-1}(\lambda). 
\end{align}
By (\ref{r2A13}), (\ref{B45}) and (\ref{B46}), we get the second inequality in (\ref{ada36}).

From (\ref{B45}) and (\ref{B46}), we also get
\begin{align}\label{r2A18}
\left\|M^{-1}_T(\lambda)\right\|_{\Gamma,p}\leq C(1+|\lambda|),
\end{align}
which get the third inequality in (\ref{ada36}).

Finally using (\ref{B44})-(\ref{r2A18}), we get the last inequality in (\ref{ada36}).

\end{proof}

If $B\in\mathcal{L}(\Omega^{\bullet}_{(2)}(\widetilde{X},\widetilde{F}|_{\widetilde{X}}))$, for any $T\geq 1$, we write $B$ as a matrix with respect to the splitting $\Omega^{\bullet}_{(2)}(\widetilde{X},\widetilde{F}|_{\widetilde{X}})=\widetilde{E}_T\oplus \widetilde{E}^{\bot}_{T}$ in the form
$$B=\left(\begin{matrix}B_{1}&B_{2}\\ B_{3}&B_{4}\end{matrix}\right).$$

\begin{definition}(Compare with \cite[Definition 9.22]{BLe})
If $B\in \mathcal{L}(\Omega^{\bullet}_{(2)}(\widetilde{X},\widetilde{F}|_{\widetilde{X}}))$, $C\in\mathcal{L}(L^2 (\widetilde{F}|_{\widetilde{B}}))$, set
\begin{align}
d_{\Gamma}(B,C)=\sum_{j=2}^{4}\|B_{j}\|_{\Gamma,1}+\|\widetilde{J}_{T}^{-1}B_{1}\widetilde{J}_{T}-C\|_{\Gamma,1}.
\end{align}
\end{definition}
Clearly if $B\in\mathcal{L}_{1}(\Omega^{\bullet}_{(2)}(\widetilde{X},\widetilde{F}|_{\widetilde{X}}))$, $C\in\mathcal{L}_{1}(L^2 (\widetilde{F}|_{\widetilde{B}}))$,
\begin{align}
|{\rm Tr}_{\Gamma}(B)-{\rm Tr}_{\Gamma}(C)|\leq d_{\Gamma}(B,C).
\end{align}

From the proof of \cite[Theorem 9.23]{BLe}, one can easily get the following analogue of \cite[Theorem 9.23]{BLe} in our case.
\begin{theorem}\label{Bt9}
There exists $T_{0}\geq 1$ such that for any $T\geq T_{0}$, $\lambda\in U_{T}$, $\lambda-\widetilde{D}'_{T}$ is invertible. For any integer $p\geq {\rm dim}X+2$, there exists $C>0$ such that if $T\geq T_{0}$, $\lambda\in U_{T}$, then
$$d_{\Gamma}\left(\left(\lambda-\widetilde{D}'_{T}\right)^{-p},\lambda^{-p}|_{L^2(\widetilde{F}|_{\widetilde{B}})}\right)\leq {C\over{\sqrt{T}}}(1+|\lambda|)^{p+1}.$$
\end{theorem}

\begin{proof}
Set 
\begin{align}
\widetilde{B}_T=\left(\lambda-\widetilde{D}'_T\right)^{-1}.
\end{align}
Then we find that
\begin{align}\nonumber
\widetilde{B}_{T,1}=M_T^{-1}(\lambda),\ \widetilde{B}_{T,2}=M_T^{-1}(\lambda)\widetilde{D}'_{T,2}\left(\lambda-\widetilde{D}'_{T,4}\right)^{-1},
\end{align}
\begin{align}
\widetilde{B}_{T,3}=\left(\lambda-\widetilde{D}'_{T,4}\right)^{-1}\widetilde{D}'_{T,3}M_T^{-1}(\lambda),\ \widetilde{B}_{T,4}=\left(\lambda-\widetilde{D}'_{T,4}\right)^{-1}\left(1+\widetilde{D}'_{T,3}\widetilde{B}_{T,2}\right).
\end{align}
If $\lambda\in U_T$, then $|\lambda|\leq c_1 \sqrt{T}$. Using Proposition \ref{r2PA2} and Theorem \ref{Bt7}, we find that if $p\geq 2{\rm dim}X+2$, $T\geq T_0$, $\lambda\in U_T$, for $2\leq j\leq 4$, then
\begin{align}\label{B53}
\left\|\widetilde{B}_{T,j}\right\|_{\Gamma,p-1}\leq C;\ \ \left\|\widetilde{B}_{T,j}\right\|_{\infty}\leq {C\over{\sqrt{T}}}. 
\end{align}
From Theorem \ref{Bt7} and from (\ref{B53}), we deduce that if $j_1,\cdots,j_p\in\{1,2,3,4\}$, if one of the $j_k's$ is not equal to $1$, then
\begin{align}\label{ada48}
\left\|\widetilde{B}_{T,j_1}\cdots \widetilde{B}_{T,j_p}\right\|_{\Gamma,1}\leq {C\over{\sqrt{T}}}(1+|\lambda|)^{p-1}.
\end{align}
Then follows from the fourth inequality in (\ref{ada36}) and (\ref{ada48}), we get the theorem. 
\end{proof}

\section{ $L^2$ analogue of \cite[Theorem 9.5]{B} in the current case}

\setcounter{equation}{0}

In this appendix, we prove the $L^2$ analogue of \cite[Theorem 9.5]{B} in the current case. We will use the notations of \cite{B} in this appendix otherwise indicated.

Although \cite[Theorem 9.1]{B} does not hold in $L^2$-case, but when $T$ large enough, the spectra gap also holds. In this paper, we only need the gap. This is the key point that we can get our theorem. So the proofs of the $L^2$-case are almost the same as in \cite{B}.

We also denote by $\widetilde{\pi}$ the canonical projection $T\widetilde{X}\to \widetilde{\bf B}$.

\begin{definition}(Compare with \cite[Definition 9.9]{B})
For $s\in S$, $\mu\in{\bf R}$, let $\widetilde{E}^\mu_s$ (resp. $\widetilde{{\bf E}}^\mu_s$, resp. $\widetilde{F}^\mu_s$) be the set of sections of $\Lambda(T^*\widetilde{X})\otimes \widetilde{F}$ over $\widetilde{X}_s$ (resp. of $\widetilde{\pi}^*((\Lambda(T^*\widetilde{X})\otimes \widetilde{F})|_{\widetilde{\bf{B}}})$ over $T\widetilde{X}|_{\widetilde{B}_s}$, resp. $\widetilde{F}$ over $\widetilde{B}_s$) which lie in the $\mu^{\rm th}$ Sobolev space, and let $\|\ \|_{\widetilde{E}^\mu_s}$ (resp. $\|\ \|_{\widetilde{\bf E}^\mu_s}$, resp. $\|\ \|_{\widetilde{F}^\mu_s}$) be the corresponding Sobolev norm. 
\end{definition}

Let $\varepsilon_0>$ be defined in \cite[Section 7.3]{B}. We take $\varepsilon\in (0,{\varepsilon_0\over 4}]$.

For $\mu\geq 0$, $T>0$, let $\widetilde{J}_T:\widetilde{F}^\mu\to \widetilde{\bf {E}}^\mu$ be the linear map defined by (\ref{ad5.124}). Let $\widetilde{{\bf E}}^\mu_T$ be the image of $\widetilde{F}^\mu$ in $\widetilde{\bf E}^{\mu}$ by $\widetilde{J}_T$. Then $\widetilde{J}_T$ is an isometric embedding of $\widetilde{F}^0$ into $\widetilde{\bf E}^0$.

Let $\widetilde{\bf E}^{0,\bot}_T$ be the orthogonal space to $\widetilde{\bf E}^0_T$ in $\widetilde{\bf E}^0$, let $\widetilde{p}_T$, $\widetilde{p}^{\bot}_T$ be the orthogonal projection operators from $\widetilde{\bf E}^0$ on $\widetilde{\bf E}^0_T$, $\widetilde{\bf E}^{0,\bot}_T$ respectively.

If $\sigma\in\widetilde{F}^\mu$, we consider $\widetilde{J}_T\sigma$ as an element of $\widetilde{E}^\mu$.

Then $\widetilde{J}_T$ is an isometric embedding from $\widetilde{F}^0$ into $\widetilde{E}^0$. Let $\widetilde{E}^\mu_T$ be the image of $\widetilde{F}^\mu$ in $\widetilde{E}^\mu$. Let $\widetilde{E}^{0,\bot}_T$ be the orthogonal bundle to $\widetilde{E}^0_T$ in $\widetilde{E}^0$.

For $\mu\geq 0$, set
\begin{align}
\widetilde{E}^{\mu,\bot}_T=\widetilde{E}^\mu\cap \widetilde{E}^{0,\bot}_T.
\end{align}
Let $\bar{p}_T$, $\bar{p}^{\bot}_T$ be the orthogonal projection operators from $\widetilde{E}^0$ on $\widetilde{E}^0_T$, $\widetilde{E}^{0,\bot}_T$.

Let $e_1,\cdots,e_{n}$ be a locally defined smooth orthogonal basis of $TX$. Denote by $\widetilde{e}_1, \cdots, \widetilde{e}_n$, the lifting to $T\widetilde{X}$.

If $U\in T\widetilde{\bf {B}}$, as in \cite[Section 8.2]{B}, one define $U^H\in T^H T\widetilde{X}$.

If $s\in\widetilde{E}^{0}$, put
\begin{align}
|s|_0=\|s\|_{\widetilde{E}^0},\ \ \langle s,s'\rangle_0=\langle s,s'\rangle_{\widetilde{E}^0}.
\end{align}

\begin{definition}
For $T\geq 1$, $s\in\widetilde{E}$, set
\begin{align}\label{A9}
|s|^2_{T,1}=\left|\bar{p}_T s\right|^2_0+T\left|\bar{p}^{\bot}_{T}s\right|_0^2+T^2 \left|\widehat{c}(\nabla\widetilde{f})\bar{p}^{\bot}_T s\right|^2_0+\sum_{1}^{n}\left|\nabla^{\Lambda(T^*\widetilde{X}){\otimes}\widetilde{F}}_{\widetilde{e}_i}\bar{p}^{\bot}_T s\right|^2_0.
\end{align}
\end{definition}

Then (\ref{A9}) defines a Hilbert norm on $\widetilde{E}^1$. We identify $\widetilde{E}^0$ to its antidual by $\langle\ ,\ \rangle_0$. Then we can identify $\widetilde{E}^{-1}$ to the antidual of $\widetilde{E}^1$. Let $|\ |_{T,-1}$ be the norm on $\widetilde{E}^{-1}$ associated to $|\ |_{T,1}$. Then we have the continuous dense embeddings with norms smaller than $1$,
\begin{align}
\widetilde{\bf E}^1\to \widetilde{\bf E}^0\to \widetilde{\bf E}^{-1}. 
\end{align}
For convenience, we introduce a metric $g^{TS}$ on $TS$. Then the definition of $|s|_0$, $|s|_{T,1}$ obviously extends to $\Lambda(T^*S)\widehat{\otimes}\widetilde{E}$.

As in \cite[Definition 2.12]{B}, we define 
\begin{align}\nonumber
\widetilde{B}_u^{\widetilde{M}}=\sqrt{u}D^{\widetilde{X}}+\nabla^{\Omega^{\bullet}_{(2)}(\widetilde{X},\widetilde{F}|_{\widetilde{X}})}-{{c\left(T^H\right)}\over{2\sqrt{u}}},
\end{align}
where $\nabla^{\Omega^{\bullet}_{(2)}(\widetilde{X},\widetilde{F}|_{\widetilde{X}})}$ is defined in (\ref{nn1}) and $T^H$ is defined in (\ref{ad1}).

Then for $u>0$, $T>0$, as \cite[(4.1)]{B}, we define
\begin{align}
{A}_{u,T}=\widetilde{B}^{\widetilde{M}}_{u^2}+T\widehat{c}(\nabla\widetilde{f}).
\end{align}

Let $f_1,\dots, f_m$ be a locally defined smooth basis of $TS$, and let $f^1,\dots,f^m$ be the corresponding dual basis of $T^*S$. Then as \cite[(7.35)]{B}, we define $\widetilde{A}_{u,T}$ as follows.

\begin{definition}
For $u>0$, $T>0$, set
\begin{align}\nonumber
\widetilde{A}_{u,T}
=\exp\left\{-f^\alpha{c\over u}\left(f_{\alpha}^{H,TX}\right)\right\}A_{u,T}\exp\left\{f^\alpha{c\over u}\left(f_{\alpha}^{H,TX}\right)\right\},
\end{align}
where $f_\alpha ^{H,TX}$ is defined as \cite[(7.32)]{B}.
\end{definition}

Put 
\begin{align}
\widetilde{A}_T=\widetilde{A}_{1,T}.
\end{align}
Let $\widetilde{A}^{(0)}_{T}$ (resp. $\widetilde{A}^{(>0)}$) be the piece of $\widetilde{A}_T$ which has degree $0$ (resp. positive degree) in $\Lambda(T^*S)$. Then
\begin{align}
\widetilde{A}_T=\widetilde{A}^{(0)}_T+\widetilde{A}^{(>0)}. 
\end{align}
As \cite[(9.41)]{B}, we also have
\begin{align}\label{apB8}
\widetilde{A}^{(0)}_{T}=D^{\widetilde{X}}+T\widehat{c}\left(\nabla\widetilde{f}\right).
\end{align}
Set
\begin{align}
\widetilde{R}_T=\left[\widetilde{A}^{(0)}_T,\widetilde{A}^{(>0)}\right]+\widetilde{A}^{(>0),2}.
\end{align}
Then $\widetilde{R}_T$ is a first order differential operator and moreover
\begin{align}
\widetilde{A}^2_T=\widetilde{A}^{(0),2}_T+\widetilde{R}_T.
\end{align}
Then by Lemma \ref{r2lA1}, proceeding as the proof of \cite[Theorem 9.14]{B}, we have
\begin{theorem}\label{thA5}
If $\varepsilon\in (0,\varepsilon_0/4]$ is small enough, there exist constants $C_1>0$, $C_2>0$, $C_3>0$ such that for $T\geq 1$, $s,s'\in\Lambda(T^*S)\widehat{\otimes}\widetilde{E}$,
\begin{align}\nonumber
\left|\widetilde{A}^{(0)}_T s\right|^2_0\geq C_1 |s|^2_{T,1}-C_2 |s|^2_0,
\end{align}
\begin{align}\nonumber
\left|\left\langle \widetilde{A}^{(0)}_T s,\widetilde{A}^{(0)}_T s'\right\rangle\right|\leq C|s|_{T,1} |s'|_{T,1},
\end{align}
\begin{align}\label{A17}
\left|\left\langle \widetilde{R}_T s,s'\right\rangle_0\right|\leq C_3 \left(|s|_{T,1}|s'|_0+|s|_0|s'|_{T,1}\right).
\end{align}
\end{theorem}

\begin{proof}
In the whole proof, $C$, $C'$, $\dots$ are positive constants, which may vary form line to line.

To establish the first inequality in (\ref{A17}), we may as well assume that $s\in \widetilde{E}$. If $s\in \widetilde{E}$, then
\begin{multline}\label{apB12}
\left|\widetilde{A}^{(0)}_T s\right|^2_0=\left|\overline{p}_T\widetilde{A}^{(0)}_T s\right|^2_0+\left|\overline{p}^{\bot}_T\widetilde{A}^{(0)}_Ts\right|^2_0\geq {1\over 2}\left|\overline{p}^{\bot}_{T}\widetilde{A}^{(0)}_T\overline{p}^{\bot}_T s\right|^2_0\\
-\left|\overline{p}_T\widetilde{A}^{(0)}_T\overline{p}^{\bot}_Ts\right|_0^2-\left|\overline{p}^{\bot}_T\widetilde{A}^{(0)}_T\overline{p}_Ts\right|^2_0,
\end{multline}
since $\overline{p}_T\widetilde{A}^{(0)}_T\overline{p}_T=0$.

If $s\in \widetilde{E}$ supported in $\mathcal{U}_{\varepsilon_0}$, 
\begin{align}\label{apB13}
C|Z||s|\leq |\widehat{c}(\nabla \widetilde{f})s|\leq C'|Z||s|.
\end{align}
By \cite[(9.87)]{BLe} and (\ref{apB13}), if $\varepsilon\in (0,\varepsilon_0/4]$ is small enough, if $s\in \widetilde{E}^{\bot}_{1,T}$ is supported in $\mathcal{U}_{2\varepsilon}$, then 
\begin{align}\label{apB14}
\left|\widetilde{A}^{(0)}_T s\right|^2_0\geq C\|s\|^2_{\widetilde{E}^1}+C'T^2\left|\widehat{c}\left(\nabla \widetilde{f}\right)s\right|^2_0+C''T |s|^2_0-C'''|s|^2_0.
\end{align}
As \cite[(9.93)]{BLe}, using Lemma \ref{r2lA1}, if $s\in\widetilde{E}$ vanishes on $\mathcal{U}_{\varepsilon}$, then
\begin{align}\label{apB15}
\left|\widetilde{A}^{(0)}_T s\right|^2_0\geq C\|s\|^2_{\widetilde{E}^1}+C'T^2 |s|^2_0-C''|s|^2_0. 
\end{align}
Using (\ref{apB14}), (\ref{apB15}), Lemma \ref{r2lA1} and proceeding as in \cite[p. 115, 116]{BLe}, we find that if $s\in \widetilde{E}^{1,\bot}_{T}$,
\begin{align}\label{apB16}
\left|\widetilde{A}^{(0)}_T s\right|^2_0\geq C\|s\|^2_{\widetilde{E}^1}+C'T^2 \left|\widehat{c}\left(\nabla \widetilde{f}\right)s\right|^2_0+C'' T|s|^2_0-C'''|s|^2_0. 
\end{align}
By Appendix $A$,
\begin{align}\nonumber
\left|\overline{p}^{\bot}_T\widetilde{A}^{(0)}_T \overline{p}_T s\right|_0\leq C\left({{\left\|\overline{p}_T s\right\|_{\widetilde{E}^1}}\over{\sqrt{T}}}+|\overline{p}_T s|_0\right),
\end{align}
\begin{align}\label{apB17}
\left|\overline{p}_T\widetilde{A}^{(0)}_T \overline{p}^{\bot}_T s\right|_0\leq C\left({{\left\|\overline{p}^{\bot}_T s\right\|_{\widetilde{E}^1}}\over{\sqrt{T}}}+|\overline{p}^{\bot}_T s|_0\right).
\end{align}
Using (\ref{apB16}) and the second inequality in (\ref{apB17}), for $T\geq 1$ large enough,
\begin{align}
\left|\overline{p}^{\bot}_T\widetilde{A}^{(0)}_T\overline{p}^{\bot}_T s\right|^2_0\geq C\left\|\overline{p}^{\bot}_Ts\right\|^2_{\widetilde{E}^1}+C'T^2 \left|\widehat{c}\left(\nabla\widetilde{f}\right)\overline{p}^{\bot}_T s\right|^2_0+C'' T\left|\overline{p}^{\bot}_T s\right|^2_0-C'''\left|\overline{p}^{\bot}_T s\right|^2_0.
\end{align}
By definition of $\overline{p}_T$, we have
\begin{align}\label{apB19}
\left\|\overline{p}_T s\right\|_{\widetilde{E}^1}\leq C'\sqrt{T}\left|\overline{p}_T s\right|_0. 
\end{align}
Using (\ref{apB12}), (\ref{apB17})-(\ref{apB19}), for $T\geq 1$ large enough, we obtain the first inequality in (\ref{A17}).

Now we prove the second inequality in (\ref{A17}). Clearly
\begin{align}
\left|\widetilde{A}^{(0)}_T s\right|_0\leq \left|\widetilde{A}^{(0)}_T \overline{p}_T s\right|_0+\left|\widetilde{A}^{(0)}_T\overline{p}^{\bot}_T s\right|_0.
\end{align}
By (\ref{apB8}),
\begin{align}\label{apB21}
\left|\widetilde{A}^{(0)}_T \overline{p}^{\bot}_T s\right|_0\leq C\left(\left\|s\right\|_{\widetilde{E}^1}+T\left|\widehat{c}\left(\nabla\widetilde{f}\right)\overline{p}^{\bot}_Ts\right|_0\right).
\end{align}
From (\ref{apB17}), (\ref{apB19})-(\ref{apB21}), we get the second inequality in (\ref{A17}).

Now we prove the third inequality in (\ref{A17}). Put
\begin{align}
\widetilde{H}=\left[D^{\widetilde{X}},\widetilde{A}^{(>0)}\right]+\widetilde{A}^{(>0),2}.
\end{align}
Then $\widetilde{H}$ is a fiberwise first order differential operator, and moreover
\begin{align}
\widetilde{R}_T=\widetilde{H}+T\left[\widehat{c}\left(\nabla\widetilde{f}\right),\widetilde{A}^{(>0)}\right].
\end{align}
Clearly, by Lemma \ref{r2lA1}, if $s,s'\in \Lambda(T^*S)\widehat{\otimes}\widetilde{E}$,
\begin{align}
\left|\left\langle \widetilde{H}s,s'\right\rangle_0\right|\leq C\left(|s|_{T,1}|s'|_0+|s|_0|s'|_{T,1}\right)+\left|\left\langle\widetilde{H}\overline{p}_Ts,\overline{p}_T s'\right\rangle_0\right|.
\end{align}
Moreover, if $U$ is a smooth section of $T\widetilde{X}$, as \cite[(9.66)]{B}, we have
\begin{align}
\left||Z|\nabla_U^{\Lambda(T^*\widetilde{X})\otimes \widetilde{F}}\overline{p}_T s\right|_0\leq C\left|\overline{p}_T s\right|_0.
\end{align}
Then by \cite[(9.65)]{B}, we have
\begin{align}
\left|\left\langle\widetilde{H}\overline{p}_Ts,\overline{p}_T s'\right\rangle_0\right|\leq C\left(|s|_{T,1}|s'|_0+|s|_0|s'|_{T,1}\right).
\end{align}
Also as \cite[(9.68)-(9.72)]{B}, we have
\begin{align}
T\left|\left\langle \left[\widehat{c}\left(\nabla\widetilde{f}\right),\widetilde{A}^{(>0)}\right]s,s'\right\rangle\right|\leq \left(|s|_{T,1}|s'|_0+|s|_0|s'|_{T,1}\right).
\end{align}
Then from the above we get the third inequality in (\ref{A17}).
\end{proof}

Now we fix $\varepsilon>0$ as in Theorem \ref{thA5}.

By Theorem \ref{r2th5.3}, there is $T_0\geq 1$ such that for $T\geq T_0$, there is a gap of ${\rm Sp}(\widetilde{A}^{(0),2}_T)$.

Take $c_2\in (0,1]$, let $D=\delta\cup \Delta$ be the contour in $\bf C$ defined by 
\begin{align}\nonumber
\delta=\left\{\lambda\in {\bf C}| |\lambda|={{c_2}\over 2}\right\}
\end{align}
and
\begin{align}\nonumber
\Delta=\left\{z=x+iy| x={3\over 4}c_2, -1\leq y\leq 1\right\}\cup \left\{z=x+iy| y=\pm 1, x\geq {3\over 4}c_2\right\}.
\end{align}

If $A\in \mathcal{L}(\widetilde{E}^0,\widetilde{E}^0)$ (resp. $A\in \mathcal{L}(\widetilde{E}^{-1},\widetilde{E}^1)$), let $\|A\|^{0,0}$ (resp. $\|A\|^{-1,1}_{T}$) be the norm of $A$ with respect to the norm $|\ |_0$ (resp. the norms $|\ |_{T,-1}$, $|\ |_{T,1}$). The following $L^2$ analogue of \cite[Theorem 9.15]{B} holds.

\begin{theorem}\label{aptB5}
There exist $T_0\geq 1$, $C>0$, $p\in {\bf N}$, such that for $T>T_0$, $\lambda\in\delta\cup \Delta$, the resolvent $(\lambda-\widetilde{A}^2_T)^{-1}$ is such that
\begin{align}\label{A18}
\left\|\left(\lambda-\widetilde{A}^2_T\right)^{-1}\right\|^{-1,1}_{T}\leq C(1+|\lambda|)^p. 
\end{align}
\end{theorem}

\begin{proof}
Recall that $\widetilde{A}^{(0)}_{T}=D^{\widetilde{X}}+T\widehat{c}(\nabla \widetilde{f})$. For $\delta>0$, $A>0$ set
\begin{align}\label{apB29}
U=\{\lambda\in {\bf C}, {\rm Re}(\lambda)\leq \delta {\rm Im}^2(\lambda)-A\}. 
\end{align}
Using the first two inequalities in (\ref{A17}) and Lemma \ref{r2lA1}, proceeding as in \cite[Theorems 11.26 and 11.27]{BLe}, we find that if $\delta$ is small enough, and $A$ is large enough, for $T\geq 1$, $\lambda\in U$,
\begin{align}\nonumber
\left\|\left(\lambda-\widetilde{A}_{T}^{(0),2}\right)^{-1}\right\|^{0,0}\leq C,
\end{align}
\begin{align}\label{A20}
\left\|\left(\lambda-\widetilde{A}_{T}^{(0),2}\right)^{-1}\right\|^{-1,1}_{T}\leq C(1+|\lambda|)^2.
\end{align}

Take $\lambda\in\Delta$, for $T\geq T_0$, $(\lambda-\widetilde{A}_{T}^{(0),2})^{-1}$ exists and moreover
\begin{align}
\left\|\left(\lambda-\widetilde{A}_{T}^{(0),2}\right)^{-1}\right\|^{0,0}\leq C.
\end{align}
If $\lambda_0\in U$, $\lambda\in\delta\cup\Delta$, $T\geq T_0$, then
\begin{align}\label{A22}
\left(\lambda-\widetilde{A}_{T}^{(0),2}\right)^{-1}=\left(\lambda_0-\widetilde{A}_{T}^{(0),2}\right)^{-1}+\left(\lambda-\widetilde{A}_{T}^{(0),2}\right)^{-1}(\lambda_0-\lambda)\left(\lambda_0-\widetilde{A}_{T}^{(0),2}\right)^{-1}.
\end{align}
From (\ref{A20})-(\ref{A22}), we get
\begin{align}\label{A23}
\left\|\left(\lambda-\widetilde{A}_{T}^{(0),2}\right)^{-1}\right\|^{-1,0}_{T}\leq C(1+|\lambda|).
\end{align}
Also 
\begin{align}\label{A24}
\left(\lambda-\widetilde{A}_{T}^{(0),2}\right)^{-1}=\left(\lambda_0-\widetilde{A}_{T}^{(0),2}\right)^{-1}+\left(\lambda_0-\widetilde{A}_{T}^{(0),2}\right)^{-1}(\lambda_0-\lambda)\left(\lambda-\widetilde{A}_{T}^{(0),2}\right)^{-1}.
\end{align}
By (\ref{A20}), (\ref{A23}), (\ref{A24}), we obtain 
\begin{align}\label{A25}
\left\|\left(\lambda-\widetilde{A}_{T}^{(0),2}\right)^{-1}\right\|^{-1,1}_{T}\leq C(1+|\lambda|)^2.
\end{align}
Moreover, if $\lambda\in\delta\cup \Delta$, then
\begin{align}\label{A26}
\left(\lambda-\widetilde{A}^2_T\right)^{-1}=\left(\lambda-\widetilde{A}_{T}^{(0),2}\right)^{-1}+\left(\lambda-\widetilde{A}_{T}^{(0),2}\right)^{-1}\widetilde{R}_T \left(\lambda-\widetilde{A}_{T}^{(0),2}\right)^{-1}+\cdots
\end{align}
and the expansion terminates after a finite numbers of terms. By Theorem \ref{thA5},
\begin{align}\label{A27}
\left\|\widetilde{R}_T\right\|_{T}^{-1,1}\leq C.
\end{align}
Using (\ref{A25})-(\ref{A27}), we get (\ref{A18}). The proof is completed. 
\end{proof}

Since $\bf {B}$ is compact, there exist a finite family of smooth functions $f_1,\cdots,f_q$ on $M$ with values in $[0,1]$, such that
\begin{align}
{\bf {B}}=\cap _{j=1}^{q}\left\{x\in M, f_j(x)=0\right\},
\end{align}
and that on $\bf {B}$, $df_1,\cdots,df_q$ span $T^* X$. We lift $f_1,\cdots,f_q$ to $\bf {B}$ and denote them by $\widetilde{f}_1,\cdots,\widetilde{f}_q$.

Similarly, there exists a finite family of smooth sections $U_1,\cdots,U_r$ of $TX$ such that for any $x\in V$, $U_1(x),\cdots, U_r(x)$ spans $(TX)_x$. Let $\widetilde{U}_1,\cdots, \widetilde{U}_r$ be the liftings of 
${U}_1,\cdots,{U}_r$. 
\begin{definition}
For $T\geq 1$, let $\widetilde{L}_T$ be the family of operators acting on $\widetilde{E}$
\begin{align}
\widetilde{L}_T
=
\left\{\nabla^{\Lambda(T^*\widetilde{X}){\otimes}\widetilde{F}}_{(1-\rho(Z/2))\widetilde{U}_i},{1\over{\sqrt{T}}}\bar{p}^{\bot}_T {^{0}}\widetilde{\nabla}^{(\Lambda(T^{*}\widetilde{X}){\otimes}\widetilde{F})|_{\widetilde{\bf{B}}}}_{\rho(Z/2)\widetilde{U}_i}\bar{p}^{\bot}_{T}, \sqrt{T}\bar{p}^{\bot}_{T}\widetilde{f}_j \bar{p}^{\bot}_{T}\right\}.
\end{align}
\end{definition}
For $k\in {\bf N}$, let $Q^k_T$ be the family of operators $Q$ acting on $\widetilde{E}$ which can be written in the form
\begin{align}
Q=Q_1\cdots Q_k,\ \ Q\in \widetilde{L}_T.
\end{align}
If $k\in{\bf N}$, we equip the Sobolev fibers $\widetilde{E}^k$ with the Hilbert norm $\|\ \|_{T,k}$ such that if $s\in\widetilde{E}$,
\begin{align}
\|s\|^2_{T,k}=\sum_{l=0}^k \sum _{Q\in\mathcal{L}^l_T}|Qs|^2_{T,0}.
\end{align}

By the same proof of \cite[Theorem 9.17]{B}, we have
\begin{theorem}\label{aptB7}
Take $k\in{\bf N}$. There exists $C_k>0$ such that for $T\geq 1$, $Q_1,\cdots, Q_k\in\widetilde{L}_T$, $s,s'\in\Lambda(T^*S)\widehat{\otimes}\widetilde{E}$,
\begin{align}
\left|\left\langle \left[Q_1,\left[Q_2,\cdots\left[Q_k,\widetilde{A}^2_T\right]\right]\right]s,s'\right\rangle_0\right|\leq C_k |s|_{T,1} |s'|_{T,1}. 
\end{align}
\end{theorem}

If $A\in \mathcal{L}(\widetilde{E}^m,\widetilde{E}^{m'})$, we denote by $|||A|||^{m,m'}_{T}$ the norm of $A$ with respect to the norms $\|\ \|_{T,m}$, $\|\ \|_{T,m'}$. Then by the same proof of \cite[Theorem 9.18]{B}, we have

\begin{theorem}\label{thA9}
For any $m\in{\bf N}$, there exist $p_m\in {\bf N}$, $C_m>0$ such that for $T\geq T_0$, $\lambda\in\Delta$,
\begin{align}
\left|\left|\left|\left(\lambda-\widetilde{A}^2_T\right)^{-1}\right|\right|\right|^{m,m+1}_{T}\leq C_m (1+|\lambda|)^{p_m}.
\end{align}
\end{theorem}
\begin{proof}
Clearly for $T\geq 1$,
\begin{align}\label{apB44}
\|s\|_{T,1}\leq C|s|_{T,1}.
\end{align}
When $m=0$, our Theorem follows from Theorem \ref{aptB5} and from (\ref{apB44}).

Using Theorems \ref{aptB5} and \ref{aptB7}, the proof of our Theorem proceed as the proof of \cite[Theorem 11.30]{BLe}.
\end{proof}

If $a\notin \Delta$, put
\begin{align}\label{A34}
F_u(a)={1\over{2\pi i}}\int_{\Delta}{{\exp(-u^2 \lambda)}\over{\lambda-a}}d\lambda. 
\end{align}
Then 
\begin{align}\nonumber
F_u(a)=\left\{\begin{matrix}\exp(-u^2 a)& {\rm if}\ a\ {\rm lies\ inside\ the\ contour}\ \Delta,\\
0& {\rm if}\ a\ {\rm lies\ outside}\ \Delta.\end{matrix}\right.
\end{align}
Put 
\begin{align}
F_u \left(\widetilde{A}^2_T\right)={1\over{2\pi i}}\int_{\Delta}{{\exp(-u^2 \lambda)}\over{\lambda-\widetilde{A}^2_T}}d\lambda.
\end{align}

\begin{definition}
Let $F_u(\widetilde{A}^2_T)(x,x')$ $x,x'\in\widetilde{X}$ be the smooth kernel associated to the operator $F_u (\widetilde{A}^2_T)$ with respect to ${{dv_{\widetilde{X}}(x')}\over{(2\pi)^{{\rm dim}\widetilde{X}}}}$.
\end{definition}

\begin{theorem}\label{At11}
For any $\alpha>0$, $m\in {\bf N}$, there exist $C>0$, $C'>0$ such that if $x\in \widetilde{V}$, $d^{\widetilde{X}}(x,\widetilde{B})\geq \alpha$, for $u\geq u_0$, $T\geq T_0$,
\begin{align}\label{A36}
\left|F_u \left(\widetilde{A}^2_T\right)(x,x')\right|\leq {{C\exp(-C'u^2)}\over{T^m}}.
\end{align}
For any $m\in {\bf N}$, there exist $C>0$, $C'>0$ such that for $y\in\widetilde{\bf{B}}$, $u\geq u_0$, $T\geq T_0$,
\begin{align}\label{A37}
\sup_{|Z|\leq {\varepsilon\over 2}\sqrt{T}}(1+|Z|)^m {1\over{T^{{\rm dim}{\widetilde{X}}}}}\left|F_u\left(\widetilde{A}^2_T\right)\left(\left(y,{Z\over\sqrt{T}}\right),\left(y,{Z\over{\sqrt{T}}}\right)\right)\right|\leq C\exp (-C'u^2).
\end{align}
For any $m\in{\bf N}$, there exist $C>0$, $C'>0$ such that for $y\in\widetilde{\bf B}$, $u\geq u_0$, $T\geq T_0$,
\begin{align}\label{A38}
\sup_{ \scriptstyle |\alpha|\leq m', |\alpha'|\leq m'\atop {\atop{ \scriptstyle |Z|\leq {\varepsilon\over 2}\sqrt{T}\atop\scriptstyle |Z'|\leq {\varepsilon\over 2}\sqrt{T}}}}\left|{{\partial^{|\alpha|+|\alpha'|}}\over{\partial x^\alpha \partial x'^{\alpha'}}}{1\over{T^{{\rm dim}{\widetilde{X}}}}}F_u\left(\widetilde{A}^2_T\right)\left(\left(y,{Z\over\sqrt{T}}\right),\left(y',{Z'\over{\sqrt{T}}}\right)\right)\right|\leq C\exp(-C'u^2).
\end{align}
\end{theorem}
\begin{proof}
Clearly for $p\in {\bf N}$,
\begin{align}\label{A39}
{1\over{2\pi i}}\int_{\Delta} {{\exp(-u^2\lambda)}\over{\lambda-\widetilde{A}^2_T}}d\lambda=(-1)^{2p-1}{{(2p-1)!}\over{2\pi i (u^2)^{2p-1}}}\int_{\Delta}{{\exp(-u^2\lambda)}\over{(\lambda-\widetilde{A}^2_T)^{2p}}}d\lambda. 
\end{align}
By Theorem \ref{thA9}, we know that there exists $C>0$, $q\in {\bf N}$ such that if $\lambda\in\Delta$, $Q\in\mathcal{L}_T^l$, $l\leq p$,
\begin{align}\label{A40}
\left\|Q\left(\lambda-\widetilde{A}^2_T\right)^{-p}\right\|^{0,0}_{T}\leq C(1+|\lambda|)^q.
\end{align} 
By introducing the obvious adjoint operator with respect to $\langle\cdot\ ,\ \cdot\rangle_0$, we also find that if $\lambda\in\Delta$, $Q'\in\mathcal{L}^l_T$, $l\leq p$,
\begin{align}\label{A41}
\left\|\left(\lambda-\widetilde{A}^2_T\right)^{-p}Q'\right\|^{0,0}_{T}\leq C(1+|\lambda|)^q.
\end{align} 
From (\ref{A40}), (\ref{A41}), we see that if $\lambda\in\Delta$, $Q\in\mathcal{L}^{l}_T$, $Q'\in\mathcal{L}^{l'}_T$, $l,l'\leq p$,
\begin{align}\label{A42}
\left\|Q\left(\lambda-\widetilde{A}^2_T\right)^{-2p}Q'\right\|^{0,0}_{T}\leq C(1+|\lambda|)^{2q}.
\end{align} 
From (\ref{A39}), (\ref{A42}), we find that if $Q\in\mathcal{L}^l_T$, $Q'\in\mathcal{L}^{l'}_{T}$, there exist $C>0$, $C'>0$ such that 
\begin{align}\label{A43}
\left\|QF_u \left(\widetilde{A}^2_T\right)Q'\right\|^{0,0}_{T}\leq C\exp\left(-C'u^2\right).
\end{align}
By (\ref{A43}) and by Sobolev inequalities (cf. \cite{H}), we get (\ref{A36}). Using (\ref{A42}) and proceeding as in \cite[proof of Theorem 13.32]{BLe}, we get (\ref{A37}), (\ref{A38}). The proof is completed. 
\end{proof}

By the same proof of \cite[Proposition 9.21]{B}, we have
\begin{proposition}\label{thA12}
There exist $C>0$, $p\in {\bf N}$ such that for $T\geq T_0$, $\lambda\in\delta \cup \Delta$, then
\begin{align}
\left\|\bar{p}^{\bot}_{T}\left(\lambda-\widetilde{A}^2_T\right)^{-1}\right\|^{0,0}\leq {C\over{\sqrt{T}}}(1+|\lambda|)^p. 
\end{align}
\end{proposition}
\begin{proof}
This follows from Theorem \ref{aptB5}.
\end{proof}

If $A$ is an bounded operator acting on $\widetilde{E}^0$, we write $A$ in matrix form with respect to the splitting $\widetilde{E}^0=\widetilde{E}^0_T\oplus \widetilde{E}^{0,\bot}_{T}$,
\begin{align}
A=\left(\begin{matrix}A_1&A_2\\ A_3& A_4\end{matrix}\right).
\end{align}

By the same proof of \cite[Proposition 9.22]{B}, we have
\begin{proposition}
There exist $C>0$, $p\in{\bf N}$, $T_0\geq 1$ such that if $T\geq T_0$, $\lambda\in\delta\cup \Delta$, the resolvent $(\lambda-\widetilde{A}^2_{T,4})^{-1}$ exists and moreover
\begin{align}\label{apB57}
\left\|\left(\lambda-\widetilde{A}^2_{T,4}\right)^{-1}\right\|^{-1,1}_{T}\leq C(1+|\lambda|)^p.
\end{align}
\end{proposition}
\begin{proof}
By Theorem \ref{thA5}, it is clearly that $\widetilde{A}^2_{T,4}$ verifies inequalities similar to (\ref{A17}). Therefore by using the notation in (\ref{apB29}), for $\delta>0$ small enough, and $A>0$ large enough, if $\lambda\in U$, then
\begin{align}\nonumber
\left\|\left(\lambda-\widetilde{A}^2_{T,4}\right)^{-1}\right\|^{0,0}\leq C,
\end{align}
\begin{align}\label{apB58}
\left\|\left(\lambda-\widetilde{A}^2_{T,4}\right)^{-1}\right\|^{-1,1}_{T}\leq C(1+|\lambda|)^2.
\end{align}
By Theorem \ref{thA5}, for $T\geq 1$ large enough, if $s\in \widetilde{E}^{1,\bot}_{T}$,
\begin{align}\label{apB59}
\left\langle \widetilde{A}^{(0),2}_Ts,s\right\rangle_0\geq CT \left|\overline{p}^{\bot}_T\right|^2_0.
\end{align}
By (\ref{apB59}), we find that there is $C>0$, $T_0\geq 1$ such that for $T\geq T_0$, $\lambda\in \Delta\cup \delta$,
\begin{align}\label{apB60}
\left\|\left(\lambda-\widetilde{A}^{(0),2}_{T,4}\right)^{-1}\right\|^{0,0}\leq C.
\end{align}
Using (\ref{apB58}), (\ref{apB60}), and proceeding as in (\ref{A20})-(\ref{A25}), we get for $T\geq T_0$, $\lambda\in \Delta\cup\delta$,
\begin{align}
\left\|\left(\lambda-\widetilde{A}^{(0),2}_{T,4}\right)^{-1}\right\|_T^{-1,1}\leq C(1+|\lambda|)^2.
\end{align}
Thus if $\lambda\in\Delta\cup\delta$,
\begin{align}\label{apB62}
\left(\lambda-\widetilde{A}^2_{T,4}\right)^{-1}=\left(\lambda-\widetilde{A}^{(0),2}_{T,4}\right)^{-1}+\left(\lambda-\widetilde{A}^{(0),2}_{T,4}\right)^{-1}\widetilde{R}_{T,4}\left(\lambda-\widetilde{A}^{(0),2}_{T,4}\right)^{-1}+\cdots.
\end{align}
By Theorem \ref{thA5}, we get
\begin{align}\label{apB63}
\left\|\widetilde{R}_{T,4}\right\|^{1,-1}_{T}\leq C.
\end{align}
By (\ref{apB62}) and (\ref{apB63}), we get (\ref{apB57}).
\end{proof}

By Proposition \ref{thA12}, proceeding as \cite[Theorem 9.23]{B}, we have
\begin{theorem}
There exist $C>0$, $C'>0$ such that for $u\geq 1$, $T\geq T_0$,
\begin{align}
\left\|\bar{p}^{\bot}_{T}F_u \left(\widetilde{A}^2_T\right)\bar{p}^{\bot}_T\right\|^{0,0}\leq {C\over{\sqrt{T}}},\ 
\left\|\bar{p}^{\bot}_{T}F_u \left(\widetilde{A}^2_T\right)\bar{p}_T\right\|^{0,0}\leq {C\over{\sqrt{T}}},\
\left\|\bar{p}_{T}F_u \left(\widetilde{A}^2_T\right)\bar{p}^{\bot}_T\right\|^{0,0}\leq {C\over{\sqrt{T}}}.\end{align}
\end{theorem}
\begin{proof}
In view of Proposition \ref{thA12}, we can proceed exactly as in \cite[p. 264-267]{BLe}.
\end{proof}

By the same proof of \cite[Theorem 9.27]{B}, we have
\begin{theorem}\label{thA15}
There exist $c>0$, $C>0$ such that for $u\geq u_0$, $T\geq 1$,
\begin{align}
\left\|F_{u}\left(\widetilde{A}^2_T\right)\right\|^{0,0}\leq {{c\exp(-Cu^2)}\over{T^{1/4}}}.
\end{align}
\end{theorem}

\begin{proof}
The proof is the same as the proof of \cite[Theorem 13.42]{BLe}.
\end{proof}

Using above results, we now prove the $L^2$-case of \cite[(9.149)]{B}. Let $Y'\subset \widetilde{B}$ be a fundamental domain and $X'\subset \widetilde{X}$ be a fundamental domain and such that $Y'\subset X'$. Then
\begin{align}\label{A51}
{\rm Tr}_{\Gamma,s}\left[N{1\over{2\pi i}}\int_{\Delta}{{\exp(-u^2 \lambda)}\over{\lambda-\widetilde{A}^2_T}}\right]=
\int_{X'}{\rm Tr}_{s}\left[N F_u \left(\widetilde{A}^2_T\right)(x,x)\right]{{dv_{\widetilde{X}}(x)}\over{(2\pi)^{{\rm dim}\widetilde{X}}}}.
\end{align}
By Theorem \ref{At11}, for any $m\in {\bf N}$
\begin{align}
\left|\int_{X'\cap \{x,d^{\widetilde{X}}(x,Y')\geq \varepsilon/4\}}{\rm Tr}_{s}\left[N F_u \left(\widetilde{A}^2_T\right)(x,x)\right]{{dv_{\widetilde{X}}(x)}\over{(2\pi)^{{\rm dim}\widetilde{X}}}}\right|
\leq {C\over{T^m}}\exp\left(-C'u^2\right).
\end{align}
Also
\begin{multline}\label{A53}
\int_{X'\cap \{x,d^{\widetilde{X}}(x,Y')\leq \varepsilon/4\}}{\rm Tr}_{s}\left[N F_u \left(\widetilde{A}^2_T\right)(x,x)\right]{{dv_{\widetilde{X}}(x)}\over{(2\pi)^{{\rm dim}\widetilde{X}}}}=\\
\sum_{y\in Y'}\int_{Z\in TX',|Z|\leq {\varepsilon\over 4}\sqrt{T}}{\rm Tr}_s\left[N {{F_u\left(\widetilde{A}^2_T\right)}\over{T^{{\rm dim}\widetilde{X}}}}\left(\left(y,{Z\over\sqrt{T}}\right),\left(y,{Z\over\sqrt{T}}\right)\right)\right]
{{dv_{T\widetilde{X}}(Z)}\over{(2\pi)^{{\rm dim}\widetilde{X}}}}.
\end{multline}

By Theorems \ref{At11} and \ref{thA15} and proceeding as in \cite[Section 13 q]{BLe}, we find that there exist $c>0$, $C>0$, $\delta\in (0,1/2]$ such that if $y\in\widetilde{\bf{B}}$, $Z\in T\widetilde{X}$, $|Z|\leq {{\varepsilon \sqrt{T}}\over{4}}$,
\begin{align}\label{A54}
\left|\left({1\over {2\pi}}\right)^{{\rm dim}\widetilde{X}}{{F_u\left(\widetilde{A}^2_T\right)}\over{T^{{\rm dim}{\widetilde{X}}}}}\left(\left(y,{Z\over\sqrt{T}}\right),\left(y,{Z\over\sqrt{T}}\right)\right)
\right|\leq {{c\exp(-Cu^2)}\over{T^\delta}}.
\end{align}

By (\ref{A37}), (\ref{A54}), we find that for any $p\in{\bf N}$, there is $c>0$, $C>0$ such that if $y\in\widetilde{{\bf B}}$, $Z\in T_y\widetilde{X}$, $|Z|\leq {{\varepsilon\sqrt{T}}\over{4}}$,
\begin{align}\label{A55}
\left|\left({1\over {2\pi}}\right)^{{\rm dim}\widetilde{X}}{{F_u\left(\widetilde{A}^2_T\right)}\over{T^{{\rm dim}{\widetilde{X}}}}}\left(\left(y,{Z\over\sqrt{T}}\right),\left(y,{Z\over\sqrt{T}}\right)\right)
\right|\leq {{c\exp(-Cu^2)}\over{(1+|Z|)^p T^{\delta/2}}}.
\end{align}
From (\ref{A55}), we deduce that there exist $c>0$, $C>0$, $\delta\in (0,1/4]$ such that 

\begin{multline}\label{A57}
\left|\sum_{y\in Y'}\int_{Z\in N_{Y'/X'},|Z|\leq {\varepsilon\over 4}\sqrt{T}}{\rm Tr}_s\left[N{{F_u\left(\widetilde{A}^2_T\right)}}\left(\left(y,{Z\over\sqrt{T}}\right),\left(y,{Z\over\sqrt{T}}\right)\right)\right]{{dv_{T\widetilde{X}}(Z)}\over{(2\pi)^{{\rm dim}\widetilde{X}}}}\right|\\
\leq {{c\exp(-Cu^2)}\over{T^\delta}}.
\end{multline}
So by (\ref{A51})-(\ref{A53}), (\ref{A57}), we obtain
\begin{align}
\left|{\rm Tr}_{\Gamma,s}\left[NF_u \left(\widetilde{A}^2_T\right)\right]\right|\leq {{c\exp(-Cu^2)}\over{T^\delta}}.
\end{align}

\vspace{\baselineskip}

\end{document}